\documentclass[10pt]{article} % For LaTeX2e
\usepackage[accepted]{tmlr}
\usepackage{amsmath,amsfonts,bm,amsthm}

\def\1{\bm{1}}

\DeclareMathAlphabet{\mathsfit}{\encodingdefault}{\sfdefault}{m}{sl}
\SetMathAlphabet{\mathsfit}{bold}{\encodingdefault}{\sfdefault}{bx}{n}

\usepackage[colorlinks=true, urlcolor=cyan, citecolor=blue]{hyperref}
\usepackage{url}
\newcommand{\mc}{\mathbb}
\newtheorem{assume}{Assumption}
\newtheorem{lemma}{Lemma}

\newtheorem{definition}{Definition}

\newtheorem{theorem}{Theorem}
\newtheorem{corollary}{Corollary}
\newtheorem{proposition}{Proposition}
\usepackage{algorithm}
\usepackage{subcaption}
\usepackage{algpseudocode}
\usepackage{thmtools, thm-restate}
\usepackage{graphicx}
\usepackage{amsthm}
\usepackage{tikz}
\definecolor{winered}{rgb}{0.8,0,0}
\newcommand*\circled[1]{\tikz[baseline=(char.base)]{
            \node[shape=circle,draw,inner sep=2pt] (char) {#1};}}

\allowdisplaybreaks
\title{Model-Free Learning with Heterogeneous Dynamical Systems: A Federated LQR Approach}

\author{\name Han Wang \email hw2786@columbia.edu\\
      \addr Department of Electrical Engineering\\
      Columbia University
      \AND
      \name Leonardo F. Toso \email leonardo.toso@columbia.edu \\
      \addr Department of Electrical Engineering\\
      Columbia University
      \AND
      \name Aritra Mitra \email amitra2@ncsu.edu\\
      \addr  Department of Electrical and Computer Engineering \\
      North Carolina State University
      \AND
      \name James Anderson \email james.anderson@columbia.edu\\
      \addr Department of Electrical Engineering\\
      Columbia University}

\def\month{02}  % Insert correct month for camera-ready version
\def\year{2026} % Insert correct year for camera-ready version
\def\openreview{\url{https://openreview.net/forum?id=WSRQeCUc3g}} % Insert correct link to OpenReview for camera-ready version

\begin{document}

\maketitle

\begin{abstract}
We study a model-free federated linear quadratic regulator (LQR) problem where $M$ agents with \emph{unknown}, \emph{distinct yet similar} dynamics collaboratively learn an optimal policy to minimize an average quadratic cost while keeping their data private. To exploit the similarity of the agents' dynamics, we propose to use federated learning (FL) to allow the agents to periodically communicate with a central server to train policies by leveraging a larger dataset from all the agents. With this setup, we seek to understand the following questions: (i) \textit{Is the learned common policy stabilizing for all agents?} (ii) \textit{How close is the learned common policy to each agent's own optimal policy?}  (iii) \textit{Can each agent learn its own optimal policy faster by leveraging data from all agents}? To answer these questions, we propose the federated and model-free algorithm  \texttt{FedLQR}. Our analysis overcomes numerous technical challenges, such as heterogeneity in the agents’ dynamics, multiple local updates, and stability concerns. We show that \texttt{FedLQR} produces a common policy that, at each iteration, is stabilizing for all agents. Moreover, we prove that when learning each agent's optimal policy, \texttt{FedLQR}  achieves a sample complexity reduction proportional to the number of agents $M$ in a low-heterogeneity regime, compared to the single-agent setting. 
\end{abstract}

\section{Introduction}

There has been significant progress in the application of model-free reinforcement learning (RL) methods to fields such as video games~\citep{mnih2015human} and robotic manipulation~\citep{rajeswaran2017learning, levine2016end,tobin2017domain}. In particular, RL has recently been used in the fine-tuning of pretrained foundation models \citep{ziegler2019fine}. Although RL has shown impressive results in simulation, it often suffers from poor sample complexity, thereby limiting its effectiveness in real-world applications~\citep{dulac2019challenges}. To resolve the sample complexity issue and accelerate the learning process, federated learning (FL) has emerged as a popular paradigm ~\citep{konevcny2016federated,mcmahan2017communication}, where multiple similar agents collaboratively learn a common model without sharing their raw data. The incentive for collaboration arises from the fact that these agents are ``similar'' in some sense and hence end up learning a ``superior'' model than if they were to learn alone. In the RL setting, Federated Reinforcement Learning (FRL) aims to learn a common value function~\citep{wang2023federated} or produce a better policy from multiple RL agents interacting with similar environments. In the survey paper by \cite{qi2021federated}, FRL has empirically shown great success in reducing the sample complexity for autonomous driving~\citep{liang2022federated}, IoT devices~\citep{lim2020federated}, and resource management in networking~\citep{yu2020deep}. 

Lately, there has been a lot of interest in applying RL techniques to classical control problems such as the linear quadratic regulator (LQR) problem~\citep{anderson2007optimal}. In the standard control setting, the dynamical model of the system is known and one seeks to obtain a controller that stabilizes the closed-loop system and provides optimal performance. RL approaches such as policy gradient~\citep{williams1992simple,sutton1999policy} (which we pursue here) differ in that they are  ``model-free'', i.e., a control policy is obtained despite not having access to the model of the dynamics. Despite the lack of convexity in even simple problems, policy gradient (PG) methods have been shown to be globally convergent for certain structured settings such as the LQR problem~\citep{fazel2018global}. While this is promising, a major challenge in applying PG methods is that in general, one does not have access to \textit{exact deterministic} policy gradients. Instead, one relies on estimating such gradients via sampling based approaches. This typically leads to noisy gradients that can suffer from high variance. As such, reducing the variance in PG estimates to achieve ``good performance" may end up requiring several samples. 

\begin{figure}[t]
\centering
  \includegraphics[width=0.65\linewidth]{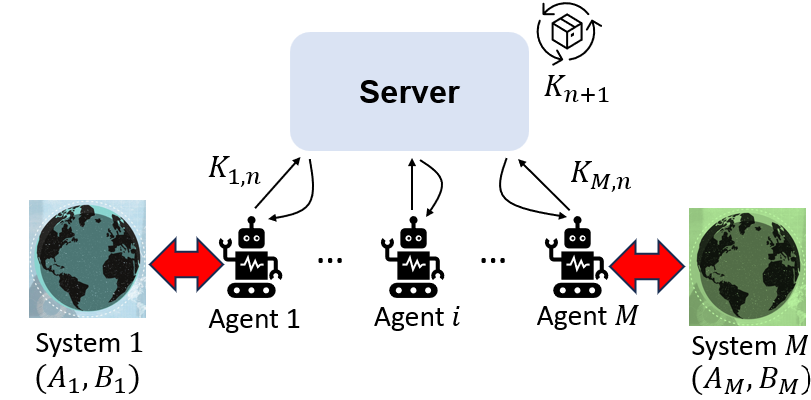}
  \caption{ Illustration of the heterogeneous federated LQR setup involving $M$ agents, where the $i$-th agent interacts with its own LTI system characterized by the pair $(A_i, B_i)$ for $i \in [M]$. The agents communicate via a central server in communication rounds $n=1, 2, \ldots, N$, and upload their locally computed control gain matrices $\{K^{(i)}_{n}\}_{i \in [M]}$ in each round $n$. The server averages the agents' gain matrices to compute a global control policy $K_{n+1}$ that is broadcast to all agents to initiate the $(n+1)$-th round.}
  \label{fig:Setup}
\end{figure} 

\textbf{Motivation.} The main premise of this paper is to draw on ideas from the FL literature to alleviate the high sample-complexity burden of PG methods~\citep{agarwal2019logarithmic, wang2019neural,liu2020improved}, with the focus being on model-free control. As a motivating example, consider a fleet of identical robots produced by the same manufacturer. Each robot can collect data from its own dynamics and learn its own optimal policy using, for instance, PG methods.  Since the fleet of robots shares similar dynamics, and more data can potentially lead to improved policy performance (via more accurate PG estimates), it is natural to ask: \textit{Can a robot accelerate the process of learning its own optimal policy by leveraging the data of the other robots in the fleet?} The answer is not as obvious as one might expect since in reality, it is  unlikely that any two robots will have \textit{exactly} the same underlying dynamics, i.e., \textit{heterogeneity in system dynamics is inevitable.} The presence of such heterogeneity makes the question posed above both interesting and non-trivial. In particular, when the heterogeneity across agents' dynamics is large, leveraging data from other agents might degrade the performance of a single agent. Indeed, large heterogeneity may make it impossible to learn a common stabilizing policy\footnote{See Section \ref{sec:heterogeneity_regime} for more details on the underlying intuition and necessity behind the low heterogeneity regime.}. Moreover, even when such a stabilizing policy exists, it may deviate from each agent's local optimal policy, rendering poor performance and discouraging participation in the FL process. Thus, to understand whether more data\footnote{In accordance with both FL \& FRL frameworks, the agents in our problem do not exchange their private data (e.g., rewards, states, etc.). Instead, each agent only transmits its policy gradient.  } helps or hurts, it is crucial to characterize the effects of heterogeneity in the federated control setting.

With this aim in mind, we study a multi-agent \textit{model-free} LQR problem based on policy gradient methods. Specifically, there are $M$ agents in our setup, each with its own \emph{distinct yet similar} linear time-invariant (LTI) dynamics. Inspired by the typical objective in FL, our goal is to find a common policy which can minimize the average of the LQR costs across agents.\footnote{While one can certainly consider other FL objective functions, we study the average LQR cost formulation here since a comprehensive analysis of this setting is missing in the current literature.} In particular, we address the following questions. 

\textbf{Q1.} \emph{Is this common policy stabilizing for all the systems?  If so, under what conditions?}

\textbf{Q2.} \emph{How far is the learned common policy from each agent's locally optimal policy?} 

\textbf{Q3.} \textit{Can an agent use the common policy as an initial guess to fine-tune and learn its own optimal policy faster (i.e., with fewer overall samples) than if it acted alone?}

\textbf{Challenges:} There are several challenges to answering the above questions. First, even for the single agent setting, the policy gradient-based LQR problem is non-convex, and requires a fairly intricate analysis~\citep{fazel2018global}. Second, a key distinction relative to  standard federated supervised learning stems from the need to guarantee \textit{stability} at every averaging iteration -- this problem is amplified in the heterogeneous multi-agent scenario we consider. 
 It remains an open problem to design an algorithm ensuring that  policies  are simultaneously stabilizing for each distinct system. Third, to reduce the communication cost, FL algorithms rely on the agents performing multiple local update steps between successive communication rounds. When agents have non-identical loss functions, these local steps lead to a ``client-drift" effect where each agent drifts towards its own local minimizer~\citep{charles,charles2}. While several works in FL have investigated this phenomenon~\citep{li2020federated, khaled1, khaled2020tighter, li, karimireddy2020scaffold, fedsplit, FedNova, acar2021, gorbunov, mitra2021linear, proxskip, quantile}, \emph{the effect of ``client-drift" on stability remains completely unexplored}. Unless accounted for, such drift effects can potentially lead to non-stabilizing controllers. 

\textbf{Our Contributions:} In response to the above challenges, we propose a policy gradient method called \texttt{FedLQR} to solve the (model-based \emph{and} model-free) federated LQR problem, and provide a rigorous finite-time analysis of its performance that accounts for the interplay between system heterogeneity, multiple local steps, client-drift effects, and stability. Our specific contributions in this regard are as follows.

   $\bullet$  \textbf{Iterative stability guarantees.}  We show via a careful inductive argument that under suitable requirements on the level of heterogeneity across systems, the learning rate schedule can be designed to ensure that \texttt{FedLQR} provides a stabilizing controller at every iteration for \emph{all} systems.  Theorem~\ref{thm:one_local_model_based} provides a proof in the model-based setting, and Theorem~\ref{thm:main_fedlqr} provides the model-free result. 

   $\bullet$ \textbf{Bounded policy gradient heterogeneity in the LQR problem.}  We prove in Lemma~\ref{lem:gradient_het} that, for each pair of agents $i, j\in [M],$ the policy gradient direction (in the model-based setting) of agent $i$ is close to that of agent $j$, if their dynamics are similar (i.e., Definition~\ref{assumption:bnd_sys_heterogeneity}).  
\textit{This is the first result to observe and characterize this bounded gradient heterogeneity phenomenon in the multi-agent LQR setting.}
   
   $\bullet$ \textbf{Quantifying the gap between  \texttt{FedLQR}'s output and each system's optimal policy.} Building on Lemma~\ref{lem:gradient_het}, we prove that when the agents' dynamics are similar, the common policy returned by \texttt{FedLQR} is close to each agent's optimal policy; see Theorem~\ref{thm:one_local_model_based}. In other words, we can leverage the federated formulation to help each agent find its own optimal policy up to some accuracy depending on the level of heterogeneity. Moreover, we prove that \texttt{FedLQR} finds a controller that is close to the average's cost optimal solution up to a heterogeneity bias (Corollary \ref{coro:global convergence}). Our work is the first to provide a result of this flavor.

     $\bullet$ \textbf{Linear speedup.} As our main contribution, we prove that in the model-free setting, \texttt{FedLQR} converges to a solution that is in a neighborhood of each agent's optimal policy, \textit{using $M$-times fewer samples} relative to when each agent just uses its own data (see Theorem~\ref{thm:main_fedlqr}). The radius of this neighborhood captures the level of heterogeneity across the agents' dynamics.  \textit{The key implication of this result is that in a low-heterogeneity regime, \texttt{FedLQR} (in the model-free setting) reduces the sample-complexity by a factor of $M$ w.r.t. the centralized setting~\citep{fazel2018global, malik2019derivative}, highlighting the benefit of collaboration.}\footnote{Throughout this paper, we use the terms ``centralized" and ``single-agent" interchangeably.} Simply put, \texttt{FedLQR} enables each agent to \textit{quickly} find an \textit{approximate} locally optimal policy; as in standard FL~\citep{collins2022fedavg}, the agent can \textit{fine-tune} this policy based on its own data. 

     In summary, we provide a new theoretical framework that quantitatively \emph{characterizes the interplay between the price of heterogeneity and the benefit of collaboration} for model-free control. 

\section{Related Work}  There has been a line of work that explores various RL algorithms for solving the model-free LQR problem \citep{fazel2018global, malik2019derivative, hambly2021policy,mohammadi2021convergence, gravell2020learning, jin2020analysis,ju2022model}. However, their analysis is limited to the single-agent setting. Most recently, \cite{ren2020federated} solve the model-free LQR tracking problem in a federated manner and achieve a linear convergence speedup with respect to
the number of agents. However, they consider a  simplified setting where all agents follow the \emph{same} dynamics. As such, the stability analysis of~\cite{ren2020federated} follow from arguments for the centralized setting. In sharp contrast, to establish the linear speedup for \texttt{FedLQR}, we need to address the key technical challenges arising from the effect of heterogeneity and local steps on the stability of distinct systems. This requires new analysis tools that we develop. {For related work on multi-agent RL (that do not specifically  look at the control setting) we point the reader to \cite{lin2021multi,zhangsurv}} and the references therein. Below we highlight the relevant work related to our problem setting.

\begin{itemize}
   \item \textbf{Federated Learning (FL):} In this work, we employ the federated learning (FL) paradigm to facilitate collaborative learning among systems without the need to share raw data with other participants or a  server~\citep{konevcny2016federated, mcmahan2017communication, konevcny2016FL, bonawitz2019towards}. Despite FL being a relatively recent approach, it has already gained significant attention and boasts a wealth of literature. 
  
   Federated averaging (\texttt{FedAvg}) stands as the pioneering and most widely adopted algorithm in FL. Originally proposed in \cite{mcmahan2017communication}, \texttt{FedAvg} has demonstrated its effectiveness in homogeneous settings \citep{stich2018local, wang2021cooperative, spiridonoff2020local, reisizadeh2020fedpaq, haddadpour2019local} where all participating clients aim to minimize the same objective function. However, ensuring convergence guarantees for \texttt{FedAvg} becomes notably more challenging in the presence of heterogeneity \citep{khaled2019first,khaled2020tighter,haddadpour2019convergence, li2019convergence}, thus necessitating additional assumptions on gradient and Hessian dissimilarity bounds \citep{li2019convergence,li2019convergence2, khaled2019first,karimireddy2020scaffold}. This difficulty arises primarily due to a ``client-drift'' effect, which is inherent to the \texttt{FedAvg} algorithm and has a detrimental impact on its convergence performance \citep{charles2020outsized, charles2021convergence}. As a result of the challenges posed by \texttt{FedAvg}, several alternative algorithms have been proposed to address its limitations. Notable examples include \texttt{FedProx} \citep{li2020federated}, \texttt{Scaffold} \citep{karimireddy2020scaffold}, \texttt{FedSplit} \citep{pathak2020fedsplit}, \texttt{FedDR} \citep{tran2021feddr}, \texttt{FedADMM} \citep{wang2022fedadmm}, \texttt{FedLin} \citep{mitra2021linear}, and \texttt{S-Local-SVRG} \citep{gorbunov}. Each of them introduces unique techniques and modifications to the original \texttt{FedAvg} algorithm, aiming to enhance convergence guarantees while handling communication cost concerns, statistical heterogeneity, client dropout, and sample complexity more effectively. 

   Applying federated learning (FL) to control systems introduces a novel research direction that comes with its own set of challenges. Control systems exhibit unique characteristics, such as non-iid and non-isotropic data, as well as system instability, which arise due to the dynamic nature of the systems. These characteristics pose specific challenges when attempting to leverage data from multiple systems for tasks such as system identification \citep{wang2022fedsysid} or control  synthesis \citep{ren2020federated}. 

    Although \cite{ren2020federated} address the model-free LQR tracking problem in a multi-agent setting, it  focuses on a significantly simpler scenario where all agents follow identical dynamics (i.e., homogeneous). In contrast, our work introduces new analysis techniques to achieve linear speedup when dealing with heterogeneous systems and multiple local updates per communication round.\footnote{In the time since this work was first made available online it has been used to develop new results on multi-task LQR design, including studies on meta-learning linear quadratic regulators \citep{toso2024meta, pan2025model, aravind2024moreau}, asynchronous federated LQR design \citep{toso2024asynchronous, zhao2025asynchronous}, learning LQR controllers from proxy systems \citep{ye2024convergence}, and domain randomization for the LQR problem \citep{fujinami2025policy}.  }

    \item \textbf{Policy Gradient (PG):} The policy gradient (PG) approach is a fundamental component of the success of reinforcement learning (RL) and plays a crucial role in policy optimization (PO). This approach directly optimizes the policy to improve system-level performances through gradient ascent steps. The concept of policy optimization has been influential in RL \citep{sutton1999policy} with some well-known algorithms such as \texttt{REINFORCE} \citep{williams1992simple}, trust-region policy optimization \texttt{TRPO} \citep{schulman2015trust}, actor-critic methods \citep{konda1999actor}, and proximal policy optimization \texttt{PPO} \citep{schulman2017proximal}. We highlight an important difference between  standard MDP  models and control models in RL. In control, one requires the policy to provide closed-loop stability, i.e., all trajectories of the system must converge for a given policy. In contrast, there are no analogous stability concerns in finite state-action MDP settings studied in RL.

     The extensive body of literature on policy optimization for reinforcement learning (RL) and its adaptability to the model-free setting paves the way for leveraging policy gradient methods in the pursuit of learning optimal control policies for classical control problems \citep{hu2022towards, perdomo2021stabilizing}. Despite the non-convex nature of the formulation involved in policy gradient methods, recent work \citep{fazel2018global, malik2019derivative, hambly2021policy, mohammadi2021convergence, gravell2020learning, jin2020analysis, ju2022model, perdomo2021stabilizing, lamperski2020computing, toso2024oracle, toso2025learning} has demonstrated global convergence in solving the model-free LQR problem via policy gradient methods in different settings. This global convergence is achieved due to certain properties of the quadratic cost function inherent in the LQR problem (e.g., smoothness and gradient dominance) as established in \cite{fazel2018global}. In contrast to the aforementioned work, which exclusively focus on the centralized control setting, our paper offers convergence guarantees for the multi-agent setting. In this context, each agent follows different dynamics, thereby distinguishing it from the simpler scenario proposed in \cite{ren2020federated}.

\item \textbf{Federated Reinforcement Learning (FRL):} The flexibility of policy gradient methods in the model-free RL setting has paved the way for a relatively recent research direction known as federated reinforcement learning (FRL), which aims to address practical implementation challenges of RL through the use of federated learning \citep{qi2021federated}. FRL focuses on learning a common value function \citep{wang2023federated, fabbro2023federated} or improving the policy by leveraging multiple RL agents interacting with similar environments~\citep{woo2023blessing, wang2024momentum}. The empirical evidence presented in the survey paper  \citep{qi2021federated} demonstrates the significant success of FRL in reducing sample complexity across various applications such as autonomous driving \citep{liang2022federated}, IoT devices \citep{lim2020federated}, resource management in networking \citep{yu2020deep}, and communication efficiency \citep{gatsis2022federated}. The FRL literature can be broadly grouped into two categories: (i) The \emph{homogeneous} setting, where all agents interact with the same MDP~\citep{khodadadian, woo2023blessing, lan2023improved, shen2023, liuMARL, fabbro2023federated, tian2024one, salgia2024sample, dal2025finite}, and (ii) the \emph{heterogeneous} setting, where the agents' MDPs can potentially differ in their reward functions and/or probability transition kernels~\citep{jinFRL, xie, wang2023federated, zhang2024finite, wang2024momentum, zhu2024towards, mangold2024, labbi2024federated, labbi2025global}. In particular, federated variants of policy gradient methods are analyzed in~\cite{xie, lan2023improved, wang2024momentum, zhu2024towards, labbi2025global, zhang2025gap}. It is important to note that all the papers in FRL discussed above consider MDPs with finite state and action spaces, where there is no notion of stability. In sharp contrast, our work considers a federated control problem involving dynamical systems with continuous state-action spaces, where we need to tackle the challenge of finding a common, stabilizing policy for all the agents at every iteration.

\end{itemize}

\section{Notation}  Given a set of matrices $\{S^{(i)}\}_{i=1}^M$, we denote $||S||_{\max} := \max_{i} ||S^{(i)}||$, and  $||S||_{\min} := \min_{i} ||S^{(i)}||$. All vector norms are Euclidean and matrix norms are spectral,  unless otherwise stated.

\section{Problem Setup}

Classical control approaches aim to design optimal controllers from a well-defined dynamical system model. The model-based LQR problem is  well understood and admits a convex solution. In this work, we consider the LQR problem but in the \textit{model-free setting}. Moreover, we consider a \emph{federated model-free LQR problem} in which there are $M$ agents, each with their own distinct but ``similar'' dynamics. Our goal is to collaboratively learn an optimal controller that minimizes an average quadratic cost. We seek to characterize the optimality of our solution as a function of the ``difference'' across the agent's dynamics. In what follows, we formally describe our problem of interest. 
\iffalse
We note that even in the absence of this heterogeneity, model-free LQR is a non-trivial problem~\citep{hu2022towards}.  
\fi

\textbf{Federated LQR:} Consider a system with $M$ agents. Associated with each agent is a linear time-invariant (LTI) dynamical system of the form
\begin{equation*}\label{eq:different_system}
    x^{(i)}_{t+1}=A^{(i)} x^{(i)}_t+B^{(i)} u^{(i)}_t, \quad x_0^{(i)} \sim \mathcal{D}, \quad \text{for all }i\in [M],
\end{equation*}
where $[M]:=1,\hdots, M$, and $A^{(i)} \in \mathbb{R}^{n_x\times n_x}$, $B^{(i)} \in \mathbb{R}^{n_x\times n_u}$. We assume each initial state $x_0^{(i)}$ is randomly generated from the same distribution $\mathcal{D}$. In the single-agent setting, the optimal LQR control policy is known to be linear and static. We denote the policy by  $u_t^{(i)} = {-} K_i^* x_t^{(i)}$, where for each agent, $K_i^\star$ solves
\begin{align}\label{eq:centralized_LQR}
K_i^*=\arg\min_{K} &\left\{  C^{(i)}(K):=\mc{E}\left[\sum_{t=0}^{\infty} x_t^{(i) \top}Q x_t^{(i) } + u_t^{(i)\top} R u_t^{(i)}\right]\right\}\notag \\
\text{s.t.}&\quad x_{t+1}^{(i)} = A^{(i)} x_t^{(i)} + B^{(i)} u_t^{(i)}, \ u_t^{(i)}= -K x_t^{(i)}, \ x_0^{(i)} \sim \mathcal{D},
\end{align}
where  $Q \in \mc{R}^{n_x \times n_x}$ and $R \in \mc{R}^{n_u \times n_u}$ are known positive definite matrices.
In our federated setting, the objective is to find an optimal common policy $\{ u_t\}_{t=0}^{\infty}$ to minimize the average cost of all the agents $C_{\text{avg}}(K) :=\frac{1}{M} \sum_{i=1}^M C^{(i)}(K)$ \textit{without knowledge of the system dynamics}, i.e., ($A^{(i)},B^{(i)}$). As mentioned above, classical  results~\citep{anderson2007optimal} from optimal control theory show that, given the system matrices $A^{(i)}$, $B^{(i)}$, $Q$ and $R$, the optimal policy can be written as a linear function of the current state. Thus, we consider a common policy of the form $u_t^{(i)} = - K x_t^{(i)}$. The objective of the federated LQR problem can be written as:
\begin{align}\label{eq:avg_LQR}
K^*=&\arg\min_{K} \left\{C_{\text{avg}}(K):=\frac{1}{M} \sum_{i=1}^M  \mc{E}\left[\sum_{t=0}^{\infty} x_t^{(i)\top}Q x_t^{(i)} + u_t^{(i)\top} R u_t^{(i)}\right]\right\}\notag \\
\text{s.t.}&\quad x_{t+1}^{(i)} = A^{(i)}x_t^{(i)} + B^{(i)} u_t^{(i)},\quad u^{(i)}_t= -Kx_t^{(i)}, x_0^{(i)} \sim \mathcal{D}.
\end{align}

The rationale for finding $K^*$ is as follows. Intuitively, when all agents have similar dynamics, $K^*$ will be close to each $K_i^*$. Thus, $K^*$ will serve to provide a good common initial guess from which each agent $i$ can then fine-tune/personalize (using only its own data) to converge \textit{exactly} to its own locally optimal controller $K_i^*$. The key here is that the initial guess $K^*$ can be obtained \textit{quickly} by using the \textit{collective data} of all the agents. We will formalize this intuition in Theorem~\ref{thm:main_fedlqr}. 

We make the standard assumption that for each agent, $(A^{(i)}, B^{(i)})$ is stabilizable. In addition, we make the following assumption  on the distribution of the initial state:
\begin{assume}\label{assume:initial_point} Given $\mu>0,$ for any $i \in [M],$ the initial state $x_0^{(i)} \sim \mathcal{D}$ and distribution $\mathcal{D}$ satisfy
\begin{equation*}\mc{E}_{x_0^{(i)}\sim \mathcal{D}}[x_0^{(i)}]=0, \;\ \mc{E}_{x_0^{(i)}\sim \mathcal{D}} [x_0^{(i)}x_0^{(i)\top}] \succ \mu \mc{I}_{d_x}, \text{ and } \lVert x_0^{(i)}\rVert \le H \text{ almost surely. }
\end{equation*} 
\end{assume}

We quantify the heterogeneity in the agent's dynamics through the following definition:
\begin{definition}
    (Bounded system heterogeneity) \label{assumption:bnd_sys_heterogeneity}
There exist positive constants $\epsilon_1$ and $\epsilon_2$ such that
$$\underset{i,j \in [M]}{\max} \lVert  A^{(i)} -A^{(j)}\rVert \leq \epsilon_1,\text{ and } 
\underset{i,j \in [M]}{\max}\lVert B^{(i)} -B^{(j)} \rVert\leq \epsilon_2.$$
\end{definition}
We assume that $\epsilon_1$ and $\epsilon_2$ are finite. Similar bounded heterogeneity assumptions are commonly made in FL~\citep{karimireddy2020scaffold,khaled2020tighter,reddi2020adaptive}. However, unlike typical FL works where one directly imposes heterogeneity assumptions on the agents' gradients, in our setting, we need to carefully characterize how heterogeneity in the system parameters $(A^{(i)},B^{(i)})$ translates to differences in the policy gradients; see Lemma~\ref{lem:gradient_het}.

\noindent\textbf{On the meaning of system similarity:} The notion of heterogeneity introduced in Definition \ref{assumption:bnd_sys_heterogeneity} measures similarity in terms of parameter deviations between system matrices. It is important to emphasize that similarity in the system parameters $(A^{(i)},B^{(i)})$ does not in general imply similar closed-loop behavior (see Propositions \ref{proposition:low_het1} and \ref{proposition:low_het2}). Even a small perturbation in system parameters can produce significant differences in the set of stabilizing controllers, and may not guarantee the existence of a common stabilizing controller across agents. Recent work \citep{stamouli2025policy} has proposed a refined similarity measure based on bisimulation-type metrics that quantify \emph{closed-loop} heterogeneity rather than parameter differences, and the main results of our paper extend to their setting.

Before providing our solution to the federated LQR problem, we first recap existing results on model-free LQR in the single-agent setting. 

\textbf{Single-agent setting:} When there is only one agent, i.e., $M=1$, let us denote the system matrix as $(A, B)$ and likewise, the cost functional by $C(K)$. If $(A,B)$ is known, the optimal controller $K^*$ can be computed by solving the discrete-time algebraic Riccati equation (ARE)~\citep{anderson2007optimal}. 

Strikingly, \cite{fazel2018global} show that policy gradient methods can find the globally optimal LQR policy $K^*$ despite the non-convexity of the problem. 
The policy gradient of the LQR problem can be expressed as:
$$
\nabla C(K)=2E_K\Sigma_K=2\left(\left(R+B^{\top} P_K B\right) K-B^{\top} P_K A\right) \Sigma_K,
$$
where $P_{K}$ is the positive definite solution to the  Lyapunov equation: $P_{K} = Q +K^{\top}RK + (A-BK)^{\top}P_{K}(A-BK),$
$
E_K:=\left(R+B^{\top} P_K B\right) K-B^{\top} P_K A
$, and $\Sigma_K:=\mathbb{E}_{x_0 \sim \mathcal{D}} [\sum_{t=0}^{\infty} x_t x_t^{\top}].$  The policy gradient method $K \leftarrow K - \eta \nabla C(K)$  will find the global optimal LQR policy, i.e., $K\rightarrow K^*$, provided that $\mc{E}_{x_0\sim \mathcal{D}}[x_0 x_0^{\top}]$ is full rank and an initial stabilizing policy is used. When the model is unknown, the analysis technique employed by~\cite{fazel2018global} is to construct near-exact gradient estimates from reward samples and show that the sample complexity of such a method is bounded polynomially in the parameters of the problem. 

In contrast to the single-agent setting, the heterogeneous, multi-agent scenario we consider here is considerably more difficult to analyze. First, designing an algorithm satisfying the iterative stability guarantees becomes a complex task. Second, since each agent in the system has its own unique dynamics and gradient estimates, it can be difficult to aggregate these directions in a manner that ensures the updating direction moves toward the average optimal policy $K^*$. Nonetheless, in the sequel, we will overcome these challenges and provide a  finite-time analysis of \texttt{FedLQR}.

\section{Necessity of the Low Heterogeneity Requirement} \label{sec:heterogeneity_regime}

In our main theorems, we require certain bounds on the parameters $\epsilon_1$ and $\epsilon_2$ that define the heterogeneity of the $M$ dynamical systems we work with. Here, we point out that, unlike standard federated learning settings, these bounds are \textit{necessary} for convergence. From a control and dynamical systems viewpoint, these bounds are perhaps intuitive: if the systems are too different, then there is no reason to believe there exists a common stabilizing controller, i.e., there is no solution to the problem~\eqref{eq:avg_LQR}. In what follows, we will formalize this point. To do so, let us define an ``instance" of our FedLQR problem via a parameter $M$ that characterizes the number of agents/systems and the set of corresponding system matrices $\{A^{(i)}, B^{(i)}\}_{i\in [M]}$.\footnote{Although technically the cost matrices $Q$ and $R$ are also part of a FedLQR problem formulation, they are not needed to establish the necessity of a low-heterogeneity requirement. As such, we do not include them here.}

We now prove a couple of simple impossibility results. Our first result shows that even when the input matrices are identical across agents, heterogeneity in the state transition matrices can lead to the non-existence of simultaneously stabilizing controllers, thereby rendering the FedLQR problem infeasible. 

\begin{proposition}\label{proposition:low_het1} There exists an instance of the FedLQR problem with $M=2$ and $\epsilon_2=0$, such that if $\epsilon_1 > 2$, then it is impossible to find a common state-feedback gain $K$ that simultaneously stabilizes both systems.
\end{proposition}

\textbf{Proof:}
Consider an instance with just two scalar systems defined by:
\[
x^{(1)}_{t+1} = \alpha x^{(1)}_t+u^{(1)}_t \quad \text{and} \quad x^{(2)}_{t+1} = -\alpha x^{(2)}_t+u^{(2)}_t,
\]
for some $\alpha>0$. By simple inspection, note that in this case $\epsilon_1=2\alpha$ and $\epsilon_2=0.$ Thus, $\epsilon_1> 2\Rightarrow \alpha>1$. Now for a controller $u_t^{(i)}=-kx_t^{(i)}$ to  stabilize  both systems, the spectral radius conditions are $|\alpha -k|<1$ and $|\alpha+k|<1$. Trivially, there exists no gain $k$ that satisfies both these requirements when $\alpha >1.$

To complement the above result, we now show that the effect of heterogeneity is not just limited to the state transition matrices. In particular, even when the state transition matrices are identical across agents, (arbitrarily small) heterogeneity in the input matrices can also lead to the non-existence of simultaneously stabilizing control gains. We formalize this below. 

\begin{proposition}\label{proposition:low_het2} There exists an instance of the FedLQR problem with $M=2$ and $\epsilon_1=0$, such that if $\epsilon_2 > 0$, then it is impossible to find a common state-feedback gain $K$ that simultaneously stabilizes both systems.
\end{proposition}

\textbf{Proof:}
Consider an instance with two scalar systems defined by:
\[
x^{(1)}_{t+1} = x^{(1)}_t+ \beta u^{(1)}_t \quad \text{and} \quad x^{(2)}_{t+1} = x^{(2)}_t - \beta u^{(2)}_t,
\]
for some $\beta$. By simple inspection, note that in this case $\epsilon_1=0$ and $\epsilon_2=2 \beta.$ 
Thus, $\epsilon_2> 0\Rightarrow \beta > 0$. Now for a controller $u_t^{(i)}=-kx_t^{(i)}$ to  stabilize  both systems, the spectral radius conditions are $|1-\beta k|<1$ and $|1+\beta k|<1$. Trivially, there exists no gain $k$ that satisfies both these requirements when $\beta > 0.$ This concludes the proof.

The above example suggests that in certain settings, we can tolerate no heterogeneity whatsoever in the input matrices. More generally, the main take-home message from this section is that the requirement of a ``low-heterogeneity regime'' is \textit{fundamental} to the problem and not merely an artifact of our analysis. We formally define a low heterogeneity regime in the sequel.

\section{The \texttt{FedLQR} algorithm}
In this section, we introduce our algorithm~\texttt{FedLQR}, formally described by Algorithm~\ref{alg:model_free_FedLQR}, to solve for $K^*$ in~\eqref{eq:avg_LQR}. First, we impose the following assumption regarding the algorithm's initial condition $K_0$:
\begin{assume}\label{assume:initial_controller}
We can access an initial stabilizing controller, $K_0$, which  stabilizes all  systems $\{(A^{(i)}, B^{(i)})\}_{i=1}^M$, i.e., the spectral radius $\rho (A^{(i)} -B^{(i)}K_0) < 1$ holds for all $i \in [M]$.
\end{assume}

\textbf{Algorithm description:}
At a high level, \texttt{FedLQR} follows the standard FL algorithmic template:  a server first initializes a global policy, $K_0$, which it sends to the agents. Each agent proceeds to execute  multiple PG updates using their local data. Once the local training is finished, agents transmit their model update to the server. The server aggregates the models and broadcasts an averaged model to the clients. The process repeats until a termination criterion is met. Prototypical FL algorithms that adhere to this structure include, for instance, \texttt{FedAvg}~\citep{khaled2020tighter} and \texttt{FedProx}~\citep{li2020federated}.

With this template in mind, we now dive into the details: \texttt{FedLQR} initializes the server and all agents with $K_{0,0}^{(i)} = K_0$ -- a controller that stabilizes all agent's dynamics.\footnote{To establish notation, $K_{n,l}^{(i)}$ refers to the local controller associated to agent $i$ after $n$ rounds of averaging, and after  $l^{\mathrm{{th}}}$ local iteration of policy gradient updates.} In each round $n$, starting from a common global policy $K_n$, each agent $i$ independently samples $n_s$ trajectories from its own system at each local iteration $l$ and performs approximate policy gradient updates  using the zeroth-order optimization procedure~\citep{fazel2018global} which we denote \texttt{ZO}; see line 7. For clarity, we present the explicit steps of using the zeroth-order method to estimate the true gradient in Algorithm~\ref{algorithm:gradient_estimation}, which will be discussed shortly. Between every communication round, each agent updates their local policy $L$ times. Such an $L$ is  chosen to balance between the benefit of information sharing and the cost of  communication. After $L$ local iterations, each agent $i$  uploads its local policy difference $\Delta_n^{(i)}$ (line 10) to the server. Once all differences are received, the server averages these differences $\{\Delta_n^{(i)}\}$ (line 12) to construct a new global policy $K_{n+1}.$ The whole process is repeated $N$ times.

\begin{algorithm}
\caption{Model-free Federated Policy Learning for the LQR (\texttt{FedLQR})} 
\label{alg:model_free_FedLQR}
\begin{algorithmic}[1]
\State \textbf{Input:} initial policy $K_0$, local step-size $\eta_l$ and global step-size $\eta_g.$
\State \textbf{Initialize} the server with $K_0$ and $\eta_g$
\State \textbf{for} $n=0, \ldots, N-1$ \textbf{ do} 
\State \quad\quad \textbf{for each system} $i \in [M]$ \textbf{ do} \
\State \quad\quad \quad\textbf{for} $l=0,\cdots, L-1$ \textbf{do} \
\State \quad\quad \quad \quad Agent $i$ initializes $K^{(i)}_{n,0}=K_n$\
\State \quad\quad \quad  \quad Agent $i$ estimates $\widehat{\nabla} C^{(i)}(K^{(i)}_{n,l}) = \texttt{ZO}(K^{(i)}_{n,l}, i)$ and updates local policy as\
\State \quad\quad \quad \quad  \quad  \quad\quad $K^{(i)}_{n,l+1} = K^{(i)}_{n,l} -\eta_l \widehat{\nabla} C^{(i)}(K^{(i)}_{n,l})$
\State \quad\quad \quad \textbf{end for}
\State \quad\quad \quad send $\Delta^{(i)}_n = K^{(i)}_{n,L}- K_{n}$ back to the server
\State \quad\quad \textbf{end for}\
\State \quad \quad Server computes and broadcasts global model $ K_{n+1} = K_n + \frac{\eta_g }{M} \sum_{i=1}^M \Delta^{(i)}_n$\
\State \textbf{end for}
\end{algorithmic}
\end{algorithm}

Zeroth-order optimization~\citep{conn2009introduction,nesterov2017random} provides a method of optimization that only requires oracle access to the function being optimized. Here, we briefly describe the details of our zeroth-order gradient estimation step\footnote{See Appendix \ref{sec:zeroth_order_optimization} for more details on  zeroth-order optimization.} in Algorithm~\ref{algorithm:gradient_estimation}. To obtain a gradient estimate at a given policy $K$, we sample trajectories from the $i$-th system $n_s$ times. At each time $s,$ we use the perturbed policy $\widehat{K}_s$ (line 3) and a randomly generated initial point $x_0 \sim \mathcal{D}$ to simulate the $i$-th closed-loop system for $\tau$ steps. Thus, we can approximately calculate the cost by adding the stage cost from the first $\tau$ time steps on this trajectory (line 4), and then estimating the gradient as in line 6.

\textbf{Discussion of Assumption~\ref{assume:initial_controller}:} Assumption~\ref{assume:initial_controller} is commonly adopted in the LQR~\citep{fazel2018global,dean2020sample,agarwal2019logarithmic,ren2020federated} and robust control literature~\citep{boyd1994linear,LuZD96,DGKF89}. In addition, there exist efficient ways to find such a stabilizing policy $K_0$; \citep{boyd1994linear} that addresses the model-based setting, while \citep{lamperski2020computing,perdomo2021stabilizing,zhao2022sample,toso2025learning} each addresses the problem in the single-agent model-free setting. Moreover, it is well-known that the sample complexity of finding an initial stabilizing policy only adds a logarithmic factor to that for solving the single-agent LQR problem \citep{zhao2022sample}. More recently, \cite{fujinami2025policy} have shown that solving a sequence of discounted LQR problems using the \emph{average} $C_{\mathrm{avg}}(K)$ cost yields a common stabilizing controller for all systems under low heterogeneity conditions.

\noindent\textbf{Choice of the number of local updates:} The parameter $L$ controls the standard communication, drift trade-off encountered in federated optimization. Executing more local steps before aggregation reduces communication frequency, but at the cost of amplifying the deviation between the aggregated update and the true global gradient direction. Conversely, choosing a smaller $L$ keeps the local iterates closer to the global descent direction and improves agreement with the centralized update, albeit requiring more communication rounds. In our setting, this trade-off becomes particularly delicate because drift between local controllers directly affects stability. We refer the reader to Lemmata \ref{lem:stability_local} and \ref{lem：drift_term} for additional details on how the number of local updates enters our analysis.

\textbf{Challenges in \texttt{FedLQR} analysis:} Although \texttt{FedLQR} is similar in spirit to \texttt{FedAvg}~\citep{li2019convergence,mcmahan2017communication} (in the supervised learning setting), it is significantly more difficult to analyze the convergence of \texttt{FedLQR} for the following reasons. 
\begin{itemize}
\item First, the problem we study is non-convex. Unlike most existing non-convex FL optimization results~\citep{karimireddy2020scaffold} which only guarantee convergence to stationary points, our work investigates whether \texttt{FedLQR} can find a globally optimal policy. 

\item Second, standard convergence analyses in FL~\citep{mcmahan2017communication, karimireddy2020scaffold,wang2020tackling,li2020federated} rely on a ``bounded gradient-heterogeneity" assumption. For the LQR problem, it is not clear a priori whether similar bounded policy gradient dissimilarity still holds. In fact, this is something we prove in Lemma~\ref{lem:gradient_het}.

\item Third, the randomness in FL usually comes from only one source: the data obtained by each agent are drawn i.i.d. from some distribution; we call this sample randomness. However, in \texttt{FedLQR}, there are randomness in the initial condition, and randomness from the smoothing matrices that show up in the gradient estimation process outlined in Algorithm~\ref{algorithm:gradient_estimation}. To reason about these different forms of randomness (that are intricately coupled), we provide a careful martingale-based analysis.
 
\item  Finally, we need to determine whether the solution given by \texttt{FedLQR} is meaningful, i.e., to decide whether the policy generated at each (local and global) iteration will stabilize all the systems. 
\end{itemize}
To tackle these difficulties, we first define a stability region in our setting comprising of $M$ heterogeneous systems as:
\begin{definition}\label{def:stable_set}(The stabilizing set) The stabilizing set is defined as $\mathcal{G}^0({\beta}) :=\cap_{i=1}^M \mathcal{G}^{(i)}({\beta})$ where
\begin{align*}
    \mathcal{G}^{(i)}({\beta}):= \left\{K: C^{(i)}(K) - C^{(i)}(K_i^*) \leq \beta \left( C^{(i)}(K_0) - C^{(i)}(K_i^*) \right)\right\}.  
\end{align*}
\end{definition}
As in~\cite{malik2019derivative}, $\mathcal{G}^0({\beta})$ is defined as the intersection of sub-level sets containing points $K$ whose cost gap is at most $\beta$ times the initial cost gap for all systems. We drop the dependence on $\beta$ when it is clear from the analysis and context. It was shown in~\cite{hu2022towards} that this is a compact set. Each sub-level set corresponds to a cost gap to agent $i$'s optimal policy $K_i^*$, which is at most $\beta$ times the initial cost gap $C^{(i)}(K_0) - C^{(i)}(K_i^*)$. Note that $\beta$ can be any positive finite constant. Since any finite cost function indicates that $K$ is a stabilizing controller, we conclude that any $K \in \mathcal{G}^{0}({\beta})$ stabilizes all the systems. Following from Assumption~\ref{assume:initial_controller}, there exists a constant $\beta$ such that $\mathcal{G}^0({\beta})$ is nonempty.  Moreover, it is worth remarking that the LQR cost function in the single-agent setting is \emph{coercive}. That is, the cost acts as a barrier function, ensuring that the policy gradient update remains within the feasible stabilizing set $\mathcal{G}^{(i)}({\beta})$. By defining the stabilizing set $\mathcal{G}^0({\beta})$ as above, the cost function $C^{(i)}(K)$ retains its coerciveness on $\mathcal{G}^{0}({\beta})$ for the federated setting considered in this paper. We also note that for our per-iteration stability analysis in Lemma \ref{lem:stability_local}, $\beta$ is properly set and kept fixed over iterations.

In order to solve the federated LQR problem and  provide convergence guarantees for \texttt{FedLQR}, we first need to recap some favorable properties of the  LQR problem in the single-agent setting that enables PG to find the globally optimal policy.

\begin{algorithm}[t]
\caption{Zeroth-order gradient estimation (\texttt{ZO})} 
\label{algorithm:gradient_estimation}
\begin{algorithmic}[1]
\State \textbf{Input:} $K$, number of trajectories $n_s$, trajectory length $\tau$, smoothing radius $r$, dimension $n_x$ and $n_u$, system index $i$.
\State \textbf{for} $s=1, \ldots, n_s$ \textbf{ do } \
\State \quad Sample a policy $\widehat{K}_s=K+U_s$, with $U_s$ drawn uniformly at random over matrices whose (Frobenius)  norm is $r$.\
\State \quad Simulate the $i$-th system for $\tau$ steps starting from $x_0 \sim \mathcal{D}$ using policy $\widehat{K}_s$. Let $\widehat{C}_s$ be the empirical estimate:
$\widehat{C}_s=\sum_{t=1}^{\tau} c_t,$
where $c_t: = x_t^{\top}\left( Q +\widehat{K}_s^{\top}R \widehat{K}_s\right)x_t$ on this trajectory.\
\State \textbf{end for}\
\State \textbf{Return} the estimate:
$\quad 
\widehat{\nabla} C(K)=\frac{1}{n_s} \sum_{s=1}^{n_s} \frac{n_x n_u}{r^2} \widehat{C}_s U_s.
$
\end{algorithmic}
\end{algorithm}

\section{Background on the centralized LQR using PG}\label{sec:background}
In the single-agent setting, it was shown that policy gradient methods (i.e., model-free) can produce the global optimal policy despite the  LQR problem being non-convex~\citep{fazel2018global}. We summarize the properties that make this possible and which we also exploit in our analysis.

\begin{lemma}\label{Lemma:lipschitz}(Local Cost and Gradient Smoothness) 
Suppose $K^{\prime}$ is such that $\|K^{\prime} -K\|  \leq h_{\Delta}(K)\ <\infty$. 
Then, the cost and gradient function satisfy:
\begin{equation*}
\begin{aligned}
\left|C\left(K^{\prime}\right)-C(K)\right| &\leq h_{\text {cost}}(K) \|K^{\prime} -K\|, \\
\left\|\nabla C\left(K^{\prime}\right)-\nabla C(K)\right\| \leq h_{\text {grad }}(K)\|\Delta\| &\text{ and }   \left\|\nabla C\left(K^{\prime}\right)-\nabla C(K)\right\|_F \leq h_{\text {grad}}(K)\|\Delta\|_F,
\end{aligned}
\end{equation*}
respectively, where $h_{\Delta}(K)$, $h_{\text {cost}}(K)$ and $h_{\text {grad}}(K)$ are positive scalars depending on $C(K)$.
\end{lemma}

\begin{lemma}\label{lemma:gradient_domination} (Gradient Domination)  Let $K^*$ be an optimal policy. Then, %following inequality
$$
C(K)-C\left(K^*\right) \leq \frac{\left\|\Sigma_{K^*}\right\|}{4\mu^2 \sigma_{\min }(R)}\|\nabla C(K)\|_F^2
$$
holds for any stabilizing controller $K$, i.e., any $K$ satisfying the spectral radius  $\rho(A-BK) < 1.$
\end{lemma}

For simplicity, we skip the explicit expressions in these lemmas for $h_{\Delta}(K)$, $h_{\text {cost}}(K)$, and $h_{\text {grad}}(K)$ as functions of the parameters of the LQR problem. Interested readers are referred to the appendix for full details. With  Definition~\ref{def:stable_set} of the stabilizing set in hand, we can define the following quantities:
 \begin{equation*}
 \bar{h}_{\text{grad}} := \sup_{K \in \mathcal{G}^{0}} {h}_{\text{grad}}(K), \ \bar{h}_{\text{cost}} := \sup_{K \in \mathcal{G}^{0}} {h}_{\text{cost}}(K), \text{ and  } \underline{h}_{\Delta} := \inf_{K \in \mathcal{G}^{0}} h_{\Delta}(K).
 \end{equation*}
With these quantities, we can transform the \emph{local} properties of the LQR problem discussed in  Lemmas~\ref{Lemma:lipschitz}--\ref{lemma:gradient_domination}  into   properties that hold over the \emph{global} stabilizing set $\mathcal{G}^0$. For convenience, we use letters with(under and over) bars such as $\bar{h}_{\text{grad}}$ to denote the global parameters. We  are now ready to present our main results of \texttt{FedLQR} in the next section. 

\begin{figure}[t]
\centering
  \includegraphics[width=0.5\linewidth]{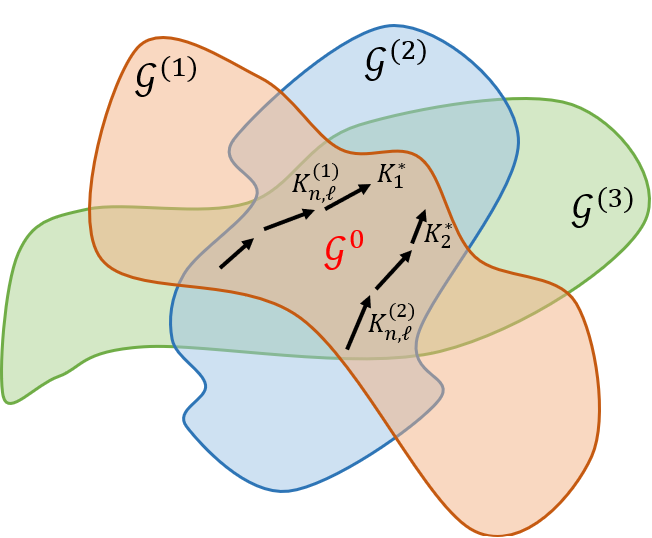}
  \caption{Illustration of the dynamics of the \texttt{FedLQR} algorithm for a setting with $M=3$ agents. Recall from Definition~\ref{def:stable_set}
that $\mathcal{G}^{(i)}$ represents the set of stabilizing policies for the $i$-th agent, and $\mathcal{G}^{0} = \bigcap_{i=1}^{3} \mathcal{G}^{(i)}$ is the common stabilizing set. In a low-heterogeneity regime, the agents' policy gradients are ``close" to one another (as formalized by Lemma~\ref{lem:gradient_het}), causing the sequence/trajectory of policies (controllers) generated by \texttt{FedLQR} to also evolve close to each other during the course of the algorithm. These trajectories $\{K^{(i)}_{n, l}\}_{i \in [3]}$ are illustrated by the black arrows. Crucially, initializing these trajectories within $\mathcal{G}^{0}$ causes them to always remain within $\mathcal{G}^{0}$, i.e., the set $\mathcal{G}^{0}$ is \emph{forward-invariant} for the trajectories generated by \texttt{FedLQR} - a fact established rigorously in Lemma~\ref{lem:stability_local}. This forward invariance property is the key to showing that the policies generated by \texttt{FedLQR} remain stabilizing for all systems.}
  \label{fig:PG_dyn}
\end{figure} 

\section{Main results}
To analyze the performance of \texttt{FedLQR} in the model-free case, we first need to examine its behavior in the model-based case. Although this is not our end goal, these results are of independent interest. 

\subsection{Model-based setting}\label{sec:model}
When $(A^{(i)},B^{(i)})$ are available, exact gradients can be computed, and so the \texttt{ZO} estimation scheme is no longer needed. In this case, the updating rule of \texttt{FedLQR} (Algorithm~\ref{alg:model_free_FedLQR}, line 12) reduces to
$$  K_{n+1} = K_n - \frac{\eta }{ML} \sum_{i=1}^M \sum_{l=0}^{L-1}\nabla C(K^{(i)}_{n,l}),$$
where $\eta := L \eta_g\eta_l.$ Intuitively, if two systems are similar, i.e., satisfy Assumption~\ref{assumption:bnd_sys_heterogeneity}, their exact policy gradient directions should not differ too much. We formalize this intuition as follows.
\begin{lemma}\label{lem:gradient_het}\textbf{(Policy gradient heterogeneity)}
For any $i, j \in [M]$ and $K \in \mathcal{G}^0$, we have:
\begin{align}\label{eq:gradient_het}
    ||\nabla C^{(i)}(K) -  \nabla C^{(j)}(K)|| \leq \epsilon_1 h^1_{\text{het}}(K) + \epsilon_2 h^2_{\text{het}}(K),
\end{align}
 where $h_{\text{het}}^1(K)$ and $h_{\text{het}}^2(K)$ are positive bounded functions depending on the parameters of the LQR problem.\footnote{For simplicity, we write $h_{\text{het}}^1, h_{\text{het}}^2$ as a function of only $K$ since only $K$ changes during the iterations while other parameters remain fixed.}
\end{lemma}
By Lemma~\ref{lem:gradient_het} (the proof of which is deferred to Appendix~\ref{sec:proof_model_based}), if $K$ belongs to a bounded set, the right-hand side of Eq.~\eqref{eq:gradient_het} is of the order  $\mathcal{O}(\epsilon_1 +\epsilon_2).$ In other words, the exact gradient direction of agent $i$ can be well-approximated by the gradient direction of agent $j$ when the heterogeneity constants $\epsilon_1$ and $\epsilon_2$ are small. This justifies why it is beneficial to use other agents' data under the
low-heterogeneity setting. Moreover, we can immediately conclude that the exact update direction of our \texttt{FedLQR} algorithm is also close to each agent's policy gradient direction based on Lemma~\ref{lem:gradient_het}. This fact is crucial for analyzing the convergence of \texttt{FedLQR} since we can map the convergence of \texttt{FedLQR} to that of the centralized LQR problem (with only one agent). However, Lemma~\ref{lem:gradient_het} alone is not sufficient to provide the final guarantees since we still need to consider the impact of multiple local updates i.e., $L>1$ (Algorithm~\ref{alg:model_free_FedLQR}, line 5) and stability concerns with heterogeneous systems -- i.e., is $\rho(A^{(i)}-B^{(i)}K_N)<1$ for all agents $i\in [M]$ \emph{and} all $n\in [N]$? Nevertheless, by overcoming these difficulties, we establish the convergence of \texttt{FedLQR} in the model-based setting as follows: 
\begin{theorem}\textbf{(Optimality in each agent's cost function)}\label{thm:one_local_model_based} When the heterogeneity level satisfies\footnote{The notation $\Bar{h}^3_{\text{het}}$ is a positive scalar depending on the parameters of the LQR problem; see Appendix~\ref{sec:proof_thm3} for full details.} $(\epsilon_1 \Bar{h}^1_{\text{het}} + \epsilon_2 \Bar{h}^2_{\text{het}})^2 \le \Bar{h}^3_{\text{het}}$ , there exist constant step-sizes $\eta_g$ and $\eta_l$ such that \texttt{FedLQR} enjoys the following performance guarantees over  $N$ rounds:
\begin{align*}
C^{(i)}(K_{N}) -& C^{(i)}(K_i^*)  \le \left(1 -\frac{\eta \mu^2 \sigma_{\min }(R)}{\left\|\Sigma_{K_i^*}\right\|} \right)^N(C^{(i)}(K_0) - C^{(i)}(K_i^*)) +c_{\text{uni},1}\times \mathcal{B}(\epsilon_1, \epsilon_2),
\end{align*}
with $\mathcal{B}(\epsilon_1, \epsilon_2):= \frac{\upsilon \left\|\Sigma_{K_i^*}\right\|}{4\mu^2 \sigma_{\min}(R)} (\epsilon_1 h^{1}_{\text{het}} + \epsilon_2 h^2_{\text{het}})^2$, where $\bar{h}^{1,2}_{\text{het}} := \sup_{K \in \mathcal{G}^{0}} {h}^{1,2}_{\text{het}}(K)$, $ \upsilon :=\min\{n_x,n_u\}$, and $c_{\text{uni},1}$ is a universal constant. Moreover, we have $K_n \in \mathcal{G}^0$ for all $n =0,\cdots,N$.
\end{theorem}
\textbf{Main Takeaways:}  Theorem~\ref{thm:one_local_model_based} reveals that the output $K_n$ of \texttt{FedLQR} can stabilize all $M$ systems at each round $n$. However,
 \texttt{FedLQR} can only converge to a ball of radius $\mathcal{B}(\epsilon_1, \epsilon_2)$ around each system's optimal controller $K_i^*$, regardless of the choice of the step-sizes. The term $\mathcal{B}(\epsilon_1, \epsilon_2)$ captures the effect of heterogeneity and becomes zero when each agent follows the same system dynamics, i.e., $\epsilon_1 =\epsilon_2 =0.$ When there is no heterogeneity, the convergence rate matches the rate of the centralized setting~\citep{fazel2018global} up to a constant factor. But, since there is no noise introduced by the zeroth-order gradient estimate, there is no expectation of obtaining a benefit from collaboration. Nonetheless, understanding the model-based setting provides valuable insights for exploring the model-free setting. The proof of Theorem~\ref{thm:one_local_model_based}  is given in Appendix~\ref{sec:proof_thm3}.  We are now ready to provide the convergence guarantees for  \texttt{FedLQR} with respect to solving Eq.~\eqref{eq:avg_LQR} -- the average cost. 

\begin{corollary}\label{coro:global convergence}
\textbf{(Optimality in average cost function)}  When the heterogeneity level satisfies\footnote{The notation $\Bar{h}^3_{\text{het}}$ is a positive scalar depending on the parameters of the LQR problem; see Appendix~\ref{sec:proof_thm3} for full details.} $(\epsilon_1 \Bar{h}^1_{\text{het}} + \epsilon_2 \Bar{h}^2_{\text{het}})^2 \le \Bar{h}^3_{\text{het}}$, after $N$ rounds, \texttt{FedLQR} enjoys the following optimality gap in average cost function across all $M$ agents:
\begin{align*}%\label{eq:distance_between_two}
C_{\text{avg}}(K_N)-C_{\text{avg}}(K^*) &\le
 \left(1 -\frac{\eta \mu^2 \sigma_{\min }(R)}{\max_i\|\Sigma_{K^*_i}\|} \right)^N \sup_{i \in [M]} (C^{(i)}(K_0) - C^{(i)}(K_i^*)) +c_{\text{uni},1}\times \mathcal{B}(\epsilon_1, \epsilon_2).
\end{align*}
\end{corollary}

The main message conveyed by Corollary~\ref{coro:global convergence} is that \texttt{FedLQR} can converge to a ball around the average optimal controller $K^*$ with a linear convergence rate. The size of the ball depends on the  system heterogeneity level, i.e., $\epsilon_1$ and $\epsilon_2$. Combining Theorem~\ref{thm:one_local_model_based} and Corollary~\ref{coro:global convergence}, we infer that \texttt{FedLQR} not only approximates each system's optimal controller $K_i^*$ but also approximately converges toward the average optimal controller $K^*$ when the underlying $M$ systems are close (in the sense of Definition \ref{assumption:bnd_sys_heterogeneity}). The primary distinction between converging to $K_i^*$ and $K^*$ lies in the linear convergence rate. Compared to converging to $K_i^*$, where the linear converge rate depends only on system $i$'s parameter $\left\|\Sigma_{K_i^*}\right\|$, the rate in converging to $K^*$ depends on all systems' parameters, i.e., $\left\{\left\|\Sigma_{K_i^*}\right\|\right\}_{i=1}^N$. See Appendix~\ref{sec:proof of coro} for a comprehensive proof.

\textbf{How to ensure \texttt{FedLQR}'s stability?}  We briefly discuss our proof technique for ensuring the iterative stability guarantees. The main idea is to leverage an inductive argument. We start from a stabilizing global policy $K_n \in \mathcal{G}^0.$ We aim to show that the next global policy $K_{n+1}$ is stabilizing; see Fig.~\ref{fig:PG_dyn} for an illustration. This is achieved by demonstrating that $K_{n+1}$ can reduce each system's cost function compared to $K_n$. To achieve this goal, we take the following steps: (1) at each iteration, initiate from the globally stabilizing controller computed at the previous iterate, (2) determine a small global step-size such that inequalities in Section~\ref{sec:background} can be applied; (3) use Lemma~\ref{lem:gradient_het} to provide a descent direction to reduce each system's cost function; (4) bound the drift term $\frac{1}{ML}\sum_{i=1}^M\sum_{l=0}^{L-1} \|K_{n,l}^{(i)}-K_n\|^2$. Step (4) can be accomplished using a  small local step-size $\eta_l$ such that each local policy is a small perturbation of the global policy $K_n$. 

\noindent\textbf{On the role of control variates for local drift mitigation:} Control-variate techniques such as proposed in \textsc{Scaffold} \citep{karimireddy2020scaffold} are designed to correct the optimization drift that arises in federated optimization under non-IID data by making each local update an approximately unbiased estimator of the global gradient. However, we emphasize that such methods do not eliminate the intrinsic gap between the minimizer of the averaged LQR cost $C_{\textrm{avg}}(K)$ and the individual optima $K_i^\star$ when the underlying dynamical systems differ. This limitation reflects a fundamental distinction between control-variate methods in classical FL and their applicability to our setting of dynamical systems. In standard FL, bias correction can ensure that the global iterates closely follow an ``idea'' gradient-descent step on a strongly convex objective, enabling the use of gradient-domination inequalities to guarantee fast convergence. In contrast, for heterogeneous LQR systems, the averaged cost $C_{\textrm{avg}}(K)$ does not, in general, correspond to the LQR cost of any single ``averaged'' dynamical system, and hence, one cannot directly appeal to gradient-domination properties like in the optimization setting to achieve fast linear rates. An important exception arises when all agents share identical dynamics but differ only through their cost matrices: in that case, $C_{\textrm{avg}}$ reduces to the LQR cost associated with $(A, B, \bar{Q}, \bar{R})$, where $\bar{Q}= (1/M) \sum_{i \in [M]} Q_i$ and $\bar{R}= (1/M) \sum_{i \in [M]} R_i$, making gradient-domination applicable and mitigating the local drift as discussed in \cite{zhu2024towards}.

Equipped with these model-based results, we are ready to present our main results of the model-free setting.

\subsection{Model-free setting}
We now analyze \texttt{FedLQR}'s convergence in the model-free setting, where the policy gradient steps are approximately computed using  zeroth-order optimization (Algorithm~\ref{algorithm:gradient_estimation}), without knowing the true dynamics, i.e., $A^{(i)},B^{(i)}$ are not available and so $\nabla C^{(i)}(K^{(i)})$ can't be directly computed. The key point in this setting is to bound the gap between the estimated gradient and the true gradient. In the centralized setting~\citep{fazel2018global}, the gap can   be made arbitrarily accurate with enough  trajectory samples $n_s$, sufficiently long trajectory length $\tau$, and small smoothing radius $r$. 

We aim to achieve a sample complexity reduction for each agent by utilizing data from other similar but non-identical systems with the help of the server. This presents a significant challenge, as averaging gradient estimates from multiple agents may not necessarily reduce the variance even for homogeneous systems due to the high correlation between local gradient estimates. This challenge is compounded in our case as the gradient estimates are not only \emph{correlated} but also come from \emph{non-identical systems}. As a result, the variance reduction and sample complexity reduction for the \texttt{FedLQR} algorithm is not obvious a priori. After addressing these challenges  using a  martingale-based analysis, we show that one can establish variance reduction for our setting as well. This is formalized in the next result:

\begin{lemma}\label{lem:variance_reduction}
\textbf{(Variance Reduction)} Suppose the smoothing radius $r$ and trajectory length $\tau$ from Algorithm~\ref{algorithm:gradient_estimation} satisfy  $r \leq h_r\left(\frac{\epsilon}{4}\right)$ and  $\tau\geq h_{\tau}\left(\frac{r \epsilon}{4 n_x n_u}\right)$, respectively.\footnote{The notation $h_r$, $h_{\tau}$, $h_{\text {sample,trunc}}$ and $h_r^{\prime}$ in Lemma~\ref{lem:variance_reduction} and Theorem~\ref{thm:main_fedlqr} are  polynomial functions of the LQR problem, depending on  $\epsilon$. For simplicity, we defer their definition to the appendix.} Moreover, suppose the  sample size satisfies:\footnote{For the convenience of comparison with  existing literature, we use the same notation as in ~\cite{fazel2018global}.}
\begin{equation}\label{eq:sample_requirement} n_{s} \geq \frac{h_{\text {sample,trunc }}\left(\frac{\epsilon}{4}, \frac{\delta}{ML}, \frac{H^2}{\mu}\right)}{ML}.\end{equation} Then, when $K_n \in \mathcal{G}^0$, with probability $1-\delta$, the estimated gradients satisfy:
$$\left\|\frac{1}{ML} \sum_{i=1}^M\sum_{l=0}^{L-1}\left[\widehat{\nabla} C^{(i)}(K^{(i)}_{n,l})-\nabla C^{(i)}(K^{(i)}_{n,l})\right]\right\|_F \leq \epsilon.
$$
\end{lemma}

We prove this result and provide the definition of the  parameters of $h_{sample,trunc}$  in Appendix~\ref{sec:proof_model_free}. The most important information conveyed by our variance reduction lemma is that each agent at each local step only needs to sample $\frac{1}{ML}$ fraction of samples required in the centralized setting. Notably, this lemma plays an important role in showing that \texttt{FedLQR} can help improve the sample efficiency. Equipped with Lemma~\ref{lem:variance_reduction}, we now present the main convergence guarantees for \texttt{FedLQR}:

\begin{theorem}~\label{thm:main_fedlqr} \textbf{(Model-free)}
Suppose the trajectory length satisfies $
\tau \geq h_{\tau}\left(\frac{r \epsilon^{\prime}}{4 n_x n_u}\right),$ the smoothing radius satisfies $r\le  h^{\prime}_r\left(\frac{\epsilon^{\prime}}{4}\right),$ and the sample size of each agent $n_s$ satisfies Eq.~\eqref{eq:sample_requirement}
with $\epsilon^{\prime} =\frac{4 \left\| \Sigma_{K_i^\star}  \right\|}{\mu^2 \sigma_{\min}(R)} \epsilon.$ When the heterogeneity level satisfies $(\epsilon_1 \Bar{h}^1_{\text{het}} + \epsilon_2 \Bar{h}^2_{\text{het}})^2 \le \Bar{h}^3_{\text{het}}$, 
 then, given any $\delta \in (0,1),$ with probability $1-\delta$, there exist constant step-sizes $\eta_g$ and $\eta_l$,  which are independent of $\epsilon^{\prime}$, such that \texttt{FedLQR} enjoys the following performance guarantees:\\

\begin{enumerate}
\item \textbf{(Stability of the global policy)} The global policy at each round $n$ is stabilizing, i.e., $K_{n} \in \mathcal{G}^0$;
\item \textbf{(Stability of the local policies)} The local policies satisfy $K_{n,l}^{(i)} \in \mathcal{G}^0$ for all $i$ and $l$;
\item \textbf{(Convergence rate)} After  $N \ge \frac{\left\|\Sigma_{K_i^*}\right\|}{\eta\mu^2 \sigma_{\min}(R)}\log\left(\frac{2(C^{(i)}(K_0) - C^{(i)}(K_i^*))}{\epsilon^{\prime}}\right)$ rounds, we have \begin{equation}\label{eq:overall_converge}
    C^{(i)}(K_N) - C^{(i)}(K_i^*) \le \epsilon^{\prime} + 2 \mathcal{B}(\epsilon_1, \epsilon_2),  \forall  i\in [M],
\end{equation}
where $\mathcal{B}(\epsilon_1, \epsilon_2)$ is defined in Theorem~\ref{thm:one_local_model_based}.
\end{enumerate}
\end{theorem}
This theorem establishes the finite-time convergence guarantees for \texttt{FedLQR}.  The first two points in Theorem~\ref{thm:main_fedlqr} provide the iterative stability guarantees of \texttt{FedLQR}, i.e., the trajectories of \texttt{FedLQR} will always stay inside the stabilizing set $\mathcal{G}^0$. The third point implies that when heterogeneity is small, i.e., $\mathcal{B}(\epsilon_1, \epsilon_2)$ is negligible, \texttt{FedLQR} converges to each system's optimal policy with a linear speedup w.r.t. the number of agents $M$, which we discuss further next.

\textbf{Discussion:} For a fixed desired precision $\epsilon,$ we denote $N$ to be the number of rounds such that the first term $\epsilon^{\prime}$ in Eq~\eqref{eq:overall_converge} is smaller than $\epsilon$. In what follows, we focus on analyzing the total sample complexity of \texttt{FedLQR} for each agent, which can be calculated by $N \times L\times n_s$. Note that $N$, in our case, is in the same order as the centralized setting. However, in terms of the sample size $n_s$ requirement at each local step, it is only a $\frac{1}{ML}$-fraction of that needed in the centralized setting, as presented in the variance reduction Lemma~\ref{lem:variance_reduction}. Therefore, in a low-heterogeneity regime, where $\mathcal{B}(\epsilon_1,\epsilon_2)$ is negligible, our \textit{\texttt{FedLQR} algorithm} reduces the sample complexity of learning the optimal LQR policy by $\tilde{\mathcal{O}}(\frac{1}{M})$ of the centralized setting \citep{fazel2018global, malik2019derivative}.\footnote{In~\cite{fazel2018global}, the sample complexity of policy gradient with one-point zeroth-order estimation is $\tilde{\mathcal{O}}(\frac{1}{\epsilon^4})$, this was later improved to $\tilde{\mathcal{O}}(\frac{1}{\epsilon^2})$ by~\cite{malik2019derivative}. We compare our results to the refined analysis in~\cite{malik2019derivative}.} Specifically, \texttt{FedLQR} improves the sample cost required by each agent from $\tilde{\mathcal{O}}(\frac{1}{\epsilon^2})$ to $\tilde{\mathcal{O}}(\frac{1}{M\epsilon^2})$ up to a small heterogeneity bias term. This result is highly desirable since the number of agents in FL is usually large; leading to a significant speedup due to collaboration.

It is important to mention that our results also capture the cost of federation embedded in the term $\mathcal{B}(\epsilon_1, \epsilon_2)$. That is when two systems exhibit significant differences from each other, leveraging data across them may not be beneficial in finding a common stabilizing policy that applies to both. \textit{In summary, Eq.~\eqref{eq:sample_requirement}--\eqref{eq:overall_converge} provide an explicit interplay between the price of heterogeneity and the benefit of collaboration}. The trade-off in Theorem~\ref{thm:main_fedlqr} is explored in the simulation study presented in the next section.

\section{Numerical Results} \label{sec:numerical_results}
The following section describes the experimental setup and results for \texttt{FedLQR} in the model-free setting\footnote{Code can be downloaded from \url{https://github.com/jd-anderson/FedLQR}}.

\subsection{System Generation}
Numerical experiments are conducted to illustrate and evaluate the effectiveness of \texttt{FedLQR}  (Algorithm \ref{alg:model_free_FedLQR}). The simulations involve different and unstable dynamical systems described by discrete-time linear time-invariant (LTI) models, as in \eqref{eq:different_system}, where each system has $n_x=3$ states and $n_u=3$ inputs. To generate different systems while respecting the bounded heterogeneity assumption (Assumption \ref{assumption:bnd_sys_heterogeneity}), the following steps are followed: 
\begin{enumerate}
    \item Given  nominal system matrices ($A_0$, $B_0$),  generate random variables $\gamma_1^{(i)} \sim \mathcal{U}(0,\epsilon_1)$ and $\gamma_2^{(i)} \sim \mathcal{U}(0,\epsilon_2)$, $\forall i \in [M]$, with $\epsilon_1$ and $\epsilon_2$ being predefined dissimilarity parameters.
    \item The  random variables generated above are combined with modification masks $Z_1 \in \mathbb{R}^{3\times 3}$ and $Z_2 \in \mathbb{R}^{3\times 3}$ to generate the different systems matrices $(A^{(i)}, B^{(i)})$ for all $i \in [M]$.
    \item The  systems $(A^{(i)},B^{(i)})$ for $0<i\le M$ are then constructed by perturbing the nominal systems according to: $A^{(i)} = A_0 + \gamma_1^{(i)}Z_1$ and $B^{(i)} = B_0 + \gamma_2^{(i)}Z_2$, where $Z_1$ and $Z_2$ are defined in step 2.
    \item The nominal matrices are included in the set of generated systems as $(A^{(1)}, B^{(1)}) = (A_0, B_0)$.
\end{enumerate}

In particular, we consider 

$$
A_0=\begin{bmatrix}
     1.20 & 0.50 & 0.40\\
     0.01 &  0.75 & 0.30\\
     0.10 &   0.02 &  1.50
    \end{bmatrix}, \;\ B_0= I_{3}, \;\ Q= 2 I_{3}, \text{ and } R= \frac{1}{2} I_{3},
$$
for the nominal system matrices and cost matrices respectively.

The  optimal controller for the nominal system $(A^{(1)},B^{(1)})$ is  
$$K^{*}_1 = \begin{bmatrix}
    1.0056  &   0.4293 &   0.3570\\
    0.0262  &  0.6239 &   0.2657\\
    0.1003 &   0.0298 &   1.2960
\end{bmatrix},$$
and was obtained by solving the discrete algebraic Riccati equation (DARE).

\subsection{Algorithm Parameters}
For the gradient estimation step in the zeroth-order algorithm (Algorithm \ref{algorithm:gradient_estimation}), we set the initial state for cost computation as a random sample from a standard normal distribution, denoted as $\mathcal{D}\stackrel{d}{=}\mathcal{N}(0,I_{3})$, for all systems $i \in [M]$. Additionally, we consider $n_s=5$ trajectories, where each trajectory has a rollout length of $\tau=15$, and we set the smoothing radius $r=0.1$ for the zeroth-order gradient estimation. 

Throughout our simulations, we consider the following initial stabilizing controller $K_0 = 1.62I_{3}$ (Line 1 in Algorithm \ref{alg:model_free_FedLQR}). Note that although the control action $u^{(i)}_t = -K_0x^{(i)}_t$ may not be optimal for any of the $M$ systems. For example, the suboptimality of $K_0$ applied to the nominal system is evidenced by its cost of $C^{(1)}(K_0)=18.4049$, compared to the optimal cost of $C^{(1)}(K^*_1) = 9.5220$, when computed from an initial state $x^{(1)}_0=[1\;\ 1 \;\ 1]^{\top}$ and time horizon $T=500$.  However, it is important to note that $K_0$ is still stabilizing all $M$ systems. Note that we will use $K_0$ as the initial controller for all of the experiments in this paper.

\subsection{Experiments}
To assess the performance of \texttt{FedLQR}, we evaluate the normalized gap between the current cost $C^{(1)}(K_n)$ of the nominal system when using the common stabilizing controller $K_n$ and its corresponding optimal cost $C^{(1)}(K^{*}_1)$. This metric is represented as $\frac{C^{(1)}(K_n) - C^{(1)}(K^{*}_1)}{C^{(1)}(K^{*}_1)}$ for each global iteration $n \in [N]$. In our experiments, we set the step sizes as $\eta_g=1\times 10^{-2}$, with an adaptive decrease of $0.05\%$ per global iteration, and $\eta=1\times 10^{-4}$. Further details regarding other parameters, such as the number of systems $M$, heterogeneity levels $(\epsilon_1,\epsilon_2)$ are provided in the figures and the subsequent discussion.

\begin{figure}
    \centering
    \begin{subfigure}[b]{0.5\textwidth}
        \centering
        \includegraphics[height=1.7in]{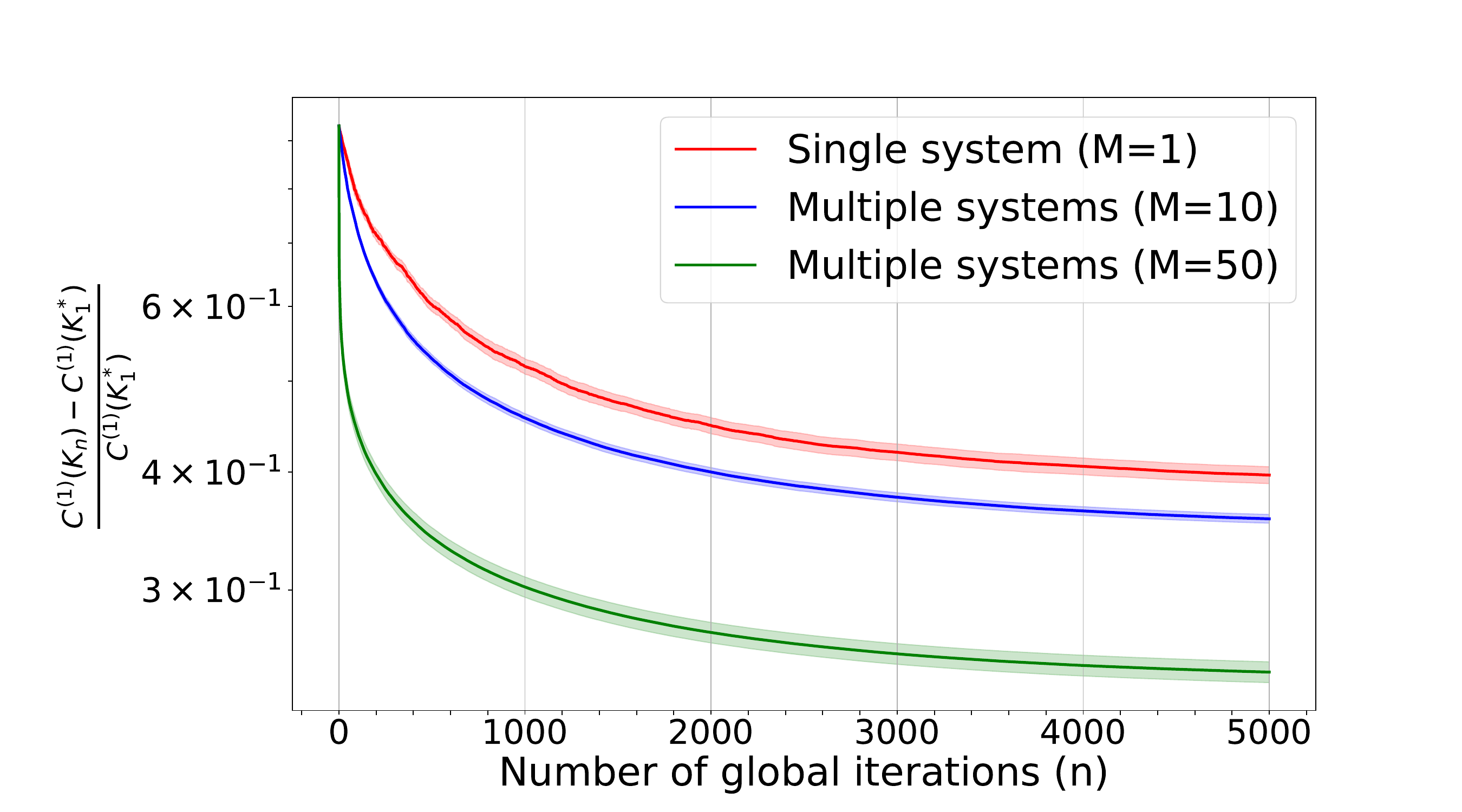}
        \caption{Varying $M$ with $\epsilon_1=\epsilon_2=0.5$.}
    \end{subfigure}%
    \begin{subfigure}[b]{0.5\textwidth}
        \centering
        \includegraphics[height=1.7in]{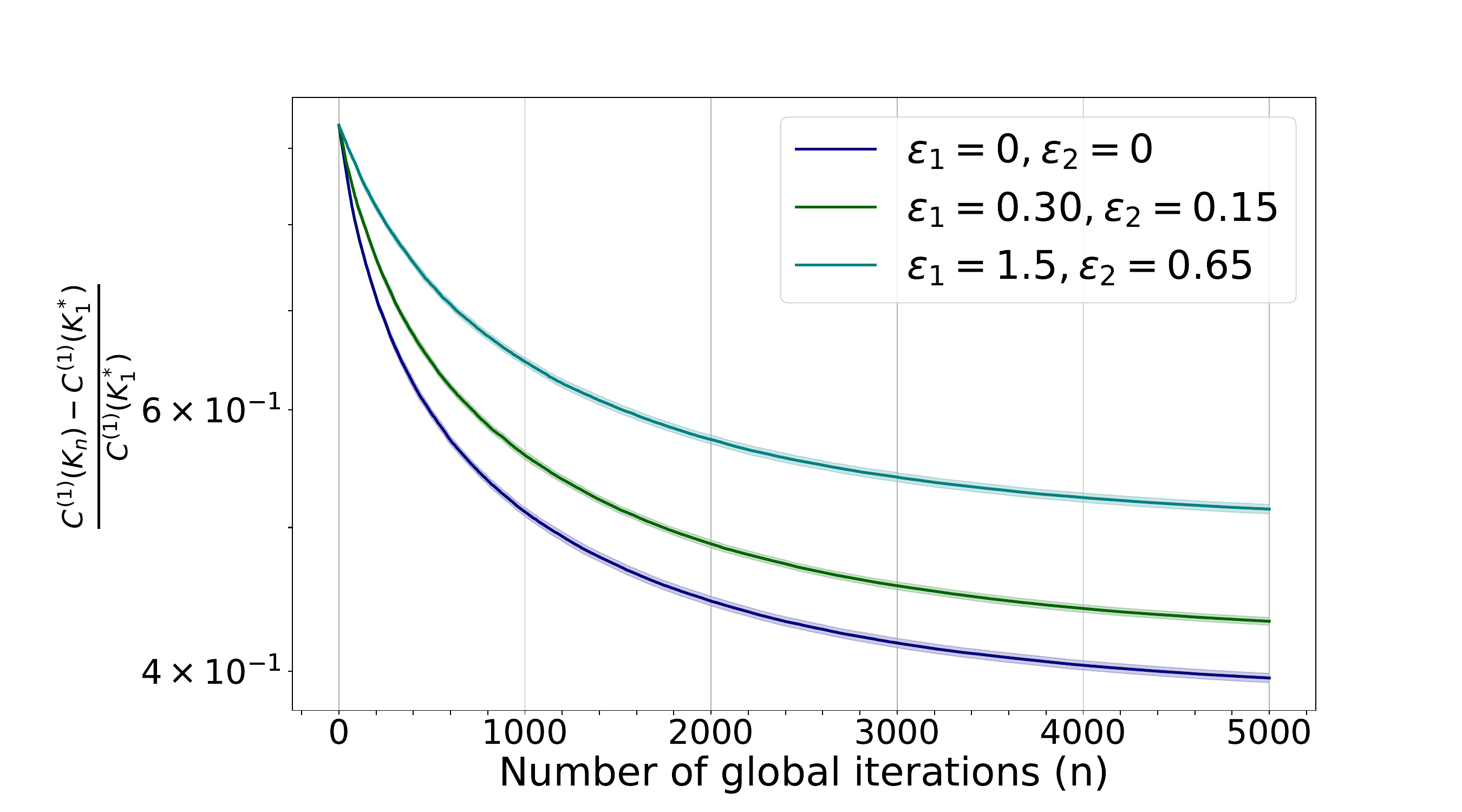}
        \caption{Varying $\epsilon_1$, $\epsilon_2$, with $M = 10$.}
    \end{subfigure}
    \caption{Gap between the current and optimal cost with respect to the number of global iterations.}\label{fig:FedLQR}
\end{figure}

Figures \ref{fig:FedLQR}-(a,b) present the normalized distance between the current cost associated with the common stabilizing controller and the optimal cost for the nominal system, plotted with respect to the number of global iterations. These figures demonstrate the impact of varying the number of systems $M$ and the heterogeneity parameters $(\epsilon_1,\epsilon_2)$ on the convergence and performance of Algorithm \ref{alg:model_free_FedLQR}.

In Figure \ref{fig:FedLQR}-(a), we specifically investigate the effect of the number of systems $M$ participating in the collaboration to compute a common controller $K^*$ on the convergence of our algorithm. In this analysis, we set the heterogeneity parameters as $\epsilon_1=0.5$ and $\epsilon_2=0.5$ and consider modification masks $Z_1 = Z_2 = I_3$. The figure reveals a noticeable reduction in the gap between the current and optimal cost as the number of participating systems $M$ increases. This numerical result aligns with our  theoretical findings, which indicate that the number of samples required to achieve reliable estimation for the cost function's gradient can be scaled down with the number of systems participating in the collaboration. Consequently, as the number of systems involved increases, there is a considerable reduction in the gap between the common computed controller and the optimal one.

Figure \ref{fig:FedLQR}-(b) illustrates the influence of the heterogeneity parameters $(\epsilon_1,\epsilon_2)$ on the convergence rate of Algorithm \ref{alg:model_free_FedLQR}. In this analysis, we set the number of systems as $M=10$, and the modification masks $Z_1=\texttt{diag}([3.5\;\ 1\;\ 0.1])$ and $Z_2=\texttt{diag}([1.5\;\ 0.1\;\ 1])$. Consistent with our theoretical findings, we observe that an increase in the dissimilarity among the systems results in a significant gap between the common and optimal controller. This discrepancy arises due to the additive effect of system heterogeneity on the convergence rate of our algorithm, as elaborated in Theorem \ref{thm:main_fedlqr}.

\section{Conclusions and Future Work}

We investigated the problem of learning a common and optimal LQR policy with the objective of minimizing an average quadratic cost. The primary focus of this paper was to thoroughly examine and provide comprehensive answers to the following questions: (i) Is the learned common policy stabilizing for all agents? (ii) How close is the learned common policy to each agent’s own optimal policy? (iii) Can each agent learn its own optimal policy faster by leveraging data from all agents? To address these questions, we proposed a federated and model-free approach,  \texttt{FedLQR}, where $M$ heterogeneous systems collaborate to learn a common and optimal policy while keeping the system's data \emph{private}. Our analysis tackles numerous technical challenges, including system heterogeneity, multiple local gradient descent updates, and stability. We have demonstrated that \texttt{FedLQR}  produces a common policy that stabilizes all systems and converges to the optimal policy (Theorem \ref{thm:main_fedlqr}) of each agent up to a heterogeneity-induced bias term. Furthermore, \texttt{FedLQR} achieves a reduction in sample complexity proportional to the number of participating agents $M$ (Lemma \ref{lem:variance_reduction}). We have also provided numerical results to effectively showcase and evaluate the performance of our \texttt{FedLQR} approach in a model-free setting. Future work will address the assumption of requiring full-state information to extend our results to the Linear Quadratic Gaussian (LQG) problem in a federated setting. We are currently investigating data-driven and system-theoretic metrics for heterogeneity, as well as personalization-based methods to mitigate the impact of heterogeneity on the performance of the proposed approach.

\section*{Acknowledgments}

The authors thank the anonymous reviewers for their thorough and constructive feedback, which significantly improved this work.  Leonardo F. Toso is funded by the Center for AI and Responsible Financial Innovation (CAIRFI) Fellowship and by the Columbia Presidential Fellowship. James Anderson is partially funded by NSF grants ECCS 2144634 and 2231350 and the Columbia Center of AI Technology in collaboration with Amazon.

\bibliography{main}

@inproceedings{fazel2018global,
  title={Global convergence of policy gradient methods for the linear quadratic regulator},
  author={Fazel, Maryam and Ge, Rong and Kakade, Sham and Mesbahi, Mehran},
  booktitle={International conference on machine learning},
  pages={1467--1476},
  year={2018},
  organization={PMLR}
}

@book{anderson2007optimal,
  title={Optimal control: linear quadratic methods},
  author={Anderson, Brian DO and Moore, John B},
  year={2007},
  publisher={Courier Corporation}
}

@article{ren2020federated,
  title={Federated {LQR}: Learning through Sharing.},
  author={Ren, Zhaolin and Zhong, Aoxiao and Zhou, Zhengyuan and Li, Na},
  journal={arXiv preprint arXiv:2011.01815v1},
  year={2020}
}

@article{gravell2020learning,
  title={Learning optimal controllers for linear systems with multiplicative noise via policy gradient},
  author={Gravell, Benjamin and Esfahani, Peyman Mohajerin and Summers, Tyler},
  journal={IEEE Transactions on Automatic Control},
  volume={66},
  number={11},
  pages={5283--5298},
  year={2020},
  publisher={IEEE}
}

@article{agarwal2019logarithmic,
  title={Logarithmic regret for online control},
  author={Agarwal, Naman and Hazan, Elad and Singh, Karan},
  journal={Advances in Neural Information Processing Systems},
  volume={32},
  year={2019}
}

@book{boyd1994linear,
  title={Linear matrix inequalities in system and control theory},
  author={Boyd, Stephen and El Ghaoui, Laurent and Feron, Eric and Balakrishnan, Venkataramanan},
  year={1994},
  publisher={SIAM}
}

@article{perdomo2021stabilizing,
  title={Stabilizing dynamical systems via policy gradient methods},
  author={Perdomo, Juan and Umenberger, Jack and Simchowitz, Max},
  journal={Advances in Neural Information Processing Systems},
  volume={34},
  pages={29274--29286},
  year={2021}
}

@article{li2019convergence,
  title={On the convergence of fedavg on non-iid data},
  author={Li, Xiang and Huang, Kaixuan and Yang, Wenhao and Wang, Shusen and Zhang, Zhihua},
  journal={arXiv preprint arXiv:1907.02189},
  year={2019}
}

@inproceedings{mcmahan2017communication,
  title={Communication-efficient learning of deep networks from decentralized data},
  author={McMahan, Brendan and Moore, Eider and Ramage, Daniel and Hampson, Seth and y Arcas, Blaise Aguera},
  booktitle={Artificial intelligence and statistics},
  pages={1273--1282},
  year={2017},
  organization={PMLR}
}

@inproceedings{karimireddy2020scaffold,
  title={Scaffold: Stochastic controlled averaging for federated learning},
  author={Karimireddy, Sai Praneeth and Kale, Satyen and Mohri, Mehryar and Reddi, Sashank and Stich, Sebastian and Suresh, Ananda Theertha},
  booktitle={International Conference on Machine Learning},
  pages={5132--5143},
  year={2020},
  organization={PMLR}
}

@article{wang2020tackling,
  title={Tackling the objective inconsistency problem in heterogeneous federated optimization},
  author={Wang, Jianyu and Liu, Qinghua and Liang, Hao and Joshi, Gauri and Poor, H Vincent},
  journal={Advances in neural information processing systems},
  volume={33},
  pages={7611--7623},
  year={2020}
}

@article{li2020federated,
  title={Federated optimization in heterogeneous networks},
  author={Li, Tian and Sahu, Anit Kumar and Zaheer, Manzil and Sanjabi, Maziar and Talwalkar, Ameet and Smith, Virginia},
  journal={Proceedings of Machine learning and systems},
  volume={2},
  pages={429--450},
  year={2020}
}

@inproceedings{malik2019derivative,
  title={Derivative-free methods for policy optimization: Guarantees for linear quadratic systems},
  author={Malik, Dhruv and Pananjady, Ashwin and Bhatia, Kush and Khamaru, Koulik and Bartlett, Peter and Wainwright, Martin},
  booktitle={The 22nd international conference on artificial intelligence and statistics},
  pages={2916--2925},
  year={2019},
  organization={PMLR}
}

@article{mnih2015human,
  title={Human-level control through deep reinforcement learning},
  author={Mnih, Volodymyr and Kavukcuoglu, Koray and Silver, David and Rusu, Andrei A and Veness, Joel and Bellemare, Marc G and Graves, Alex and Riedmiller, Martin and Fidjeland, Andreas K and Ostrovski, Georg and others},
  journal={nature},
  volume={518},
  number={7540},
  pages={529--533},
  year={2015},
  publisher={Nature Publishing Group}
}

@article{rajeswaran2017learning,
  title={Learning complex dexterous manipulation with deep reinforcement learning and demonstrations},
  author={Rajeswaran, Aravind and Kumar, Vikash and Gupta, Abhishek and Vezzani, Giulia and Schulman, John and Todorov, Emanuel and Levine, Sergey},
  journal={arXiv preprint arXiv:1709.10087},
  year={2017}
}

@article{levine2016end,
  title={End-to-end training of deep visuomotor policies},
  author={Levine, Sergey and Finn, Chelsea and Darrell, Trevor and Abbeel, Pieter},
  journal={The Journal of Machine Learning Research},
  volume={17},
  number={1},
  pages={1334--1373},
  year={2016},
  publisher={JMLR. org}
}

@inproceedings{tobin2017domain,
  title={Domain randomization for transferring deep neural networks from simulation to the real world},
  author={Tobin, Josh and Fong, Rachel and Ray, Alex and Schneider, Jonas and Zaremba, Wojciech and Abbeel, Pieter},
  booktitle={2017 IEEE/RSJ international conference on intelligent robots and systems (IROS)},
  pages={23--30},
  year={2017},
  organization={IEEE}
}

@article{dulac2019challenges,
  title={Challenges of real-world reinforcement learning},
  author={Dulac-Arnold, Gabriel and Mankowitz, Daniel and Hester, Todd},
  journal={arXiv preprint arXiv:1904.12901},
  year={2019}
}

@article{konevcny2016federated,
  title={Federated optimization: Distributed machine learning for on-device intelligence},
  author={Kone{\v{c}}n{\`y}, Jakub and McMahan, H Brendan and Ramage, Daniel and Richt{\'a}rik, Peter},
  journal={arXiv preprint arXiv:1610.02527},
  year={2016}
}

@incollection{liang2022federated,
  title={Federated transfer reinforcement learning for autonomous driving},
  author={Liang, Xinle and Liu, Yang and Chen, Tianjian and Liu, Ming and Yang, Qiang},
  booktitle={Federated and Transfer Learning},
  pages={357--371},
  year={2022},
  publisher={Springer}
}

@article{lim2020federated,
  title={Federated reinforcement learning for training control policies on multiple IoT devices},
  author={Lim, Hyun-Kyo and Kim, Ju-Bong and Heo, Joo-Seong and Han, Youn-Hee},
  journal={Sensors},
  volume={20},
  number={5},
  pages={1359},
  year={2020},
  publisher={MDPI}
}

@article{yu2020deep,
  title={When deep reinforcement learning meets federated learning: Intelligent multitimescale resource management for multiaccess edge computing in 5G ultradense network},
  author={Yu, Shuai and Chen, Xu and Zhou, Zhi and Gong, Xiaowen and Wu, Di},
  journal={IEEE Internet of Things Journal},
  volume={8},
  number={4},
  pages={2238--2251},
  year={2020},
  publisher={IEEE}
}

@article{qi2021federated,
  title={Federated reinforcement learning: Techniques, applications, and open challenges},
  author={Qi, Jiaju and Zhou, Qihao and Lei, Lei and Zheng, Kan},
  journal={arXiv preprint arXiv:2108.11887},
  year={2021}
}

@article{zhao2022sample,
  title={On the sample complexity of stabilizing linear systems via policy gradient methods},
  author={Zhao, Feiran and Fu, Xingyun and You, Keyou},
  journal={arXiv preprint arXiv:2205.14335},
  year={2022}
}

@article{mohammadi2021convergence,
  title={Convergence and sample complexity of gradient methods for the model-free linear--quadratic regulator problem},
  author={Mohammadi, Hesameddin and Zare, Armin and Soltanolkotabi, Mahdi and Jovanovi{\'c}, Mihailo R},
  journal={IEEE Transactions on Automatic Control},
  volume={67},
  number={5},
  pages={2435--2450},
  year={2021},
  publisher={IEEE}
}

@article{hu2022towards,
  title={Towards a theoretical foundation of policy optimization for learning control policies},
  author={Hu, Bin and Zhang, Kaiqing and Li, Na and Mesbahi, Mehran and Fazel, Maryam and Ba{\c{s}}ar, Tamer},
  journal={arXiv preprint arXiv:2210.04810},
  year={2022}
}

@inproceedings{khaled2020tighter,
  title={Tighter theory for local SGD on identical and heterogeneous data},
  author={Khaled, Ahmed and Mishchenko, Konstantin and Richt{\'a}rik, Peter},
  booktitle={International Conference on Artificial Intelligence and Statistics},
  pages={4519--4529},
  year={2020},
  organization={PMLR}
}

@article{reddi2020adaptive,
  title={Adaptive federated optimization},
  author={Reddi, Sashank and Charles, Zachary and Zaheer, Manzil and Garrett, Zachary and Rush, Keith and Kone{\v{c}}n{\`y}, Jakub and Kumar, Sanjiv and McMahan, H Brendan},
  journal={arXiv preprint arXiv:2003.00295},
  year={2020}
}

@article{mitra2021linear,
  title={Linear convergence in federated learning: Tackling client heterogeneity and sparse gradients},
  author={Mitra, Aritra and Jaafar, Rayana and Pappas, George J and Hassani, Hamed},
  journal={Advances in Neural Information Processing Systems},
  volume={34},
  pages={14606--14619},
  year={2021}
}

@book{conn2009introduction,
  title={Introduction to derivative-free optimization},
  author={Conn, Andrew R and Scheinberg, Katya and Vicente, Luis N},
  year={2009},
  publisher={SIAM}
}

@article{nesterov2017random,
  title={Random gradient-free minimization of convex functions},
  author={Nesterov, Yurii and Spokoiny, Vladimir},
  journal={Foundations of Computational Mathematics},
  volume={17},
  pages={527--566},
  year={2017},
  publisher={Springer}
}

@article{hambly2021policy,
  title={Policy gradient methods for the noisy linear quadratic regulator over a finite horizon},
  author={Hambly, Ben and Xu, Renyuan and Yang, Huining},
  journal={SIAM Journal on Control and Optimization},
  volume={59},
  number={5},
  pages={3359--3391},
  year={2021},
  publisher={SIAM}
}

@article{jin2020analysis,
  title={On the analysis of model-free methods for the linear quadratic regulator},
  author={Jin, Zeyu and Schmitt, Johann Michael and Wen, Zaiwen},
  journal={arXiv preprint arXiv:2007.03861},
  year={2020}
}

@article{ju2022model,
  title={A model-free first-order method for linear quadratic regulator with $\tilde{O}(1/\varepsilon)$ sampling complexity},
  author={Ju, Caleb and Kotsalis, Georgios and Lan, Guanghui},
  journal={arXiv preprint arXiv:2212.00084},
  year={2022}
}

@article{
wang2023federated,
title={Federated {TD} Learning with Linear Function Approximation under Environmental Heterogeneity},
author={Han Wang and Aritra Mitra and Hamed Hassani and George J. Pappas and James Anderson},
journal={Transactions on Machine Learning Research},
issn={2835-8856},
year={2024},
url={https://openreview.net/forum?id=hdQspgyFrk},
note={}
}

@inproceedings{charles2,
  title={{Convergence and Accuracy Trade-Offs in Federated Learning and Meta-Learning}},
  author={Charles, Zachary and Kone{\v{c}}n{\`y}, Jakub},
  booktitle={International Conference on Artificial Intelligence and Statistics},
  pages={2575--2583},
  year={2021},
  organization={PMLR}
}

@article{charles,
  title={{On the outsized importance of learning rates in local update methods}},
  author={Charles, Zachary and Kone{\v{c}}n{\`y}, Jakub},
  journal={arXiv preprint arXiv:2007.00878},
  year={2020}
}

@article{khaled1,
  title={{First analysis of local gd on heterogeneous data}},
  author={Khaled, Ahmed and Mishchenko, Konstantin and Richt{\'a}rik, Peter},
  journal={arXiv preprint arXiv:1909.04715},
  year={2019}
}

@inproceedings{gorbunov,
  title={{Local SGD: Unified theory and new efficient methods}},
  author={Gorbunov, Eduard and Hanzely, Filip and Richt{\'a}rik, Peter},
  booktitle={International Conference on Artificial Intelligence and Statistics},
  pages={3556--3564},
  year={2021},
  organization={PMLR}
}

@article{fedsplit,
  title={{FedSplit: An algorithmic framework for fast federated optimization}},
  author={Pathak, Reese and Wainwright, Martin J},
  journal={arXiv preprint arXiv:2005.05238},
  year={2020}
}

@article{FedNova,
  title={{Tackling the objective inconsistency problem in heterogeneous federated optimization}},
  author={Wang, Jianyu and Liu, Qinghua and Liang, Hao and Joshi, Gauri and Poor, H Vincent},
  journal={Advances in Neural Information Processing Systems},
  volume={33},
  year={2020}
}

@article{li,
  title={{On the convergence of fedavg on non-iid data}},
  author={Li, Xiang and Huang, Kaixuan and Yang, Wenhao and Wang, Shusen and Zhang, Zhihua},
  journal={arXiv preprint arXiv:1907.02189},
  year={2019}
}

@article{acar2021,
  title={Federated learning based on dynamic regularization},
  author={Acar, Durmus Alp Emre and Zhao, Yue and Navarro, Ramon Matas and Mattina, Matthew and Whatmough, Paul N and Saligrama, Venkatesh},
  journal={arXiv preprint arXiv:2111.04263},
  year={2021}
}

@article{proxskip,
  title={{ProxSkip: Yes! Local Gradient Steps Provably Lead to Communication Acceleration! Finally!}},
  author={Mishchenko, Konstantin and Malinovsky, Grigory and Stich, Sebastian and Richt{\'a}rik, Peter},
  journal={arXiv preprint arXiv:2202.09357},
  year={2022}
}

@inproceedings{quantile,
  title={A superquantile approach to federated learning with heterogeneous devices},
  author={Laguel, Yassine and Pillutla, Krishna and Malick, Jer{\^o}me and Harchaoui, Zaid},
  booktitle={2021 55th Annual Conference on Information Sciences and Systems (CISS)},
  pages={1--6},
  year={2021},
  organization={IEEE}}

@article{sutton1999policy,
  title={Policy gradient methods for reinforcement learning with function approximation},
  author={Sutton, Richard S and McAllester, David and Singh, Satinder and Mansour, Yishay},
  journal={Advances in neural information processing systems},
  volume={12},
  year={1999}
}

@article{liu2020improved,
  title={An improved analysis of (variance-reduced) policy gradient and natural policy gradient methods},
  author={Liu, Yanli and Zhang, Kaiqing and Basar, Tamer and Yin, Wotao},
  journal={Advances in Neural Information Processing Systems},
  volume={33},
  pages={7624--7636},
  year={2020}
}

@article{williams1992simple,
  title={Simple statistical gradient-following algorithms for connectionist reinforcement learning},
  author={Williams, Ronald J},
  journal={Reinforcement learning},
  pages={5--32},
  year={1992},
  publisher={Springer}
}

@article{wang2019neural,
  title={Neural policy gradient methods: Global optimality and rates of convergence},
  author={Wang, Lingxiao and Cai, Qi and Yang, Zhuoran and Wang, Zhaoran},
  journal={arXiv preprint arXiv:1909.01150},
  year={2019}
}

@article{collins2022fedavg,
  title={Fedavg with fine tuning: Local updates lead to representation learning},
  author={Collins, Liam and Hassani, Hamed and Mokhtari, Aryan and Shakkottai, Sanjay},
  journal={arXiv preprint arXiv:2205.13692},
  year={2022}
}

@article{LuZD96,
  title={Stabilization of uncertain linear systems: An LFT approach},
  author={Lu, W-M and Zhou, Kemin and Doyle, John C},
  journal={IEEE Transactions on Automatic Control},
  volume={41},
  number={1},
  pages={50--65},
  year={1996},
  publisher={IEEE}
}

@ARTICLE{DGKF89,
  author={Doyle, John C and Glover, Keith and Khargonekar, Pramod and Francis, Bruce},
  journal={IEEE Transactions on Automatic Control}, 
  title={State-space solutions to standard ${H}_2$ and ${H}_\infty$ control problems}, 
  year={1989},
  volume={34},
  number={8},
  pages={831-847},
  doi={10.1109/9.29425}}

@article{zhangsurv,
  title={Multi-agent reinforcement learning: A selective overview of theories and algorithms},
  author={Zhang, Kaiqing and Yang, Zhuoran and Ba{\c{s}}ar, Tamer},
  journal={Handbook of reinforcement learning and control},
  pages={321--384},
  year={2021},
  publisher={Springer}
}

@article{lin2021multi,
  title={Multi-agent reinforcement learning in stochastic networked systems},
  author={Lin, Yiheng and Qu, Guannan and Huang, Longbo and Wierman, Adam},
  journal={Advances in Neural Information Processing Systems},
  volume={34},
  pages={7825--7837},
  year={2021}
}

@article{tropp2011freedman,
  title={Freedman's inequality for matrix martingales},
  author={Tropp, Joel},
  year={2011},
  journal={Electron. Commun. Probab.}
}

@article{konevcny2016FL,
  title={Federated learning: Strategies for improving communication efficiency},
  author={Kone{\v{c}}n{\`y}, Jakub and McMahan, H Brendan and Yu, Felix X and Richt{\'a}rik, Peter and Suresh, Ananda Theertha and Bacon, Dave},
  journal={arXiv preprint arXiv:1610.05492},
  year={2016}
}

@article{bonawitz2019towards,
  title={Towards federated learning at scale: System design},
  author={Bonawitz, Keith and Eichner, Hubert and Grieskamp, Wolfgang and Huba, Dzmitry and Ingerman, Alex and Ivanov, Vladimir and Kiddon, Chloe and Kone{\v{c}}n{\`y}, Jakub and Mazzocchi, Stefano and McMahan, Brendan and others},
  journal={Proceedings of machine learning and systems},
  volume={1},
  pages={374--388},
  year={2019}
}

@inproceedings{li2019convergence2,
  title={On the convergence of stochastic gradient descent with adaptive stepsizes},
  author={Li, Xiaoyu and Orabona, Francesco},
  booktitle={The 22nd international conference on artificial intelligence and Statistics},
  pages={983--992},
  year={2019},
  organization={PMLR}
}

@article{khaled2019first,
  title={First analysis of local GD on heterogeneous data},
  author={Khaled, Ahmed and Mishchenko, Konstantin and Richt{\'a}rik, Peter},
  journal={arXiv preprint arXiv:1909.04715},
  year={2019}
}

@article{stich2018local,
  title={Local SGD converges fast and communicates little},
  author={Stich, Sebastian U},
  journal={arXiv preprint arXiv:1805.09767},
  year={2018}
}

@article{wang2021cooperative,
  title={Cooperative sgd: A unified framework for the design and analysis of local-update sgd algorithms},
  author={Wang, Jianyu and Joshi, Gauri},
  journal={The Journal of Machine Learning Research},
  volume={22},
  number={1},
  pages={9709--9758},
  year={2021},
  publisher={JMLRORG}
}

@article{spiridonoff2020local,
  title={Local sgd with a communication overhead depending only on the number of workers},
  author={Spiridonoff, Artin and Olshevsky, Alex and Paschalidis, Ioannis Ch},
  journal={arXiv preprint arXiv:2006.02582},
  year={2020}
}

@inproceedings{reisizadeh2020fedpaq,
  title={Fedpaq: A communication-efficient federated learning method with periodic averaging and quantization},
  author={Reisizadeh, Amirhossein and Mokhtari, Aryan and Hassani, Hamed and Jadbabaie, Ali and Pedarsani, Ramtin},
  booktitle={International Conference on Artificial Intelligence and Statistics},
  pages={2021--2031},
  year={2020},
  organization={PMLR}
}

@article{haddadpour2019local,
  title={Local sgd with periodic averaging: Tighter analysis and adaptive synchronization},
  author={Haddadpour, Farzin and Kamani, Mohammad Mahdi and Mahdavi, Mehrdad and Cadambe, Viveck},
  journal={Advances in Neural Information Processing Systems},
  volume={32},
  year={2019}
}

@article{haddadpour2019convergence,
  title={On the convergence of local descent methods in federated learning},
  author={Haddadpour, Farzin and Mahdavi, Mehrdad},
  journal={arXiv preprint arXiv:1910.14425},
  year={2019}
}

@article{charles2020outsized,
  title={On the outsized importance of learning rates in local update methods},
  author={Charles, Zachary and Kone{\v{c}}n{\`y}, Jakub},
  journal={arXiv preprint arXiv:2007.00878},
  year={2020}
}

@inproceedings{charles2021convergence,
  title={Convergence and accuracy trade-offs in federated learning and meta-learning},
  author={Charles, Zachary and Kone{\v{c}}n{\`y}, Jakub},
  booktitle={International Conference on Artificial Intelligence and Statistics},
  pages={2575--2583},
  year={2021},
  organization={PMLR}
}

@article{pathak2020fedsplit,
  title={FedSplit: An algorithmic framework for fast federated optimization},
  author={Pathak, Reese and Wainwright, Martin J},
  journal={Advances in neural information processing systems},
  volume={33},
  pages={7057--7066},
  year={2020}
}

@article{tran2021feddr,
  title={FedDR--Randomized Douglas-Rachford splitting algorithms for nonconvex federated composite optimization},
  author={Tran Dinh, Quoc and Pham, Nhan H and Phan, Dzung and Nguyen, Lam},
  journal={Advances in Neural Information Processing Systems},
  volume={34},
  pages={30326--30338},
  year={2021}
}

@inproceedings{wang2022fedadmm,
  title={{FedADMM}: A federated primal-dual algorithm allowing partial participation},
  author={Wang, Han and Marella, Siddartha and Anderson, James},
  booktitle={2022 IEEE 61st Conference on Decision and Control (CDC)},
  pages={287--294},
  year={2022},
  organization={IEEE}
}

@article{wang2022fedsysid,
  title={{FedSysID}: A Federated Approach to Sample-Efficient System Identification},
  author={Wang, Han and Toso, Leonardo F and Anderson, James},
  journal={arXiv preprint arXiv:2211.14393},
  year={2022}
}

@inproceedings{schulman2015trust,
  title={Trust region policy optimization},
  author={Schulman, John and Levine, Sergey and Abbeel, Pieter and Jordan, Michael and Moritz, Philipp},
  booktitle={International conference on machine learning},
  pages={1889--1897},
  year={2015},
  organization={PMLR}
}

@article{konda1999actor,
  title={Actor-critic algorithms},
  author={Konda, Vijay and Tsitsiklis, John},
  journal={Advances in neural information processing systems},
  volume={12},
  year={1999}
}

@article{schulman2017proximal,
  title={Proximal policy optimization algorithms},
  author={Schulman, John and Wolski, Filip and Dhariwal, Prafulla and Radford, Alec and Klimov, Oleg},
  journal={arXiv preprint arXiv:1707.06347},
  year={2017}
}

@inproceedings{lamperski2020computing,
  title={Computing stabilizing linear controllers via policy iteration},
  author={Lamperski, Andrew},
  booktitle={2020 59th IEEE Conference on Decision and Control (CDC)},
  pages={1902--1907},
  year={2020},
  organization={IEEE}
}

@inproceedings{gatsis2022federated,
  title={Federated Reinforcement Learning at the Edge: Exploring the Learning-Communication Tradeoff},
  author={Gatsis, Konstantinos},
  booktitle={2022 European Control Conference (ECC)},
  pages={1890--1895},
  year={2022},
  organization={IEEE}
}

@article{fabbro2023federated,
  title={Federated TD Learning over Finite-Rate Erasure Channels: Linear Speedup under Markovian Sampling},
  author={Fabbro, Nicol{\`o} Dal and Mitra, Aritra and Pappas, George J},
  journal={arXiv preprint arXiv:2305.08104},
  year={2023}
}

@article{dean2020sample,
  title={On the sample complexity of the linear quadratic regulator},
  author={Dean, Sarah and Mania, Horia and Matni, Nikolai and Recht, Benjamin and Tu, Stephen},
  journal={Foundations of Computational Mathematics},
  volume={20},
  number={4},
  pages={633--679},
  year={2020},
  publisher={Springer}
}

@inproceedings{bach2016highly,
  title={Highly-smooth zero-th order online optimization},
  author={Bach, Francis and Perchet, Vianney},
  booktitle={Conference on Learning Theory},
  pages={257--283},
  year={2016},
  organization={PMLR}
}

@article{duchi2015optimal,
  title={Optimal rates for zero-order convex optimization: The power of two function evaluations},
  author={Duchi, John C and Jordan, Michael I and Wainwright, Martin J and Wibisono, Andre},
  journal={IEEE Transactions on Information Theory},
  volume={61},
  number={5},
  pages={2788--2806},
  year={2015},
  publisher={IEEE}
}

@article{polyak1987introduction,
  title={Introduction to optimization. optimization software},
  author={Polyak, Boris T},
  journal={Inc., Publications Division, New York},
  volume={1},
  pages={32},
  year={1987}
}

@inproceedings{
zhang2024finite,
title={Finite-Time Analysis of On-Policy Heterogeneous Federated Reinforcement Learning},
author={Chenyu Zhang and Han Wang and Aritra Mitra and James Anderson},
booktitle={The Twelfth International Conference on Learning Representations},
year={2024},
url={https://openreview.net/forum?id=D2eOVqPX9g}
}

@article{toso2024meta,
  title={Meta-Learning Linear Quadratic Regulators: A Policy Gradient MAML Approach for the Model-free {LQR}},
  author={Toso, Leonardo F and Zhan, Donglin and Anderson, James and Wang, Han},
  journal={arXiv preprint arXiv:2401.14534},
  year={2024}
}

@article{toso2025learning,
  title={Learning Stabilizing Policies via an Unstable Subspace Representation},
  author={Toso, Leonardo F and Ye, Lintao and Anderson, James},
  journal={arXiv preprint arXiv:2505.01348},
  year={2025}
}

@article{ziegler2019fine,
  title={Fine-tuning language models from human preferences},
  author={Ziegler, Daniel M and Stiennon, Nisan and Wu, Jeffrey and Brown, Tom B and Radford, Alec and Amodei, Dario and Christiano, Paul and Irving, Geoffrey},
  journal={arXiv preprint arXiv:1909.08593},
  year={2019}
}

@article{xie,
  author={Xie, Zhijie and Song, Shenghui},
  journal={IEEE Journal on Selected Areas in Communications}, 
  title={FedKL: Tackling Data Heterogeneity in Federated Reinforcement Learning by Penalizing KL Divergence}, 
  year={2023},
  volume={41},
  number={4},
  pages={1227-1242}}

@article{lan2023improved,
  title={Improved communication efficiency in federated natural policy gradient via {ADMM}-based gradient updates},
  author={Lan, Guangchen and Wang, Han and Anderson, James and Brinton, Christopher and Aggarwal, Vaneet},
  journal={arXiv preprint arXiv:2310.19807},
  year={2023}
}

@inproceedings{wang2024momentum,
  title={Momentum for the Win: Collaborative Federated Reinforcement Learning across Heterogeneous Environments},
  author={Wang, Han and He, Sihong and Zhang, Zhili and Miao, Fei and Anderson, James},
  booktitle={International Conference on Machine Learning},
  year={2024},
 organization={PMLR}
}

@inproceedings{zhu2024towards,
  title={Towards fast rates for federated and multi-task reinforcement learning},
  author={Zhu, Feng and Heath, Robert W and Mitra, Aritra},
  booktitle={2024 IEEE 63rd Conference on Decision and Control (CDC)},
  pages={2658--2663},
  year={2024},
  organization={IEEE}
}

@article{labbi2025global,
  title={On Global Convergence Rates for Federated Policy Gradient under Heterogeneous Environment},
  author={Labbi, Safwan and Mangold, Paul and Tiapkin, Daniil and Moulines, Eric},
  journal={arXiv preprint arXiv:2505.23459},
  year={2025}
}

@inproceedings{khodadadian,
  title={Federated Reinforcement Learning: Linear Speedup Under Markovian Sampling},
  author={Khodadadian, Sajad and Sharma, Pranay and Joshi, Gauri and Maguluri, Siva Theja},
  booktitle={International Conference on Machine Learning},
  pages={10997--11057},
  year={2022},
  organization={PMLR}
}

@article{liuMARL,
  title={Distributed {TD} $(0)$ with almost no communication},
  author={Liu, Rui and Olshevsky, Alex},
  journal={IEEE Control Systems Letters},
  volume={7},
  pages={2892--2897},
  year={2023},
  publisher={IEEE}
}

@article{tian2024one,
  title={One-Shot Averaging for Distributed {TD} ($\lambda$) Under Markov Sampling},
  author={Tian, Haoxing and Paschalidis, Ioannis Ch and Olshevsky, Alex},
  journal={IEEE Control Systems Letters},
  year={2024},
  publisher={IEEE}
}

@inproceedings{woo2023blessing,
  title={The blessing of heterogeneity in federated {Q}-learning: Linear speedup and beyond},
  author={Woo, Jiin and Joshi, Gauri and Chi, Yuejie},
  booktitle={International Conference on Machine Learning},
  pages={37157--37216},
  year={2023},
  organization={PMLR}
}

@inproceedings{salgia2024sample,
  title={The sample-communication complexity trade-off in federated q-learning},
  author={Salgia, Sudeep and Chi, Yuejie},
  booktitle={Advances in Neural Information Processing Systems},
  year={2024}
}

@article{dal2025finite,
  title={Finite-Time Analysis of Over-the-Air Federated TD Learning},
  author={Dal Fabbro, Nicol{\`o} and Mitra, Aritra and Heath, Robert W and Schenato, Luca and Pappas, George J},
  journal={IEEE Transactions on Wireless Communications},
  year={2025},
  publisher={IEEE}
}

@inproceedings{jinFRL,
  title={Federated Reinforcement Learning with Environment Heterogeneity},
  author={Jin, Hao and Peng, Yang and Yang, Wenhao and Wang, Shusen and Zhang, Zhihua},
  booktitle={International Conference on Artificial Intelligence and Statistics},
  pages={18--37},
  year={2022},
  organization={PMLR}
}

@article{shen2023,
  title={Towards Understanding Asynchronous Advantage Actor-critic: Convergence and Linear Speedup},
  author={Shen, Han and Zhang, Kaiqing and Hong, Mingyi and Chen, Tianyi},
  journal={IEEE Transactions on Signal Processing},
  year={2023},
  publisher={IEEE}
}

@article{mangold2024,
  title={Scafflsa: Taming heterogeneity in federated linear stochastic approximation and td learning},
  author={Mangold, Paul and Samsonov, Sergey and Labbi, Safwan and Levin, Ilya and Alami, Reda and Naumov, Alexey and Moulines, Eric},
  journal={Advances in Neural Information Processing Systems},
  year={2024}
}

@article{pan2025model,
  title={Model-agnostic meta-policy optimization via zeroth-order estimation: A linear quadratic regulator perspective},
  author={Pan, Yunian and Li, Tao and Zhu, Quanyan},
  journal={arXiv preprint arXiv:2503.00385},
  year={2025}
}

@inproceedings{aravind2024moreau,
  title={A Moreau envelope approach for LQR meta-policy estimation},
  author={Aravind, Ashwin and Toghani, Mohammad Taha and Uribe, C{\'e}sar A},
  booktitle={2024 IEEE 63rd Conference on Decision and Control (CDC)},
  pages={415--420},
  year={2024},
  organization={IEEE}
}

@inproceedings{toso2024asynchronous,
  title={Asynchronous heterogeneous linear quadratic regulator design},
  author={Toso, Leonardo F and Wang, Han and Anderson, James},
  booktitle={2024 IEEE 63rd Conference on Decision and Control (CDC)},
  pages={801--808},
  year={2024},
  organization={IEEE}
}

@article{zhao2025asynchronous,
  title={Asynchronous Parallel Policy Gradient Methods for the Linear Quadratic Regulator},
  author={Zhao, Feiran and Sha, Xingyu and You, Keyou},
  journal={IEEE Transactions on Automatic Control},
  year={2025},
  publisher={IEEE}
}

@inproceedings{ye2024convergence,
  title={On the Convergence of Policy Gradient for Designing a Linear Quadratic Regulator by Leveraging a Proxy System},
  author={Ye, Lintao and Mitra, Aritra and Gupta, Vijay},
  booktitle={2024 IEEE 63rd Conference on Decision and Control (CDC)},
  pages={6016--6021},
  year={2024},
  organization={IEEE}
}

@article{fujinami2025policy,
  title={Policy Gradient for LQR with Domain Randomization},
  author={Fujinami, Tesshu and Lee, Bruce D and Matni, Nikolai and Pappas, George J},
  journal={arXiv preprint arXiv:2503.24371},
  year={2025}
}

@inproceedings{toso2024oracle,
  title={Oracle complexity reduction for model-free lqr: A stochastic variance-reduced policy gradient approach},
  author={Toso, Leonardo F and Wang, Han and Anderson, James},
  booktitle={2024 American Control Conference (ACC)},
  pages={4032--4037},
  year={2024},
  organization={IEEE}
}

@article{zhang2025gap,
  title={{Gap-Dependent Bounds for Federated $ Q $-learning}},
  author={Zhang, Haochen and Zheng, Zhong and Xue, Lingzhou},
  journal={arXiv preprint arXiv:2502.02859},
  year={2025}
}

@article{labbi2024federated,
  title={{Federated ucbvi: Communication-efficient federated regret minimization with heterogeneous agents}},
  author={Labbi, Safwan and Tiapkin, Daniil and Mancini, Lorenzo and Mangold, Paul and Moulines, Eric},
  journal={arXiv preprint arXiv:2410.22908},
  year={2024}
}

@article{stamouli2025policy,
  title={{Policy Gradient Bounds in Multitask LQR}},
  author={Stamouli, Charis and Toso, Leonardo F and Tsiamis, Anastasios and Pappas, George J and Anderson, James},
  journal={arXiv preprint arXiv:2509.19266},
  year={2025}
}
\bibliographystyle{tmlr}

\newpage
\appendix
\section{Appendix}

\subsection{Appendix Roadmap}
This appendix is organized as follows. Sections \ref{sec:auxiliary_results_1} and \ref{sec:auxiliary_results_2} present auxiliary norm inequalities and lemmas that play a key role in proving the main results of this paper. The proof of our main results in the model-based setting is provided in Section \ref{sec:proof_model_based}, while Section \ref{sec:proof_model_free} is dedicated to the corresponding results in the model-free setting. Additional details on the zeroth-order optimization method are provided in Section \ref{sec:zeroth_order_optimization}.

\subsubsection{Notation Recap}

For convenience we briefly recap and summarize our notation. We use $\|S\|_{max}$ to denote the maximum spectral norm taken over the family of matrices $S^{(1)},\hdots, S^{(M)}$.  All norms for matrices and vectors are spectral and Euclidean respectively, unless otherwise stated. The integer sequence $\{1,2,\hdots,N\}$ is denoted as $[N]$. The spectral radius of a square matrix is denoted by $\rho(\cdot)$.

\begin{table}[tbhp]
\caption{Notations}
\centering
\begin{tabular}{|c|cl|} 
 \hline
 \textbf{Symbol} & & \textbf{Meaning}\\ \hline
 &&\\
$M$ & & number of systems \\
 &&\\
$L$ & & number of local updates (counter: $l$) \\
 &&\\
$N$ & & number of rounds of averaging (counter: $n$) \\
 &&\\
$K_n$ & & averaged controller at round $n$ \\
 &&\\
$K_i^*$ & & optimal controller for system $(A^{(i)},B^{(i)})$ \\
 &&\\
$K_{n,l}^{(i)}$ & & controller for system $i$ after $l$ local iterations and $n$ averaging rounds \\ \hline
% \caption{}
 \end{tabular}
 \end{table}

\subsection{Useful Norm Inequalities} \label{sec:auxiliary_results_1}
\begin{itemize}
 \item  Given any two matrices $A, B$ of the same dimensions,  for any $\xi>0$, we have
\begin{equation}\label{eq:youngs}
\|A+B\|_F^2 \leq(1+\xi)\|A\|_F^2+\left(1+\frac{1}{\xi}\right)\|B\|_F^2 .
\end{equation}
\item  Given any two matrices $A, B$ of the same dimensions,  for any $\xi>0$, we have
\begin{equation}\label{eq:youngs_inner}
\langle A, B\rangle \le \frac{\xi}{2}\lVert A\rVert_F^2 +\frac{1}{2\xi}\lVert B \rVert_F^2.
\end{equation}
This inequality goes by the name of Young's inequality.
 \item  Given $m$ matrices $A_1, \ldots, A_m$ of the same dimensions, the following is a simple application of Jensen's inequality:
\begin{align}\label{eq:sum_expand}
\left\|\sum_{i=1}^m A_i\right\|^2 &\leq m \sum_{i=1}^m\left\|A_i\right\|^2, \notag \\
\left\|\sum_{i=1}^m A_i\right\|_F^2 &\leq m \sum_{i=1}^m\left\|A_i\right\|_F^2.
\end{align}
\item  Given any two vectors $x, y \in \mathbb{R}^d$,  for any constant $\zeta>0$, we have
\begin{equation}\label{eq:vector_youngs}
\|x+y\|^2 \leq(1+\zeta)\|x\|^2+\left(1+\frac{1}{\zeta}\right)\|y\|^2 .
\end{equation}
\item Given any two vectors $x, y \in \mathbb{R}^d$,  for any constant $\zeta>0$, we have
\begin{equation}\label{eq:vector_youngs_inner}
\langle x, y\rangle \le \frac{\zeta}{2}\lVert x\rVert^2 +\frac{1}{2\zeta}\lVert y \rVert^2.
\end{equation}
\end{itemize}
\subsection{Useful Lemmas and Constants}\label{sec:auxiliary_results_2}

\begin{lemma}For each $i \in [M],$  we have that: 
\begin{align}
    ||\Sigma^{(i)}_K|| \leq \frac{{C}^{(i)}(K)}{\sigma_{\text{min}}(Q)}, 
   \quad ||P^{(i)}_K||\leq \frac{C^{(i)} (K)}{\mu}. \label{eq:upper_bound_sigma_bar_CK}
\end{align}
\end{lemma}
\begin{proof} The proof of this lemma is detailed in \citep[Lemma 13]{fazel2018global}. 
\end{proof}

\begin{lemma}\label{lemma:uniform_bounds} (Uniform bounds for $\nabla C^{(i)}(K)$ and $||K||$) For each agent $i\in [M]$, the gradient $\nabla C^{(i)}(K)$ and  $\|K\| $ can be bounded as follows:
$$
\|\nabla C^{(i)}(K)\| \leq\|\nabla C^{(i)}(K)\|_F \leq h_1(K) \quad \text {and} \quad\|K\| \leq h_2(K),
$$
where $h_1(K)$, and $h_2(K)$ are some positive scalars depending on the function $C(K)$.
\end{lemma}
\begin{proof} In this lemma, $h_1(K)$, and $h_2(K)$ are the functions defined as:
$$
\begin{aligned}
& h_0(K):=\sqrt{\frac{\left\|R_K\right\|_{\max}\left(C_{\max}(K)-C_{\min}\left(K\right)\right)}{\mu}}, \\
& h_1(K):= \frac{C_{\max}(K) h_0(K)}{\sigma_{\min}(Q)}, \quad h_2(K):=\frac{h_0(K)+\left\|B^{\top} P_K A\right\|_{\max}}{\sigma_{\min}(R)},
\end{aligned}
$$
where $\|R_K\|_{\max}:=\max_{i}\left\|R+B^{(i)\top} P^{(i)}_K B^{(i)}\right\|.$ By using \citep[Lemma 13]{fazel2018global}, we have 
\begin{align*}
\|\nabla C^{(i)}(K)\|^2 &\leq \operatorname{Tr}\left(\Sigma^{(i)}_K E^{(i)\top}_K E^{(i)}_K \Sigma^{(i)}_K\right) \leq\left\|\Sigma^{(i)}_K\right\|^2 \operatorname{Tr}\left(E^{(i)\top}_K E^{(i)}_K\right) \notag \\
& \leq\left(\frac{C^{(i)}(K)}{\sigma_{\min }(Q)}\right)^2 \operatorname{Tr}\left(E^{(i)\top}_K E^{(i)}_K\right).
\end{align*}
By \citep[Lemma 11]{fazel2018global}, we obtain 
$$
\operatorname{Tr}\left(E^{(i)\top}_K E^{(i)}_K\right) \leq \frac{\left\|R+B^{(i)\top} P^{(i)}_K B^{(i)}\right\|\left(C^{(i)}(K)-C^{(i)}\left(K_i^*\right)\right)}{\mu},
$$
which proves the first claim: 
\begin{align*}
\|\nabla C^{(i)}(K)\| & \leq \frac{C^{(i)}(K)}{\sigma_{\min}(Q)}\sqrt{\frac{\left\|R+B^{(i)\top} P^{(i)}_K B^{(i)}\right\|\left(C^{(i)}(K)-C^{(i)}\left(K_i^*\right)\right)}{\mu}}\notag \\
&\leq  \frac{C_{\max}(K)}{\sigma_{\min}(Q)}\sqrt{\frac{\left\|R+B^{(i)\top} P^{(i)}_K B^{(i)}\right\|_{\max}\left(C_{\max}(K)-C_{\min}\left(K\right)\right)}{\mu}}.
\end{align*}

On the other hand, by exploiting \citep[Lemma 11]{fazel2018global} we can also write
\begin{align*}
\|K\| & \leq\left\|\left(R+B^{(i)\top} P^{(i)}_K B^{(i)}\right)^{-1}\right\|\left\|\left(R+B^{(i)\top} P^{(i)}_K B^{(i)}\right) K\right\| \notag \\
& \leq \frac{1}{\sigma_{\min }(R)}\left\|\left(R+B^{(i)\top} P^{(i)}_K B^{(i)}\right) K\right\| \notag \\
& \leq \frac{1}{\sigma_{\min }(R)}\left(\left\|\left(R+B^{(i)\top} P^{(i)}_K B^{(i)}\right) K-B^{(i)\top} P^{(i)}_K A^{(i)}\right\|+\left\|B^{(i)\top} P^{(i)}_K A^{(i)}\right\|\right) \notag \\
& =\frac{\left\|E^{(i)}_K\right\|}{\sigma_{\min }(R)}+\frac{\left\|B^{(i)\top} P^{(i)}_K A^{(i)}\right\|}{\sigma_{\min }(R)} \notag \\
& \leq \frac{\sqrt{\operatorname{Tr}\left(E^{(i)\top}_K E^{(i)}_K\right)}}{\sigma_{\min }(R)}+\frac{\left\|B^{(i)\top} P^{(i)}_K A^{(i)}\right\|}{\sigma_{\min }(R)} \notag \\
& =\frac{\sqrt{\left(C^{(i)}(K)-C^{(i)}\left(K_i^*\right)\right)\left\|R+B^{(i)\top} P^{(i)}_K B^{(i)}\right\|}}{\sqrt{\mu} \sigma_{\min }(R)}+\frac{\left\|B^{(i)\top} P^{(i)}_K A^{(i)}\right\|}{\sigma_{\min }(R)},
\end{align*}
which completes the proof for the second claim. 
\end{proof}

%\vspace*{3em}
It is worth noting that the local cost and gradient smoothness in Lemma~\ref{Lemma:lipschitz} not only hold for the single-agent setting but also hold for the multi-agent setting. Moreover, we will make use of the following matrix concentration inequality.

\begin{definition}[Matrix Martingale] Consider a random process $\left\{{Y}_k: k=0,1,2, \ldots\right\}$ with matrices $Y_k$ of finite dimension, we say that the process is martingale if $Y_0 = 0$, $\mc{E}_{k-1} Y_k = Y_{k-1}$, and $\mc{E}\|Y_k\| < \infty,$ where $\mc{E}_{k-1}$ is the expectation conditioned on a filtration $\mathcal{F}_{k-1}$.

\end{definition}

\begin{lemma}\label{lem:Freedman}
(Rectangular Matrix Freedman~\citep{tropp2011freedman}). Consider a matrix martingale $\left\{{Y}_k: k=0,1,2, \ldots\right\}$ whose values are matrices with dimension $d_1 \times d_2$, and denote the difference sequence as $\left\{{X}_k = Y_{k} - Y_{k-1}: k=1,2,3, \ldots\right\}$. Assume that the difference sequence is uniformly bounded:
$$
\left\|{X}_k\right\| \leq R \quad \text { almost surely } \quad \text { for } k=1,2,3, \ldots .
$$
Define two predictable quadratic variation processes for this martingale:
$$
\begin{aligned}
& {W}_{\mathrm{col}, k}:=\sum_{j=1}^k \mathbb{E}_{j-1}\left({X}_j {X}_j^*\right) \text { and } \\
&{W}_{\text {row }, k}:=\sum_{j=1}^k \mathbb{E}_{j-1}\left({X}_j^* {X}_j\right) \quad \text { for } k=1,2,3, \ldots 
\end{aligned}
$$
Then, for all $t \geq 0$ and $\sigma^2>0$,
$$
\mathbb{P}\left\{\exists k \geq 0:\left\|{Y}_k\right\| \geq t \text { and } \max \left\{\left\|{W}_{\text {col }, k}\right\|,\left\|{W}_{\text {row }, k}\right\|\right\} \leq \sigma^2\right\} \leq\left(d_1+d_2\right) \cdot \exp \left\{-\frac{-t^2 / 2}{\sigma^2+R t / 3}\right\}.
$$
\end{lemma}

\subsubsection{Proof of Lemma~\ref{Lemma:lipschitz}}
\begin{proof}  In this proof, we aim to show
\begin{align*}
\left|C^{(i)}\left(K^{\prime}\right)-C^{(i)}(K)\right| &\leq h_{\text {cost}}(K) \|K^{\prime} -K\|, \\
\left\|\nabla C^{(i)}\left(K^{\prime}\right)-\nabla C^{(i)}(K)\right\| &\leq h_{\text {grad }}(K)\|\Delta\|, \quad \text{and } \\  \left\|\nabla C^{(i)}\left(K^{\prime}\right)-\nabla C^{(i)}(K)\right\|_F &\leq h_{\text {grad}}(K)\|\Delta\|_F,
\end{align*}
hold for all agents $i\in [M]$, and  $K^{\prime}$ satisfying $\|K^{\prime} -K\|  \leq h_{\Delta}(K)\ <\infty.$ The term
$h_{\Delta}(K)$ is the polynomial defined as
$$
h_{\Delta}(K):=\frac{\sigma_{\min}(Q) \mu}{4 ||B||_{\max} C_{\max}(K)\left(\left\|A-BK\right\|_{\max}+1\right)},
$$
the term $h_{\text {cost }}(K)$ and $h_{\text {grad}}(K)$ are defined as
\begin{align*}
&h_{\text {cost }}(K):=\frac{4 \operatorname{Tr}\left(\Sigma_0\right)C_{\max}(K)\|R\|}{\mu\sigma_{\min}\left(Q\right)}\left(\|K\|+\frac{h_{\Delta}(K)}{2}+\|B\|_{\max}\|K\|^2 \left(\left\|A -BK\right\|_{\max}+1\right) \frac{C_{\max}(K)}{\mu\sigma_{\min}(Q)}\right) 
\end{align*}
\begin{align*}
h_{\text {grad }}(K):=& 
 4\left(\frac{C_{\max}(K)}{\sigma_{\min}(Q)}\right)\Big[\|R\|+\|B\|_{\max}\left(\|A\|_{\max}+||B||_{\max}\left(\|K\|+h_{\Delta}(K)\right)\right)\nonumber\\
 &\times\left(\frac{h_{\text {cost }}(K) C_{\max}(K)}{\operatorname{Tr}\left(\Sigma_0\right)}\right)+ \|B\|^2_{\max}\frac{C_
{\max}(K)}{\mu}\Big]\\
&+8\left(\frac{C_{\max}(K)}{\sigma_{\min}(Q)}\right)^2\left(\frac{\|B\|_{\max}\left(\left\|A-BK\right\|_{\max}+1\right)}{\mu}\right) h_0(K).
\end{align*}
For the single-agent (i.e., $M=1$) setting, the proof is explained in detail in \citep[Lemma 24 and 25]{fazel2018global}. For the multi-agent setting (i.e., $M>1$), we can complete the proof by taking the maximum over the clients $i \in [M]$ of all the system-dependent parameters, such as $\lVert B\rVert_{\max}.$ 
\end{proof}

\subsubsection{Proof of Lemma~\ref{lemma:gradient_domination}}
\begin{proof} For the single-agent (i.e., $M=1$) setting, the proof is explained in~\citep[Lemma 11]{fazel2018global}. For the multi-agent setting (i.e., $M>1$), it is easy to see that
$$
C^{(i)}(K)-C^{(i)}\left(K_i^*\right) \leq \frac{\left\|\Sigma_{K_i^*}\right\|}{4\mu^2 \sigma_{\min }(R)}\|\nabla C^{(i)}(K)\|_F^2 \leq \frac{\max_i\left\| \Sigma_{K_i^*}\right\|}{4\mu^2 \sigma_{\min }(R)}\|\nabla C^{(i)}(K)\|_F^2
$$ holds for any stabilizing controller $K$ and any agent $i \in [M].$ 
\end{proof}

\subsection{The model-based setting} 

We first introduce the following operators on a symmetric matrix $X$,
\begin{align}\label{eq:linear_operators}
\mathcal{T}^{(i)}_K(X)&:=\sum_{t=0}^{\infty}(A^{(i)}-B^{(i)} K)^t X\left[(A^{(i)}-B^{(i)} K)^{\top}\right]^t, \notag \\
    \mathcal{F}^{(i)}_K(X)&:=(A^{(i)}-B^{(i)} K) X(A^{(i)}-B^{(i)} K)^{\top}.
\end{align}

We also define the induced norms of $\mathcal{T}$ and $\mathcal{F}$ as

\begin{align*}
    \left\|\mathcal{T}_K\right\|=\sup _X \frac{\left\|\mathcal{T}_K(X)\right\|}{\|X\|}, \quad \left\|\mathcal{F}_K\right\|=\sup _X \frac{\left\|\mathcal{F}_K(X)\right\|}{\|X\|}.
\end{align*}

\begin{lemma}
When $(A^{(i)}-B^{(i)} K)$ has spectral radius smaller than 1, we have 
$$
\mathcal{T}_K^{(i)}=\left(\mathrm{I}-\mathcal{F}_K^{(i)}\right)^{-1}
$$
holds for each $i \in [M].$
\end{lemma}
\begin{proof} The proof is detailed in \citep[Lemma 18]{fazel2018global}. 
\end{proof}

\begin{lemma} \label{lemma:bound_FK_TK}
If \footnote{This lemma has a similar flavor to that of Lemma 20 in~\citep{fazel2018global}. It is worthwhile to mention that the inequality~\eqref{eq:operator_requirement} imposes certain conditions on heterogeneity. Note that the constant $\frac{1}{2}$ can be changed into any finite constant. Thus, this heterogeneity requirement can be subsumed by that in Eq.\eqref{eq:het_requirement}.}
\begin{align}\label{eq:operator_requirement}
\left\|\mathcal{T}^{(i)}_K\right\|\left\|\mathcal{F}^{(i)}_K-\mathcal{F}^{(j)}_{K}\right\| \leq \frac{1}{2}
\end{align}
holds for any system $i, j \in [M],$
then we have
$$
\begin{aligned}
\left\|\left(\mathcal{T}^{(i)}_K-\mathcal{T}^{(j)}_{K}\right)(X)\right\| & \leq 2\left\|\mathcal{T}^{(i)}_K\right\|\left\|\mathcal{F}^{(i)}_K-\mathcal{F}^{(j)}_{K}\right\|\left\|\mathcal{T}^{(i)}_K(X)\right\| \\
& \leq 2\left\|\mathcal{T}^{(i)}_K\right\|^2\left\|\mathcal{F}^{(i)}_K-\mathcal{F}_{K}^{(j)}\right\|\|\ X \|.
\end{aligned}
$$
\end{lemma}
\begin{proof} Define $\mathcal{A}=\mathrm{I}-\mathcal{F}^{(i)}_K$, and $\mathcal{B}=\mathcal{F}^{(i)}_{K}-\mathcal{F}^{(j)}_K$. In this case $\mathcal{A}^{-1}=\mathcal{T}^{(i)}_K$ and $(\mathcal{A}-\mathcal{B})^{-1}=\mathcal{T}^{(j)}_{K}$. Hence, the condition $\left\|\mathcal{T}^{(i)}_K\right\|\left\|\mathcal{F}^{(i)}_K-\mathcal{F}^{(j)}_{K}\right\| \leq \frac{1}{2}$ translates to the condition $\left\|\mathcal{A}^{-1}\right\|\|\mathcal{B}\| \leq \frac{1}{2}$. 

First, we observe that
\begin{align}\label{eq:identical_transform}
\left(\mathcal{A}^{-1}-(\mathcal{A}-\mathcal{B})^{-1}\right)(X)&=\left(\mathrm{I}-\left(\mathrm{I}-\mathcal{A}^{-1} \circ \mathcal{B}\right)^{-1}\right)\left(\mathcal{A}^{-1}(X)\right)\notag\\
&=\left(\mathrm{I}-\left(\mathrm{I}-\mathcal{A}^{-1} \circ \mathcal{B}\right)^{-1}\right)\left(\mathcal{T}^{(i)}_K(X)\right),
\end{align}
where $ f\circ g$ denotes the composition $f(g(x))$.
Since $$\left(\mathrm{I}-\mathcal{A}^{-1} \circ \mathcal{B}\right)^{-1}=\mathrm{I}+\mathcal{A}^{-1} \circ \mathcal{B} \circ\left(\mathrm{I}-\mathcal{A}^{-1} \circ \mathcal{B}\right)^{-1},$$ we have:
\begin{equation}\label{eq:operator_inverse}
\left\|\left(\mathrm{I}-\mathcal{A}^{-1} \circ \mathcal{B}\right)^{-1}\right\| \leq 1+\left\|\mathcal{A}^{-1} \circ \mathcal{B}\right\|\left\|\left(\mathrm{I}-\mathcal{A}^{-1} \circ \mathcal{B}\right)^{-1}\right\| \leq 1+\frac{1}{2}\left\|\left(\mathrm{I}-\mathcal{A}^{-1} \circ \mathcal{B}\right)^{-1}\right\|
\end{equation}
Now rearranging terms in Eq.\eqref{eq:operator_inverse}, we obtain $\left\|\left(\mathrm{I}-\mathcal{A}^{-1} \circ \mathcal{B}\right)^{-1}\right\| \leq 2$.  Therefore, we have 
\begin{align*}
\left\|\mathrm{I}-\left(\mathrm{I}-\mathcal{A}^{-1} \circ \mathcal{B}\right)^{-1}\right\|=\left\|\mathcal{A}^{-1} \circ \mathcal{B} \circ\left(\mathrm{I}-\mathcal{A}^{-1} \circ \mathcal{B}\right)^{-1}\right\| &\leq\left\|\mathcal{A}^{-1}\right\|\|\mathcal{B}\|\left\|\left(\mathrm{I}-\mathcal{A}^{-1} \circ \mathcal{B}\right)^{-1}\right\|\\
&\leq 2\left\|\mathcal{A}^{-1}\right\|\|\mathcal{B}\|,
\end{align*}
and so
\begin{equation}\label{eq:operator_inequality}
\left\|\mathrm{I}-\left(\mathrm{I}-\mathcal{A}^{-1} \circ \mathcal{B}\right)^{-1}\right\| \leq 2\left\|\mathcal{A}^{-1}\right\|\|\mathcal{B}\|=2\left\|\mathcal{T}^{(i)}_K\right\|\left\|\mathcal{F}^{(i)}_K-\mathcal{F}^{(j)}_{K}\right\|.
\end{equation}

Then, we have 
\begin{align*}
\left\|\left(\mathcal{T}^{(i)}_K-\mathcal{T}^{(j)}_{K}\right)(X)\right\|& = \left\|\left(\mathcal{A}^{-1}-(\mathcal{A}-\mathcal{B})^{-1}\right)(X) \right\| \\
&\stackrel{(a)}{\leq}\left\|\left(\mathrm{I}-\left(\mathrm{I}-\mathcal{A}^{-1} \circ \mathcal{B}\right)^{-1}\right)\right\|\left\|\mathcal{T}^{(i)}_K(X)\right\| \notag\\
&\stackrel{(b)}{\leq} 2\left\|\mathcal{T}^{(i)}_K\right\|\left\|\mathcal{F}^{(i)}_K-\mathcal{F}^{(j)}_{K}\right\|\left\|\mathcal{T}^{(i)}_K(X)\right\| \notag\\
& \leq 2\left\|\mathcal{T}^{(i)}_K\right\|\left\|\mathcal{F}^{(i)}_K-\mathcal{F}^{(j)}_{K}\right\|\left\|\mathcal{T}^{(i)}\right\| \left\| X\right\|,
\end{align*}
where $(a)$ is due to Eq.\eqref{eq:identical_transform} and $(b)$ is due to Eq.\eqref{eq:operator_inequality}.
This completes the proof of Lemma~\ref{lemma:bound_FK_TK}. 
\end{proof}

\label{sec:proof_model_based}
\subsubsection{Proof of Lemma~\ref{lem:gradient_het}}
\begin{proof} First, we know that $\nabla C^{(i)}(K)$ and $\nabla C^{(j)}(K)$ are given by,
    \begin{align*}
        \nabla C^{(i)}(K)=2E^{(i)}_K\Sigma^{(i)}_{K}, \quad \text{and} \quad \nabla C^{(j)}(K)=2E^{(j)}_K{\Sigma}^{(j)}_{K}
        \end{align*}
where,
\begin{align*}
    E^{(i)}_K &= (R+B^{(i)\top} P^{(i)}_K B^{(i)})K - B^{(i)\top} P^{(i)}_K A^{(i)}, \nonumber 
  %  E^{(j)}_K &= (R+B^{(j)\top} P^{(j)}_K B^{(j)})K - B^{(j)\top} P^{(j)}_K A^{(j)},
\end{align*}
and
\begin{align*}
     \Sigma^{(i)}_K= _{x^{(i)}_0\sim D}\sum_{t=0}^\infty x^{(i)}_t x^{(i)\top}_t. %\;\ \Sigma^{(j)}_K= _{x^{(j)}_0\sim D}\sum_{t=0}^\infty x^{(j)}_t x^{(j)\top}_t.
\end{align*}

Thus, we can write,
\begin{align*}
    || \nabla C^{(i)}(K) - \nabla C^{(j)}(K) || &=  ||2E^{(i)}_K\Sigma^{(i)}_K - 2E^{(j)}_K\Sigma^{(j)}_K ||\\
    &\leq 2( ||E^{(i)}_K -E^{(j)}_K||\underbrace{||\Sigma^{(i)}_K||}_{\beta_1} + \underbrace{||E^{(j)}_K||}_{\beta_2}||\Sigma^{(i)}_K -\Sigma^{(j)}_K||).
\end{align*}

From Eq.~\eqref{eq:upper_bound_sigma_bar_CK} we can upper bound $||\Sigma^{(i)}_K||$ as follows:
\begin{align*}
    ||\Sigma^{(i)}_K|| \leq \frac{{C}^{(i)}(K)}{\sigma_{\text{min}}(Q)}.
    %\label{eq:upper_bound_sigma_bar_CK}
\end{align*}
With the definition of $E^{(j)}_K = RK + B^{(j)\top} P^{(j)}_K B^{(j)}K - B^{(j)\top} P^{(j)}_KA^{(j)}$, we can use triangle inequality to write,
\begin{align*}
    || E^{(j)}_K || &\leq  ||RK|| + ||B^{(j)}|| ||P^{(j)}_K|| ||B^{(j)}K|| + ||B^{(j)}||||P^{(j)}_K|| ||A^{(j)}|| \nonumber\\
    &\leq ||RK|| + \frac{||B^{(j)}|| C^{(j)}(K)}{{\mu}}(||B^{(j)}K|| + ||A^{(j)}||),
    %\label{eq:upper_bound_EK_bar}
\end{align*}
\noindent where $||P^{(j)}_K||\leq \frac{C^{(j)} (K)}{\mu}$ from Eq.~\eqref{eq:upper_bound_sigma_bar_CK}, with $\mu = \sigma_{\min}({\Sigma}^{(j)}_0)$. 

With the notation that we introduced previously, we can write 
\begin{align*}
    \beta_1 = ||\Sigma^{(i)}_K|| \leq ||\Sigma_K||_{\max} \leq \frac{{C_{\max}}(K)}{\sigma_{\min}(Q)},
\end{align*}
and, 
\begin{align*}
    \beta_2 = || E^{(j)}_K ||\leq ||E_K||_{\max} \leq ||R|| ||K|| + \frac{||B||_{\max} C_{\max}(K)}{{\mu}}(||B||_{\max} + ||A||_{\max}),
\end{align*}
where $C_{\max}(K) := \max_{i}C^{(i)}(K)$.

Next we will derive an upper bound for $||E^{(i)}_K - E^{(j)}_K||$.

\textbf{Upper bound for $||E^{(i)}_K - E^{(j)}_K||$:} We can first use the definition of  $E_K^{(i)}$ and $E^{(j)}_K$ to write, 

\begin{align*}
    E^{(i)}_K &- E^{(j)}_K = B^{(j)\top} P^{(j)}_K(A^{(j)} - B^{(j)}K) - B^{(i)\top} P^{(i)}_K(A^{(i)} - B^{(i)}K)\nonumber\\
    &=-B^{(i)\top} P^{(i)}_K(A^{(i)} - B^{(i)}K) + B^{(i)\top} P^{(i)}_K(A^{(j)}- B^{(j)}K) - B^{(i)\top} P^{(i)}_K(A^{(j)}- B^{(j)}K)\nonumber \\
    &+ B^{(i)\top} P^{(j)}_K(A^{(j)}- B^{(j)}K) - B^{(i)\top} P^{(j)}_K(A^{(j)}- B^{(j)}K) + B^{(j)\top} P^{(j)}_K(A^{(j)} - B^{(j)}K).
\end{align*}

Then, by using triangle inequality, we obtain the following expression:
\begin{align*}
    ||E^{(i)}_K - E^{(j)}_K || &\leq ||\underbrace{B^{(i)\top} P^{(i)}_K(A^{(i)} - B^{(i)}K)-B^{(i)\top} P^{(i)}_K(A^{(j)}- B^{(j)}K)}_{H_1}||\nonumber\\
    &+ ||\underbrace{B^{(i)\top} P^{(i)}_K(A^{(j)}- B^{(j)}K) - B^{(i)\top} P^{(j)}_K(A^{(j)}- B^{(j)}K)}_{H_2}||\nonumber\\
    &+||\underbrace{B^{(i)\top} P^{(j)}_K(A^{(j)}- B^{(j)}K) - B^{(j)\top}P^{(j)}_K(A^{(j)} - B^{(j)}K)}_{H_3}||.
\end{align*}
Incorporating the heterogeneity bounds from assumption~\ref{assumption:bnd_sys_heterogeneity} gives 
\begin{align*}
    ||H_1||\leq ||B^{(i)}||||P^{(i)}_K||(\epsilon_1 + \epsilon_2||K||),
\end{align*}
to which we apply the max-norm definition to arrive at
\begin{align}\label{eq:bound_H_1}
    ||H_1||\leq ||B||_{\max}(\epsilon_1 + \epsilon_2||K||) ||P_K||_{\max}.
\end{align}
Similarly, we can also derive upper bounds for $||H_2||$ and $||H_3||$, as follows, 
\begin{align}\label{eq:bound_H_2}
    ||H_2||\leq ||B^{(i)}|| ||P^{(i)}_K - P^{(j)}_K|| ||A^{(j)} - B^{(j)}K|| \leq ||B||_{\max} ||P^{(i)}_K - P^{(j)}_K|| ||A - BK||_{\max}
\end{align}
and 
\begin{align}\label{eq:bound_H_3}
    ||H_3||\leq \epsilon_2||A^{(i)} - B^{(i)}K||||P^{(j)}_K|| \leq \epsilon_2||A - BK||_{\max}||P_K||_{\max}.
\end{align}
To bound $H_2$, we need to derive an upper bound for $||P^{(i)}_K - P^{(j)}_K||$. For this purpose, we have that for any fixed system $i \in [M]$
\begin{align*}
 ||P^{(i)}_K - P^{(j)}_K||=\left\|\mathcal{T}^{(i)}_{K}\left(Q+K^{\top} R K\right)-\mathcal{T}^{(j)}_K\left(Q+K^{\top} R K\right)\right\|.
\end{align*}
Thus, by using Lemma~\ref{lemma:bound_FK_TK}, we can write, 
\begin{align*}
 ||P^{(i)}_K - P^{(j)}_K||\leq 2\left\|\mathcal{T}^{(i)}_K\right\|^2\left\|\mathcal{F}^{(i)}_K-\mathcal{F}_{K}^{(j)}\right\|\left\| Q + K^\top R K\right\|,\\
\end{align*}
where $||\mathcal{T}^{(i)}_K|| \leq \frac{C^{(i)}(K)}{\sigma_{\min}(Q)\mu} \leq \frac{C_{\max}(K)}{\sigma_{\min}(Q)\mu}$ (detailed in Lemma 17 of \citep{fazel2018global}). With the following upper bound for $\left\|\mathcal{F}^{(i)}_K-\mathcal{F}_{K}^{(j)}\right\|$: 
\begin{align*}
    ||(\mathcal{F}^{(i)}_K-\mathcal{F}_{K}^{(j)})(X)||&= ||(A^{(i)}-B^{(i)} K) X(A^{(i)}-B^{(i)} K)^{\top} \\ 
    &- (A^{(j)}-B^{(j)} K) X(A^{(j)}-B^{(j)} K)^{\top}||\nonumber \\
    &\leq 2(\epsilon_1 + \epsilon_2||K||)||X||||A - BK||_{\max},
\end{align*}
we have 
\begin{align}\label{eq:bound_diff_P}
 ||P^{(i)}_K - P^{(j)}_K||\leq 4\left(\frac{C_{\max}(K)}{\sigma_{\min}(Q)\mu}\right)^2 (\epsilon_1 + \epsilon_2||K||)||A - BK||_{\max}(||Q|| + ||R|| ||K||^2),
\end{align}

Plugging in Eq.~\eqref{eq:bound_diff_P} into $H_2$ and adding the upper bounds of $H_1$ (Eq.~\ref{eq:bound_H_1}),  $H_2$ (Eq.~\ref{eq:bound_H_2}) and $H_3$ (Eq.~\ref{eq:bound_H_3}) together, we have
%\begin{center}
%\fbox{\begin{minipage}{15em}
\begin{align*}
    ||E^{(i)}_K - E^{(j)}_K|| &\leq g_1(\epsilon_1,\epsilon_2,K),
\end{align*}
%\end{minipage}}
%\end{center}
where $g_1$ is a linear in $\epsilon_1, \epsilon_2$ and polynomial in the remaining problem data. Specifically,
\begin{align*}
&g_1(\epsilon_1,\epsilon_2,K) := \epsilon_1\left(\frac{||B||_{\max}C_{\max}(K)}{\mu}\left[ 1 + 4\left(\frac{C_{\max}(K)}{\sigma_{\min}(Q)\mu} \right)\left(||A-BK||_{\max}\right)^2  \left(||Q|| + ||R|| ||K||^2\right)  \right]  \right)\notag\\
&+ \epsilon_2\left(\frac{||B||_{\max}||K||C_{\max}(K)}{\mu}\left[ 1 + 4\left(\frac{C_{\max}(K)}{\sigma_{\min}(Q)\mu} \right)\left(||A-BK||_{\max}\right)^2\left(||Q|| + ||R|| ||K||^2\right)\right]\right. \\
&\left.+||A-BK||_{\max} \right).
\end{align*}  

In what follows, we will derive an upper bound for $||\Sigma^{(i)}_K -\Sigma^{(j)}_K||$:

\textbf{Upper bound for $||\Sigma^{(i)}_K -\Sigma^{(j)}_K||$:} From the previous definitions in Eq.\eqref{eq:linear_operators} and Lemma~\ref{lemma:bound_FK_TK}, we have, 
\begin{align*}
    ||\Sigma^{(i)}_K -\Sigma^{(j)}_K|| &= ||\mathcal{T}^{(i)}_K(\Sigma_0) - \mathcal{T}^{(j)}_K(\Sigma_0)||\leq 2\left\|\mathcal{T}^{(i)}_K\right\|^2\left\|\mathcal{F}^{(i)}_K-\mathcal{F}_{K}^{(j)}\right\|\|\ \Sigma_0 \| \nonumber \\
    &\leq 4\left(\frac{C_{\max}}{\sigma_{\min}(Q)\mu}\right )^2(\epsilon_1 + \epsilon_2||K||)||A-BK||_{\max}||\Sigma_0||
\end{align*}
where $\Sigma_0=\mathbb{E}_{x_0^{(i)} \sim \mathcal{D}}\left[x_0^{(i)} x_0^{(i)\top}\right].$

Thus, we have the following upper bound for $||\Sigma^{(i)}_K -\Sigma^{(j)}_K||$,

%\begin{center}
%\fbox{\begin{minipage}{15em}
\begin{align*}
     ||\Sigma^{(i)}_K -\Sigma^{(j)}_K|| \leq g_2(\epsilon_1,\epsilon_2,K) 
\end{align*}
%\end{minipage}}
%\end{center}
with,
\begin{align*}
    &g_2(\epsilon_1,\epsilon_2,K):=  \epsilon_1 \left(\frac{C_{\max}(K)}{\sigma_{\min}(Q)\mu}\right )^2\left(4||A-BK||_{\max}||\Sigma_0|| \right)+ \epsilon_2 ||K||\left(\frac{C_{\max}(K)}{\sigma_{\min}(Q)\mu}\right )^2\left(4||A-BK||_{\max}||\Sigma_0|| \right).
\end{align*}

Therefore, we can finally write an upper bound for $|| \nabla C^{(i)}(K) -  \nabla C^{(j)}(K)||$, which is:
\begin{align*}
    ||\nabla C^{(i)}(K) -  \nabla C^{(j)}(K)|| \leq f(\epsilon_1,\epsilon_2,K)
\end{align*}
where, 
\begin{align*}
    f(\epsilon_1,\epsilon_2,K)= 2( \beta_1g_1(\epsilon_1,\epsilon_2,K) + \beta_2g_2(\epsilon_1,\epsilon_2,K)).
\end{align*}
After some rearrangement, we have that
\begin{align*}
     f(\epsilon_1,\epsilon_2,K)= \epsilon_1 h^1_{\text{het}}(K) + \epsilon_1 h^2_{\text{het}}(K),
\end{align*}
where $h^1_{\text{het}} = h_{1f} + h_{2f}$ and $h^2_{\text{het}}= h_{3f} + h_{4f}$, and

 \begin{align*}
    &h_{1f}= \frac{2||B||_{\max}(C_{\max}(K))^2}{\sigma_{\min}(Q)\mu}\left[ 1 + 4\left(\frac{C_{\max}(K)}{\sigma_{\min}(Q)\mu} \right)\left(||A-BK||_{\max}\right)^2\left(||Q|| + ||R|| ||K||^2\right)\right],  \\
     &h_{2f}=\frac{2}{\mu}\left(\frac{C_{\max}(K)}{\sigma_{\min}(Q)}\right )^3\left(4||A-BK||_{\max}||\Sigma_0|| \right),
\end{align*}

\begin{align*}
     &h_{3f}= 2\left(||R|| ||K|| + \frac{||B||_{\max} C_{\max}(K)}{{\mu}}(||B||_{\max} + ||A||_{\max}) \right)\\
     &\times \left(\frac{||B||_{\max}||K||C_{\max}(K)}{\mu}\left[ 1 + 4\left(\frac{C_{\max}(K)}{\sigma_{\min}(Q)\mu} \right)\left(||A-BK||_{\max}\right)^2\left(||Q|| + ||R|| ||K||^2\right)\right] \right.\\
     &\left.+||A-BK||_{\max} \right),\\
      &h_{4f} = 8\left(||R|| ||K|| + \frac{||B||_{\max} C_{\max}(K)}{{\mu}}(||B||_{\max} + ||A||_{\max}\right) ||K||\left(\frac{C_{\max}(K)}{\sigma_{\min}(Q)\mu}\right )^2 ||A-BK||_{\max}||\Sigma_0||.
\end{align*}    
\end{proof}

\subsubsection{Proof of Theorem~\ref{thm:one_local_model_based}}\label{sec:proof_thm3}

In this theorem,  we consider the setting where $
(\epsilon_1 \Bar{h}^1_{\text{het}} + \epsilon_2 \Bar{h}^2_{\text{het}})^2 \le\Bar{h}^3_{\text{het}}$ with  
\begin{equation}\label{eq:het_requirement}
\Bar{h}^3_{\text{het}}:=\min_{j \in [M]} \left\{\frac{\mu^2 \sigma_{\min }(R)\left(C^{(j)}(K_0) - C^{(j)}(K_j^*)\right)}{4||\Sigma_{K_j^*}||\min\{n_x, n_u\}}\right\}.\end{equation}

\paragraph{Outline:} To prove Theorem~\ref{thm:one_local_model_based}, we first introduce some lemmas: Lemma \ref{lem:stability_local} establishes stability of the local policies; Lemma \ref{lem：drift_term} provides the drift analysis; Lemma \ref{lem:model_based_prp} quantifies the per-round progress of our \texttt{FedLQR} algorithm. As a result, we are able to present the iterative stability guarantees and convergence analysis of \texttt{FedLQR} in the model-based setting.

\begin{lemma} \label{lem:stability_local} (Stability of the local policies)  Suppose $K_n \in \mathcal{G}^0.$ If the local step-size satisfies $\eta_l \le \min \{\frac{\underline{h}_{\Delta}}{\bar{h}_1},\frac{1}{4\bar{h}_{\text{grad}}}\}$ and the heterogeneity level satisfies $(\epsilon_1 \Bar{h}^1_{\text{het}} + \epsilon_2 \Bar{h}^2_{\text{het}})^2 \le\Bar{h}^3_{\text{het}}$, then $K_{n,l}^{(i)} \in \mathcal{G}^0$ holds for all $i \in [M]$ and $l \in [L]$.
\end{lemma}
\begin{proof}
Since $K_n \in\mathcal{G}^0,$ based on the local Lipschitz property in Lemma~\ref{Lemma:lipschitz}, we have:
\begin{align}\label{eq:stability_decrease}
C^{(j)}(K_{n,1}^{(i)}) - C^{(j)}(K_n)&\le \left\langle \nabla C^{(j)}(K_n), K_{n,1}^{(i)} -K_n   \right\rangle + \frac{h_{\text{grad}}(K_n)}{2}\Big\lVert K_{n,1}^{(i)} -K_n \Big\rVert_F^2 \quad \notag\\
& \le -\left\langle \nabla C^{(j)}(K_n), \eta_l{\nabla} C^{(i)}(K_n) \right\rangle + \frac{h_{\text{grad}}(K_n)}{2}\Big\lVert \eta_l {\nabla} C^{(i)}(K_n) \Big\rVert_F^2
\end{align}
holds for any $i, j \in [M],$ if $\Big\lVert \eta_l {\nabla} C^{(i)}(K_n) \Big\rVert_F  \le \underline{h}_{\Delta} \le h_{\Delta}(K_n),$ which holds when 
\begin{align*}
\Big\lVert \eta_l {\nabla} C^{(i)}(K_n) \Big\rVert_F \stackrel{(a)}{\le}\eta_l h_1(K_n) \le \eta_l  \bar{h}_1\stackrel{(b)}{\le} \underline{h}_{\Delta},
\end{align*}
where $(a)$ comes from Lemma~\ref{lemma:uniform_bounds} and $(b)$ holds because of the requirement on $\eta_l$ in the statement of the lemma. We note that although $h_{\text{grad}}(K_n)$ depends on $K_n$, our analysis is always restricted to the common stabilizing sublevel set $\mathcal{G}^0(\beta)$, on which, $h_{\text{grad}}(K_n)$  is bounded from above by some $\bar{h}_{\text{grad}} < \infty$ (see $(\phi,\rho)$-local smoothness property in \cite{malik2019derivative}).

Following the analysis in Eq~\eqref{eq:stability_decrease}, we have
\begin{align}\label{eq:stability_decrease_II}
C^{(j)}(K_{n,1}^{(i)}) - C^{(j)}(K_n)&\le - \eta_l\left\langle \nabla C^{(j)}(K_n), \nabla C^{(j)}(K_n) \right\rangle \notag \\ 
&-\eta_l\underbrace{\left\langle \nabla C^{(j)}(K_n), \nabla C^{(i)}(K_n)-\nabla C^{(j)}(K_n)  \right\rangle}_{T_{1}} \notag\\
& + \frac{h_{\text{grad}}(K_n)}{2}\Big\lVert \eta_l{\nabla} C^{(i)}(K_n) \Big\rVert_F^2.
\end{align}
Now $T_{1}$ can be bounded as
\begin{align}
T_{1} &\le \eta_l \Big\lVert \nabla C^{(j)}(K_n)\Big\rVert_F \Big\lVert \nabla C^{(i)}(K_n)- \nabla C^{(j)}(K_n) \Big\rVert_F\notag\\
& \le \eta_l \sqrt{\min\{n_x, n_u\}} \Big\lVert \nabla C^{(j)}(K_n)\Big\rVert_F \Big\lVert \nabla C^{(i)}(K_n)- \nabla C^{(j)}(K_n) \Big\rVert \notag\\
&\stackrel{(c)}{\le}  \eta_l \sqrt{\min\{n_x, n_u\}} \Big\lVert \nabla C^{(j)}(K_n)\Big\rVert_F (\epsilon_1\bar{h}^{1}_{\text{het}} +\epsilon_2\Bar{h}^{2}_{\text{het}}),
\end{align}
where $(c)$ is due to Lemma~\ref{lem:gradient_het}.
Plugging in the upper bound of $T_1$ into ~\eqref{eq:stability_decrease}, we have:
\begin{align*}
&C^{(j)}(K_{n,1}^{(i)}) - C^{(j)}(K_n)\le - \eta_l\left\langle \nabla C^{(j)}(K_n), \nabla C^{(j)}(K_n) \right\rangle\notag\\
&+\eta_l \sqrt{\min\{n_x, n_u\}} \Big\lVert \nabla C^{(j)}(K_n)\Big\rVert_F (\epsilon_1\bar{h}^{1}_{\text{het}} +\epsilon_2\Bar{h}^{2}_{\text{het}}) + \frac{h_{\text{grad}}(K_n)}{2}\Big\lVert \eta_l{\nabla} C^{(i)}(K_n) \Big\rVert_F^2\notag\\
&\stackrel{(d)}{\le} - \eta_l\left\langle \nabla C^{(j)}(K_n), \nabla C^{(j)}(K_n) \right\rangle+\eta_l \sqrt{\min\{n_x, n_u\}} \Big\lVert \nabla C^{(j)}(K_n)\Big\rVert_F (\epsilon_1\bar{h}^{1}_{\text{het}} +\epsilon_2\Bar{h}^{2}_{\text{het}})  \notag\\
& + h_{\text{grad}}(K_n)\Big\lVert \eta_l{\nabla} C^{(j)}(K_n) \Big\rVert_F^2 +h_{\text{grad}}(K_n)\Big\lVert \eta_l{\nabla} C^{(i)}(K_n) - \eta_l{\nabla} C^{(j)}(K_n) \Big\rVert_F^2\notag\\
&\stackrel{(e)}{\le} - \eta_l\left\langle \nabla C^{(j)}(K_n), \nabla C^{(j)}(K_n) \right\rangle+\eta_l \sqrt{\min\{n_x, n_u\}} \Big\lVert \nabla C^{(j)}(K_n)\Big\rVert_F (\epsilon_1\bar{h}^{1}_{\text{het}} +\epsilon_2\Bar{h}^{2}_{\text{het}})  \notag\\
&   + \eta_l^2 h_{\text{grad}}(K_n)\Big\lVert \nabla C^{(j)}(K_n) \Big\rVert_F^2 +\eta_l^2 h_{\text{grad}}(K_n) \min\{n_x, n_u\}(\epsilon_1\bar{h}^{1}_{\text{het}} +\epsilon_2\Bar{h}^{2}_{\text{het}})^2\notag\\
& \stackrel{(f)}{\le} - \eta_l\left\langle \nabla C^{(j)}(K_n), \nabla C^{(j)}(K_n) \right\rangle+\frac{\eta_l}{4} \Big\lVert \nabla C^{(j)}(K_n) \Big\rVert_F^2+ \eta_l \min\{n_x, n_u\}(\epsilon_1\bar{h}^{1}_{\text{het}} +\epsilon_2\Bar{h}^{2}_{\text{het}})^2 \notag\\
& + \eta_l^2 \bar{h}_{\text{grad}}\Big\lVert \nabla C^{(j)}(K_n) \Big\rVert_F^2 +\eta_l^2 \bar{h}_{\text{grad}} \min\{n_x, n_u\}(\epsilon_1\bar{h}^{1}_{\text{het}} +\epsilon_2\Bar{h}^{2}_{\text{het}})^2\notag\\
& = - \eta_l\left\langle \nabla C^{(j)}(K_n), \nabla C^{(j)}(K_n) \right\rangle+(\frac{\eta_l}{4} +\eta_l^2 \bar{h}_{\text{grad}}) \Big\lVert \nabla C^{(j)}(K_n) \Big\rVert_F^2 \notag\\
& + (\eta_l + \eta_l^2 \bar{h}_{\text{grad}})\min\{n_x, n_u\}(\epsilon_1\bar{h}^{1}_{\text{het}} +\epsilon_2\Bar{h}^{2}_{\text{het}})^2\notag\\
&\stackrel{(g)}{\le} -\frac{\eta_l}{2} \Big\lVert \nabla C^{(j)}(K_n) \Big\rVert_F^2  + 2\eta_l \min\{n_x, n_u\}(\epsilon_1\bar{h}^{1}_{\text{het}} +\epsilon_2\Bar{h}^{2}_{\text{het}})^2\notag,
\end{align*}
which implies
\begin{align}\label{eq:stability_decrease_continue}
C^{(j)}(K_{n,1}^{(i)}) - C^{(j)}(K^*) &\stackrel{(h)}{\le}  \left(1-\frac{2\eta_l \mu^2 \sigma_{\min }(R)}{\left\|\Sigma_{K_j^*}\right\|}\right) (C^{(j)}(K_0) - C^{(j)}(K_j^*))\notag\\
&+2\eta_l \min\{n_x, n_u\}(\epsilon_1\bar{h}^{1}_{\text{het}} +\epsilon_2\Bar{h}^{2}_{\text{het}})^2,
\end{align}
where $(d)$ is due to Eq.~\eqref{eq:sum_expand}; $(e)$ is due to Lemma~\ref{lem:gradient_het}; $(f)$ is due to Eq.\eqref{eq:vector_youngs_inner} with $\zeta =\frac{1}{2};$ $(g)$ is due to the choice of step-size such that $\eta_l^2 \bar{h}_{\text{grad}}\le \frac{\eta_l}{4}$, which holds when $\eta_l \le \frac{1}{4\bar{h}_{\text{grad}}};$ and $(h)$ is due to Lemma~\ref{lemma:gradient_domination} and the fact that $K_n \in \mathcal{G}^0.$ If $\epsilon_1$ and $\epsilon_2$ are small enough  that $$(\epsilon_1 \Bar{h}^1_{\text{het}} + \epsilon_2 \Bar{h}^2_{\text{het}})^2 \le \min_{j \in [M]} \left\{\frac{\mu^2 \sigma_{\min }(R)\left(C^{(j)}(K_0) - C^{(j)}(K_j^*)\right)}{4||\Sigma_{K_j^*}||\min\{n_x, n_u\}}\right\},$$ we have that
\begin{align*}
C^{(j)}(K_{n,1}^{(i)}) &- C^{(j)}(K_n)\le C^{(j)}(K_0) - C^{(j)}(K_j^*),
\end{align*}
holds for any $j \in [M]$. 

The above inequality implies $K_{n,1}^{(i)}\in \mathcal{G}^0$ as long as $K_n \in \mathcal{G}^0.$ Then we can use the induction method to obtain that $K_{n,2}^{(i)}\in \mathcal{G}^0$ since $K_{n,1}^{(i)}\in \mathcal{G}^0$. As a result, an identical argument can be used from $K_{n,1}^{(i)}$ to $K_{n,2}^{(i)}$. Therefore, by repeating this step for $L$ times, we have that all the local polices $K_{n,l}^{(i)} \in \mathcal{G}^0$ holds for all $i \in [M]$ and $l = 1, \cdots, L$, when the global policy $K_n \in \mathcal{G}^0.$ 
\end{proof}

\begin{lemma}\label{lem：drift_term} 
(Drift term analysis) If $\eta_l \le\min\left\{\frac{1}{4\Bar{h}_{\text {grad}}}, \frac{1}{2}, \frac{\underline{h}_{\Delta}}{\bar{h}_1}, \frac{\log 2}{L(3\Bar{h}_{\text {grad}}+1)}\right\}$ and $K_n \in \mathcal{G}^0$, the difference between the local policy and global policy can be bounded as follows $\forall i \in [M]$ and $l \in [L]$:
\begin{align*}
\Big\lVert K_{n, l}^{(i)}-K_{n}\Big\rVert_F^2 \le 2\eta_l L  \Big\lVert\nabla C^{(i)}(K_{n})\Big\rVert_F^2 = \frac{2\eta}{\eta_g}  \Big\lVert\nabla C^{(i)}(K_{n})\Big\rVert_F^2. 
\end{align*}
\end{lemma}
\begin{proof} We have
\begin{align}
\Big\lVert K_{n, l}^{(i)}-K_{n}\Big\rVert_F^2 &= \Big\lVert K^{(i)}_{n,l-1}-K_{n} - \eta_l\nabla {C^{(i)}(K^{(i)}_{n,l-1})}\Big\rVert_F^2 \notag\\
& = \Big\lVert K^{(i)}_{n,l-1}-K_{n} \Big\rVert_F^2 -2\eta_l  \left[\left\langle \nabla {C^{(i)}(K^{(i)}_{n,l-1})}, K^{(i)}_{n,l-1}-K_{n}\right\rangle\right] \notag\\
& + \Big\lVert \eta_l\nabla {C^{(i)}(K^{(i)}_{n,l-1})}\Big\rVert_F^2 \notag\\
&= \Big\lVert K^{(i)}_{n,l-1}-K_{n} \Big\rVert_F^2  
-2\eta_l \left[\left\langle \nabla C^{(i)}(K^{(i)}_{n,l-1})-\nabla C^{(i)}(K_{n}), K^{(i)}_{n,l-1}-K_{n}\right\rangle\right] \notag\\
&- 2\eta_l  \left[\left\langle\nabla C^{(i)}(K_{n}), K^{(i)}_{n,l-1}-K_{n}\right\rangle\right]
+ \Big\lVert \eta_l\nabla {C^{(i)}(K^{(i)}_{n,l-1})}\Big\rVert_F^2 \notag\\
& \le \Big\lVert K^{(i)}_{n,l-1}-K_{n} \Big\rVert_F^2 
 +2\eta_l \Big\lVert\nabla C^{(i)}(K^{(i)}_{n,l-1})-\nabla C^{(i)}(K_{n})\Big\rVert_F \Big\lVert K^{(i)}_{n,l-1}-K_{n}\Big\rVert_F\notag\\
 &+ 2\eta_l  \Big\lVert\nabla C^{(i)}(K_{n})\Big\rVert_F\Big\lVert K^{(i)}_{n,l-1}-K_{n}\Big\rVert_F+ \Big\lVert \eta_l\nabla {C^{(i)}(K^{(i)}_{n,l-1})}\Big\rVert_F^2\notag\\
& \stackrel{(a)}{\le} \Big\lVert K^{(i)}_{n,l-1}-K_{n} \Big\rVert_F^2
+2\eta_l h_{\text {grad}}(K_n) \Big\lVert K^{(i)}_{n,l-1}- K_{n}\Big\rVert_F \Big\lVert K^{(i)}_{n,l-1}-K_{n}\Big\rVert_F\notag\\
&+\eta_l  \Big\lVert\nabla C^{(i)}(K_{n})\Big\rVert_F^2 +\eta_l\Big\lVert K^{(i)}_{n,l-1}-K_{n}\Big\rVert_F^2+ \Big\lVert \eta_l\nabla {C^{(i)}(K^{(i)}_{n,l-1})}\Big\rVert_F^2\notag\\
& \le \left(1+2\eta_l h_{\text {grad}}(K_n)+  \eta_l\right)\Big\lVert K^{(i)}_{n,l-1}-K_{n} \Big\rVert_F^2 +(\eta_l +2\eta_l^2) \Big\lVert\nabla C^{(i)}(K_{n})\Big\rVert_F^2 \notag\\
& + 2\eta_l^2\Big\lVert \nabla C^{(i)}(K^{(i)}_{n,l-1}) - \nabla C^{(i)}(K_{n})\Big\rVert_F^2 \notag\\
&\stackrel{(b)}{\le} \left(1+2\eta_l h_{\text {grad}}(K_n)+ \eta_l\right)\Big\lVert K^{(i)}_{n,l-1}-K_{n} \Big\rVert_F^2 +(\eta_l+2\eta_l^2) \Big\lVert\nabla C^{(i)}(K_{n})\Big\rVert_F^2 \notag\\
&+ 2 \eta_l^2 h^2_{\text {grad}}(K_n)\Big\lVert K^{(i)}_{n,l-1}-K_{n} \Big\rVert_F^2 \notag\\
&\stackrel{(c)}{\le} \left(1+2\eta_l \Bar{h}_{\text {grad}}+  \eta_l+2\eta_l^2 \Bar{h}_{\text {grad}}^2\right)\Big\lVert K^{(i)}_{n,l-1}-K_{n} \Big\rVert_F^2 +(\eta_l+ 2\eta_l^2) \Big\lVert\nabla C^{(i)}(K_{n})\Big\rVert_F^2 \notag\\
& \stackrel{(d)}{\le} \left(1+3\eta_l \bar{h}_{\text {grad}}+  \eta_l\right)\Big\lVert K^{(i)}_{n,l-1}-K_{n} \Big\rVert_F^2 +2\eta_l \Big\lVert\nabla C^{(i)}(K_{n})\Big\rVert_F^2,
\label{eq:drift_term_model_based}
\end{align}
where $(a)$ and $(b)$ are due to Lemma~\ref{Lemma:lipschitz}; $(c)$ is due to the fact that $K_n \in \mathcal{G}^0$; $(d)$ is due to the choice of step-size such that $2\eta_l^2 \Bar{h}_{\text {grad}}^2 \le \eta_l\Bar{h}_{\text {grad}}$ and $2\eta_l^2 \le \eta_l$, which hold when $\eta_l \le\min\{\frac{1}{2\Bar{h}_{\text {grad}}}, \frac{1}{2}\}.$ Therefore, we have
\begin{align*}
\Big\lVert K_{n, l}^{(i)}-K_{n}\Big\rVert_F^2 
&\le (1+ 3\eta_l\Bar{h}_{\text {grad}} +\eta_l)\Big\lVert K^{(i)}_{n,l-1}-K_{n} \Big\rVert_F^2+2\eta_l \Big\lVert\nabla C^{(i)}(K_{n})\Big\rVert_F^2 \notag\\
& \le (1+ 3\eta_l\Bar{h}_{\text {grad}} +\eta_l)^{l}\underbrace{\Big\lVert K^{(i)}_{n,0}-K_{n} \Big\rVert_F^2}_{= 0}+2\sum_{j=0}^{l-1}\left(1+ 3\eta_l\Bar{h}_{\text {grad}} +\eta_l \right)^j \eta_l \Big\lVert\nabla C^{(i)}(K_{n})\Big\rVert_F^2 \notag\\
& \le \frac{\left(1+ 3\eta_l\Bar{h}_{\text {grad}} +\eta_l \right)^l -1}{\left(1+ 3\eta_l\Bar{h}_{\text {grad}} +\eta_l \right)-1} 2\eta_l \Big\lVert\nabla C^{(i)}(K_{n})\Big\rVert_F^2 \notag\\
& \stackrel{(a)}{\le} 2\times\frac{1+ l(3\eta_l\Bar{h}_{\text {grad}} +\eta_l)-1}{3\Bar{h}_{\text {grad}} +1}  \Big\lVert\nabla C^{(i)}(K_{n})\Big\rVert_F^2  \notag\\
& \le 2\eta_l L \Big\lVert\nabla C^{(i)}(K_{n})\Big\rVert_F^2,
\end{align*}
where, for $(a)$, we used the fact that $(1+x)^{\tau+1} \leq 1+2 x(\tau+1)$ holds for $x \leq \frac{\log 2}{\tau}$. In other words, $\left(1+ 3\eta_l\Bar{h}_{\text {grad}} +\eta_l \right)^l \le 1+ l(3\eta_l\Bar{h}_{\text {grad}} +\eta_l)$ holds when $3\eta_l\Bar{h}_{\text {grad}} +\eta_l \le \frac{\log 2}{l},$ i.e., when $\eta_l \le \frac{\log 2}{L(3\Bar{h}_{\text {grad}}+1)}.$ \end{proof}

\begin{lemma}\label{lem:model_based_prp}
(Per round progress) Suppose $K_n \in\mathcal{G}^0.$ If we choose the local step-size as 
$$\eta_l =\frac{1}{2} \min\left\{\frac{1}{4\Bar{h}_{\text {grad}}}, \frac{1}{2}, \frac{\underline{h}_{\Delta}}{\bar{h}_1}, \frac{\log 2}{L(3\Bar{h}_{\text {grad}}+1)}, \frac{1}{80L\bar{h}^2_{\text {grad}}}\right\},$$ 
choose $\eta = \frac{1}{2}\min\{\frac{\underline{h}_{\Delta}}{\bar{h}_1},1,\frac{2}{3\Bar{h}_{\text {grad}}} \}$, and the global step-size as $\eta_g = \frac{\eta}{L \eta_l}$, then, for all $i \in [M]$, it holds that
\begin{equation}\label{eq:proof_model_based_prp}
\begin{aligned}
C^{(i)}(K_{n+1}) - C^{(i)}(K_n) 
&\le - \frac{\eta \mu^2 \sigma_{\min }(R)}{\left\|\Sigma_{K_i^*}\right\|_{\max}} (C^{(i)}(K_{n}) - C^{(i)}(K_i^*) ) +3\eta \min\{n_x, n_u\} (\epsilon_1\bar{h}^{1}_{\text{het}} +\epsilon_2\bar{h}^{2}_{\text{het}})^2.
\end{aligned}
\end{equation}
\end{lemma}
\begin{proof}
\begin{align*}
&C^{(i)}(K_{n+1}) - C^{(i)}(K_n) 
\stackrel{(a)}{\le} \langle \nabla C^{(i)}(K_n), K_{n+1}-K_n \rangle +\frac{h_{\text {grad}}(K_n)}{2}\lVert K_{n+1}-K_n\rVert_F^2 \notag\\
&= -\left\langle\nabla C^{(i)}(K_n), \frac{\eta}{ML}\sum_{j=1}^M \sum_{l=0}^{L-1} \nabla C^{(j)}(K_{n, l}^{(j)})\right\rangle +\frac{h_{\text {grad}}(K_n)}{2} \Big\lVert \frac{\eta}{ML}\sum_{j=1}^M \sum_{l=0}^{L-1}{\nabla C^{(j)}(K_{n, l}^{(j)})}\Big\rVert_F^2 \notag\\
& =-\left\langle\nabla C^{(i)}(K_n), \frac{\eta}{ML}\sum_{j=1}^M \sum_{l=0}^{L-1} \left[\nabla C^{(j)}(K_{n, l}^{(j)})-\nabla C^{(i)}(K_{n})\right]\right\rangle -\eta \Big\lVert \nabla C^{(i)}(K_n)\Big\rVert_F^2\notag\\
&+\frac{h_{\text {grad}}(K_n)}{2} \Big\lVert \frac{\eta}{ML}\sum_{j=1}^M \sum_{l=0}^{L-1}{\nabla C^{(j)}(K_{n, l}^{(j)})}\Big\rVert_F^2 \notag\\
& =-\left\langle\nabla C^{(i)}(K_n), \frac{\eta}{ML}\sum_{j=1}^M \sum_{l=0}^{L-1} \left[\nabla C^{(j)}(K_{n, l}^{(j)})-\nabla C^{(j)}(K_{n})+\nabla C^{(j)}(K_{n})-\nabla C^{(i)}(K_{n})\right]\right\rangle \notag\\
&-\eta \Big\lVert \nabla C^{(i)}(K_n)\Big\rVert_F^2+\frac{h_{\text {grad}}(K_n)}{2} \Big\lVert \frac{\eta}{ML}\sum_{j=1}^M \sum_{l=0}^{L-1}{\nabla C^{(j)}(K_{n, l}^{(j)})}\Big\rVert_F^2 \notag\\
& \le \eta \Big\lVert\nabla C^{(i)}(K_n)\Big\rVert_F \Big\lVert\frac{1}{ML}\sum_{j=1}^M \sum_{l=0}^{L-1} \left[\nabla C^{(j)}(K_{n, l}^{(j)})-\nabla C^{(j)}(K_{n})\right]\Big\rVert_F\notag\\ &+\frac{\eta}{M} \sum_{j=1}^M\Big\lVert\nabla C^{(i)}(K_n)\Big\rVert_F\Big\lVert \nabla C^{(j)}(K_{n})-\nabla C^{(i)}(K_{n})\Big\rVert_F\notag\\
&-\eta \Big\lVert \nabla C^{(i)}(K_n)\Big\rVert_F^2
+\frac{h_{\text {grad}}(K_n)}{2} \Big\lVert \frac{\eta}{ML}\sum_{j=1}^M \sum_{l=0}^{L-1}{\nabla C^{(j)}(K_{n, l}^{(j)})}\Big\rVert_F^2 \notag\\
& \stackrel{(b)}{\le}\eta \Big\lVert\nabla C^{(i)}(K_n)\Big\rVert_F   \left[\frac{h_{\text {grad}}(K_n)}{ML}\sum_{j=1}^M \sum_{l=0}^{L-1}  \Big\lVert K_{n, l}^{(j)}-K_{n}\Big\rVert_F\right] \notag\\
&+\frac{\eta}{4}\Big\lVert\nabla C^{(i)}(K_n)\Big\rVert_F^2 +\frac{\eta}{M}\sum_{j=1}^M\Big\lVert\nabla C^{(j)}(K_n)-\nabla C^{(i)}(K_n)\Big\rVert_F^2-\eta \Big\lVert \nabla C^{(i)}(K_n)\Big\rVert_F^2\notag\\
&+ h_{\text {grad}}(K_n)\frac{3\eta^2}{2ML}\sum_{j=1}^M \sum_{l=0}^{L-1}\Big\lVert {\nabla C^{(j)}(K_{n, l}^{(j)})}-{\nabla C^{(j)}(K_n)}\Big\rVert_F^2 \notag\\
&+ \frac{3\eta^2 h_{\text{grad}}(K_n)}{2M} \sum_{j=1}^M\Big\lVert\nabla C^{(j)}(K_n)-\nabla C^{(i)}(K_n)\Big\rVert_F^2 +\frac{3\eta^2h_{\text{grad}}(K_n)}{2}\Big\lVert\nabla C^{(i)}(K_n)\Big\rVert_F^2\notag\\
& \stackrel{(c)}{\le} \frac{\eta}{4}\Big\lVert \nabla C^{(i)}(K_n)\Big\rVert_F^2+\frac{\eta \bar{h}^2_{\text {grad}}}{ML}\sum_{i=1}^M \sum_{l=0}^{L-1} \Big\lVert K_{n, l}^{(i)}-K_{n}\Big\rVert_F^2 \notag\\
& +\frac{\eta}{4}\Big\lVert\nabla C^{(i)}(K_n)\Big\rVert_F^2 +\left(\eta+ \frac{3\eta^2\bar{h}_{\text {grad}}}{2}\right)\frac{1}{M}\sum_{j=1}^M\Big\lVert\nabla C^{(j)}(K_n)-\nabla C^{(i)}(K_n)\Big\rVert_F^2\notag\\
&-\eta \Big\lVert \nabla C^{(i)}(K_n)\Big\rVert_F^2
+ \frac{3\eta^2\bar{h}^2_{\text {grad}}}{2ML}\sum_{j=1}^M \sum_{l=0}^{L-1}\Big\lVert K_{n, l}^{(j)}-K_n\Big\rVert_F^2 + \frac{\eta}{8}\Big\lVert \nabla C^{(i)}(K_n)\Big\rVert_F^2\notag\\
&\stackrel{(d)}{\le} - \frac{3\eta}{8}\Big\lVert \nabla C^{(i)}(K_n)\Big\rVert_F^2 + \frac{5\eta^2\bar{h}^2_{\text {grad}}}{\eta_g M}\sum_{j=1}^M \Big\lVert \nabla C^{(j)}(K_n)\Big\rVert_F^2 \notag \\ 
&+ 2\eta \min\{n_x, n_u\} (\epsilon_1\bar{h}^{1}_{\text{het}} +\epsilon_2\bar{h}^{2}_{\text{het}})^2\notag\\
&\stackrel{(e)}{\le} -\frac{3\eta}{8}\Big\lVert \nabla C^{(i)}(K_n)\Big\rVert_F^2 +\frac{10\eta^2\bar{h}^2_{\text {grad}}}{\eta_g M}\sum_{j=1}^M \Big\lVert \nabla C^{(j)}(K_n)-\nabla C^{(i)}(K_n)\Big\rVert_F^2\notag\\
&+\frac{10\eta^2\bar{h}^2_{\text {grad}}}{\eta_g } \Big\lVert \nabla C^{(i)}(K_n)\Big\rVert_F^2+ 2\eta \min\{n_x, n_u\} (\epsilon_1\bar{h}^{1}_{\text{het}} +\epsilon_2\bar{h}^{2}_{\text{het}})^2\notag\\
& \stackrel{(f)}{\le} -\frac{\eta}{4}\Big\lVert \nabla C^{(i)}(K_n)\Big\rVert_F^2+3\eta \min\{n_x, n_u\} (\epsilon_1\bar{h}^{1}_{\text{het}} +\epsilon_2\bar{h}^{2}_{\text{het}})^2\notag\\
&\le - \frac{\eta \mu^2 \sigma_{\min }(R)}{\left\|\Sigma_{K_i^*}\right\|} (C^{(i)}(K_{n}) - C^{(i)}(K_i^*) )+3\eta \min\{n_x, n_u\} (\epsilon_1\bar{h}^{1}_{\text{het}} +\epsilon_2\bar{h}^{2}_{\text{het}})^2.
\end{align*}

In the above steps, (a) is due to the choice of step-size $\eta$ such that 

$$\|K_{n+1}- K_n\| = \|\frac{\eta}{ML}\sum_{i=1}^M \sum_{l=0}^{L-1} \nabla C^{(i)}(K_{n, l}^{(i)})\| \le \eta \bar{h}_1 \le \underline{h}_{\Delta},$$
holds when $\eta \le \frac{\underline{h}_{\Delta}}{\bar{h}_1}.$ For $(b)$, we use the Lipschitz property of the gradient (Lemma \ref{Lemma:lipschitz}) in the first line, and use Eq.\eqref{eq:vector_youngs_inner} with $\zeta =\frac{1}{2}$ in the second line, and for the third and forth lines we use Eq.~\eqref{eq:sum_expand}; $(c)$ is due to Lemma~\ref{Lemma:lipschitz} and $\frac{3\eta^2\bar{h}_{\text {grad}}}{2}\le \frac{\eta}{8};$ $(d)$ is due to Lemma~\ref{lem:gradient_het}, Lemma~\ref{lem：drift_term} and the choice of step-size such that $\frac{3\eta^2\bar{h}_{\text {grad}}}{2}\le \frac{\eta}{8} \le \eta$; $(e)$ is due to Eq.\eqref{eq:sum_expand}; and for $(f)$, we use the fact that $\frac{10\eta^2\bar{h}^2_{\text {grad}}}{\eta_g } \le \frac{\eta}{8} \le \eta,$ which holds when $\eta_l \le \frac{1}{80L\bar{h}^2_{\text {grad}}}$. We use the gradient domination property (Lemma~\ref{lemma:gradient_domination}) in the last equality. 
\end{proof}

With this lemma, we are now ready to provide the convergence guarantees for \texttt{FedLQR} in the model-based setting. 

\paragraph{Proof of the iterative stability guarantees of \texttt{FedLQR}:}  Here we leverage the method of induction to prove 
\texttt{FedLQR}'s iterative stability guarantees. First, we start from an initial policy $K_0 \in \mathcal{G}^0$. At round $n$, we assume $K_n \in \mathcal{G}^0.$ According to Lemma~\ref{lem:stability_local}, we can show that all the local policies $K_{n,l}^{(i)} \in \mathcal{G}^0$. Furthermore, by choosing the step-sizes properly in Lemma~\ref{lem:model_based_prp}, we have that 
\begin{equation}
\begin{aligned}
C^{(i)}(K_{n+1}) - C^{(i)}(K_n) 
&\le - \frac{\eta \mu^2 \sigma_{\min }(R)}{\left\|\Sigma_{K_i^*}\right\|_{\max}} (C^{(i)}(K_{n}) - C^{(i)}(K_i^*) )\notag\\&+3\eta \min\{n_x, n_u\} (\epsilon_1\bar{h}^{1}_{\text{het}} +\epsilon_2\bar{h}^{2}_{\text{het}})^2.
\end{aligned}
\end{equation}
for any $i \in [M]$.  Then, for any fixed system $i \in [M]$, with $(\epsilon_1 \Bar{h}^1_{\text{het}} + \epsilon_2 \Bar{h}^2_{\text{het}})^2 \le\Bar{h}^3_{\text{het}},$ we have that 
\begin{equation}
\begin{aligned}
C^{(i)}(K_{n+1}) - C^{(i)}(K_i^*) 
&\le\left(1- \frac{\eta \mu^2 \sigma_{\min }(R)}{\left\|\Sigma_{K_i^*}\right\|}\right) (C^{(i)}(K_{n}) - C^{(i)}(K_i^*))\notag \\
&+3\eta \min\{n_x, n_u\} (\epsilon_1\bar{h}^{1}_{\text{het}} +\epsilon_2\bar{h}^{2}_{\text{het}})^2 \label{eq:loss_decrease}\\
&\le\left(1- \frac{\eta \mu^2 \sigma_{\min }(R)}{\left\|\Sigma_{K_i^*}\right\|}\right) (C^{(i)}(K_{0}) - C^{(i)}(K_i^*))\notag\\
&+\frac{\eta \mu^2 \sigma_{\min }(R)}{\left\|\Sigma_{K_i^*}\right\|}(C^{(i)}(K_{0}) - C^{(i)}(K_i^*))\notag\\
& \le C^{(i)}(K_{0}) - C^{(i)}(K_i^*).
\end{aligned}
\end{equation}

With this, we can easily see that the global policy $K_{n+1}$ at the next round $n+1$ is also stabilizing, i.e., $K_{n+1} \in \mathcal{G}^0,$ by using the definition of $\mathcal{G}^0$ (Definition \ref{def:stable_set}). Therefore, we can complete proving \texttt{FedLQR}'s iterative stability property by inductive reasoning.

\paragraph{Proof of \texttt{FedLQR}'s convergence:} From Eq.\eqref{eq:proof_model_based_prp}, we have 
\begin{align*}
C^{(i)}(K_{n+1}) - C^{(i)}(K_i^*) 
&\le\left(1- \frac{\eta \mu^2 \sigma_{\min }(R)}{\left\|\Sigma_{K_i^*}\right\|}\right) (C^{(i)}(K_{n}) - C^{(i)}(K_i^*))\\ 
&+3\eta \min\{n_x, n_u\} (\epsilon_1\bar{h}^{1}_{\text{het}} +\epsilon_2\bar{h}^{2}_{\text{het}})^2.  
\end{align*}

Under the assumptions in Lemma \ref{lem:model_based_prp}, \texttt{FedLQR} thus 
enjoys the following convergence guarantee after $N$ rounds:
\begin{align*}
C^{(i)}(K_{N}) - C^{(i)}(K_i^*)  &\le \left(1 -\frac{\eta \mu^2 \sigma_{\min }(R)}{\left\|\Sigma_{K_i^*}\right\|} \right)^N(C^{(i)}(K_0) - C^{(i)}(K_i^*))\\
&+ \frac{3\min\{n_x, n_u\} \left\|\Sigma_{K_i^*}\right\|}{\mu^2 \sigma_{\min}(R)} (\epsilon_1 h^1_{\text{het}} + \epsilon_1 h^2_{\text{het}})^2.
\end{align*}
Thus, we finish the proof of Theorem~\ref{thm:one_local_model_based} with $c_{\text{uni},1}=12$ and $\mathcal{B}(\epsilon_1, \epsilon_2):= \frac{\upsilon \left\|\Sigma_{K_i^*}\right\|}{4\mu^2 \sigma_{\min}(R)} (\epsilon_1 h^1_{\text{het}} + \epsilon_1 h^2_{\text{het}})^2$. $\hfill \Box$

\subsubsection{Proof of 
Corollary 8.3}\label{sec:proof of coro}

\begin{proof}
From Theorem~\ref{thm:one_local_model_based}, we have that
\begin{align*}
C^{(i)}(K_{N}) - C^{(i)}(K_i^*)  &\le \left(1 -\frac{\eta \mu^2 \sigma_{\min }(R)}{\left\|\Sigma_{K_i^*}\right\|} \right)^N(C^{(i)}(K_0) - C^{(i)}(K_i^*))\\
&+ \frac{3\min\{n_x, n_u\} \left\|\Sigma_{K_i^*}\right\|}{\mu^2 \sigma_{\min}(R)} (\epsilon_1 h^1_{\text{het}} + \epsilon_1 h^2_{\text{het}})^2
\end{align*}
holds for all $i \in [M].$ Based on the fact that $K_i^*$ is system $i$'s optimal LQR controller, i.e., $C^{(i)}(K_i^*) \le C^{(i)}(K^*),$ we have that
\begin{align}\label{eq:midd}
C^{(i)}(K_{N}) - C^{(i)}(K^*)  &\le \left(1 -\frac{\eta \mu^2 \sigma_{\min }(R)}{\left\|\Sigma_{K_i^*}\right\|} \right)^N(C^{(i)}(K_0) - C^{(i)}(K_i^*))\notag\\
&+ \frac{3\min\{n_x, n_u\} \left\|\Sigma_{K_i^*}\right\|}{\mu^2 \sigma_{\min}(R)} (\epsilon_1 h^1_{\text{het}} + \epsilon_1 h^2_{\text{het}})^2
\end{align}
We can finish the proof of Corollary~\ref{coro:global convergence} by averaging Eq.~\eqref{eq:midd} across all $M$ systems.
\end{proof}

\newpage
\subsection{Zeroth-order optimization} \label{sec:zeroth_order_optimization}

To prepare for the model-free setting where the controllers only have access to the system's trajectories, we first quickly recap the basic idea behind zeroth-order optimization. Say our goal is to minimize a loss function $f(x)$, where $x\in\mathbb{R}^d$. When one has access to exact deterministic gradients of this loss function via an oracle, the standard approach for minimization would be to query the gradient oracle at each iteration, and run gradient descent. Concretely, one would run the following iterative scheme: $x_{t+1}=x_{t}-\eta \nabla f(x_t)$, where $\eta$ is a suitably chosen learning-rate/step-size. While such first-order optimization schemes have a rich history, there has also been a growing interest in understanding the behavior of derivative-free (zeroth-order) methods that can \textit{only query function values}, as opposed to the gradients. Two immediate reasons (among many) for studying zeroth-order optimization are as follows: (i) in practice, one may only have access to a black-box procedure that cannot evaluate gradients; and (ii) computing gradients might prove to be too computationally-expensive. 

Given two or more function evaluations, the basic idea behind zeroth-order algorithms is to construct an estimate of the true gradient for evaluating and updating model parameters. For instance, a typical zeroth-order scheme with single-point function evaluation would take the following form~\citep{polyak1987introduction}:
$$  x_{t+1}=x_t - \eta_t \left( \frac{f(x_t+\mu_t u)-f(x_t)}{\mu_t}\right)u.$$

In the expression above, $\{\eta_t\}$ is the learning-rate sequence, $\{\mu_t\}$ is a sequence typically chosen in a way such that $\mu_t \rightarrow 0$, and $u$ is a random vector distributed uniformly over the unit sphere. For details about the convergence of zeroth-order optimization algorithms such as the one above, we refer the interested reader to~\citep{nesterov2017random,duchi2015optimal,bach2016highly}. 

We now turn to briefly describing the model-free setup for our LQR problem. \citep{fazel2018global} propose a zeroth-order-based algorithm (Algorithm 1 in \citep{fazel2018global}) to compute an estimation $\widehat{\nabla C(K)}$ and $\widehat{\Sigma_K}$ for both ${\nabla C(K)}$ and ${\Sigma_K}$, for a given $K$. Algorithm 1 in \citep{fazel2018global} exploits a multiple-trajectory-based technique that uses a Gaussian perturbed cost function (i.e., producing a Gaussian smoothing function) to estimate ${\nabla C(K)}$ from cost function perturbed values. That is, given the cost function $C(K)$, we can  define its perturbed function as,
$$
C_r(K)=\mathbb{E}_{U \sim \mc{B}_r}[C(K+U)]
$$
where $\mc{B}_r$ is the uniform distribution over all matrices with Frobenius norm at most $r$ and $U$ is a random matrix with proper dimension and generated from $\mc{B}_r$.
For small $r$, the smooth cost $C_r(K)$ is a good approximation to the original cost $C(K)$. Due to the Gaussian smoothing, the gradient has a particularly simple functional form~\citep{gravell2020learning}:
$$\nabla C_r(K)=\frac{n_x n_u}{r^2} \mathbb{E}_{U \sim \mathbb{B}_r}[C(K+U) U].$$
Therefore, this expression implies a straightforward method to obtain an unbiased estimate of $\nabla C_r(K)$, through obtaining the infinite-horizon rollouts. However, in practice, we can only obtain the finite-horizon rollouts to approximate the gradient. Thanks to~\citep{fazel2018global}, they showed that the approximation error of the exact gradient can be reduced to  arbitrary accuracy if the number of sample trajectories $n_s$ and the length of each rollout $\tau$ are sufficiently large, and the smoothing radius $r$ is small enough.

We also emphasize that the zeroth-order estimator obtained from $C_r(K)$ is a biased approximation of the true gradient $\nabla C(K)$, with the bias arising from the smoothing step. The magnitude of this bias is controlled by the choice of the smoothing radius $r$: smaller values of $r$ yield a closer approximation to the original cost landscape at the cost of higher variance. In our analysis, we make the estimation error sufficiently small by selecting $r$, the number of samples $n_s$, and the rollout horizon $\tau$ appropriately.

It is important to emphasize that our results can be readily extended to the two-point zeroth-order estimation scheme, where the estimation variance can be reduced, as discussed in \cite{malik2019derivative}. However, addressing the two-point zeroth-order estimation is beyond the scope of this work, and we refer the reader to \citep{toso2024asynchronous} for an instantiation of the results in the asynchronous aggregation setting under the two-point gradient estimation scheme.

\newpage
\subsection{The model-free setting} \label{sec:proof_model_free}

For notational brevity we rewrite $\widehat{\nabla C^{(i)}(K)}$ as $\tilde{\nabla} C^{(i)}(K)$ where 
\begin{align*}
\widehat{\nabla C^{(i)}(K)}=\tilde{\nabla} C^{(i)}(K) &:=\frac{1}{n_{s}} \sum_{s=1}^{n_{s}} \frac{n_x n_u}{r^2} \tilde{C}^{(i),(\tau)}\left(K+U^{(i)}_s\right) U^{(i)}_s, \nonumber 
\end{align*}
and introduce two new gradient-based terms:
\begin{align*}
\nabla^{\prime} C^{(i)}(K)&:=\frac{1}{n_{s}} \sum_{s=1}^{n_{s}} \frac{n_x n_u}{r^2} C^{(i),(\tau)}\left(K+U^{(i)}_s\right) U^{(i)}_s, \nonumber\\
\hat{\nabla} C^{(i)}(K)&:=\frac{1}{n_{s}} \sum_{s=1}^{n_{s}} \frac{n_x n_u}{r^2} C^{(i)}\left(K+U^{(i)}_s\right) U^{(i)}_s,
\end{align*}
where $
\tilde{C}^{(i),(\tau)}\left(K+U^{(i)}_s\right):=\sum_{t=0}^{\tau-1}\left(x^{(i)\top}_{t} Q x^{(i)}_{t}+u^{(i) \top}_{t} R u^{(i)}_{t}\right)$ with $x^{(i)}_{t} = (K+U^{(i)}_s)u^{(i)}_{t}$, $
{C}^{(i),(\tau)}\left(K+U^{(i)}_s\right):=\mc{E}_{x_0^{(i)}\sim \mathcal{D}}\sum_{t=0}^{\tau-1}\left(x^{(i)\top}_{t} Q x^{(i)}_{t}+u^{(i) \top}_{t} R u^{(i)}_{t}\right)$ and $$
{C}^{(i)}\left(K+U^{(i)}_s\right):=\mc{E}_{x_0^{(i)}\sim \mathcal{D}}\sum_{t=0}^{\infty}\left(x^{(i)\top}_{t} Q x^{(i)}_{t}+u^{(i) \top}_{t} R u^{(i)}_{t}\right).$$
\subsubsection{Auxiliary Lemmas}\label{sec:aux}
\begin{lemma}\label{lemma:finite_horizon_approximation} (Approximating $C^{(i)}(K)$ and $\Sigma^{(i)}_K$ with finite horizon) Suppose $K$ is such that $C^{(i)}(K)$ is finite. Define the finite horizon estimates, 
\begin{align*}
\Sigma_K^{(i),(\tau)}:=\mathbb{E}\left[\sum_{t=0}^{\tau-1} x^{(i)}_t x_t^{(i)\top}\right] \quad \text{and}\quad C^{(i),(\tau)}(K):=\mathbb{E}\left[\sum_{t=0}^{\tau-1} x_t^{(i)\top} Q x^{(i)}_t+u_t^{(i)\top} R u^{(i)}_t\right],
\end{align*}
for all systems $i \in [M]$. Now, let $\epsilon$ be an arbitrarily small constant such that
\begin{align*}
\tau \geq h_{\tau}^1(\epsilon)& :=\max_{i \in [M]}\left\{\frac{n_x \cdot (C^{(i)}(K))^2}{\epsilon \mu (\sigma_{\text{min}}(Q))^2 }\right\}= \frac{n_x \cdot (C_{\max}(K))^2}{\epsilon \mu (\sigma_{\text{min}}(Q))^2 }, \end{align*} 
such that
$$\left\|\Sigma_K^{(i),(\tau)}-\Sigma^{(i)}_K\right\| \leq \epsilon.$$ 

If
\begin{align*}
\tau \geq h_{\tau}^2 (\epsilon)&:=\max_{i \in [M]}\left\{\frac{n_x \cdot (C^{(i)}(K))^2(\|Q\|+\|R\|\|K\|^2)}{\epsilon \mu (\sigma_{\text{min}}(Q))^2 }\right\} \\
&=\frac{n_x \cdot (C_{\max}(K))^2(\|Q\|+\|R\|\|K\|^2)}{\epsilon \mu (\sigma_{\text{min}}(Q))^2 },
\end{align*} we have
\begin{align*}
\ \left|C^{(i),(\tau)}(K)-C^{(i)}(K)\right| \leq \epsilon,
\end{align*}
where $C_{\max}(K):= \max_{i\in [M]} C^{(i)}(K)$.
\end{lemma}
\begin{proof} The proof for this lemma is detailed in \citep[Lemma 23]{fazel2018global}. 
\end{proof}

\begin{lemma}\label{lemma:infinite} (Estimating $\nabla C^{(i)}(K)$ with finitely many infinite-horizon rollouts) Given an arbitrary tolerance $\epsilon$ and probability $\delta$, suppose the radius $r$ satisfies 
$$
r \leq h_r\left(\frac{\epsilon}{2}\right):=\min \left\{\underline{h}_{\Delta}, \frac{\bar{C}_{\max}}{\bar{h}_{\text {cost }}}, \frac{\epsilon}{2\Bar{h}_{\text {grad }}}\right\},
$$
and the number of samples $n_s$ satisfies,
$$
\begin{aligned}
n_{s}  \geq h_{\text {sample }}\left(\frac{\epsilon}{2}, \delta \right) :&=\frac{8\sigma_{\hat{\nabla}}^2 \min (n_x, n_u)}{\epsilon^2} \log \left[\frac{n_x+n_u}{\delta}\right] \\
\sigma_{\hat{\nabla}}^2 & :=\left(\frac{2 n_x n_u \bar{C}_{\max}}{r}\right)^2+\left(\frac{\epsilon}{2}+\bar{h}_1\right)^2
\end{aligned}
$$
Then with a high probability of at least $1-\delta$, the estimate
$$
\hat{\nabla} C^{(i)}(K)=\frac{1}{n_{s}} \sum_{s=1}^{n_{s}} \frac{n_x n_u}{r^2} C^{(i)}\left(K+U^{(i)}_s\right) U^{(i)}_s
$$
satisfies
$$\|\hat{\nabla} C^{(i)}(K)-\nabla C^{(i)}(K)\|_F \leq \epsilon$$ 
for any system $i \in [M]$ and $K\in \mathcal{G}^0$.
\end{lemma}
\begin{proof}The proof for this lemma is detailed in Lemma B.6 of \citep{gravell2020learning}. It is worthwhile to mention that, in \citep{gravell2020learning}, the number of samples $n_s$ satisfies 
$$
\begin{aligned}
n_{s} & \ge\left[\underbrace{\frac{8\sigma_{\hat{\nabla}}^2 \min (n_x, n_u)}{\epsilon^2}}_{T_1}+\underbrace{\frac{8 \min (n_x, n_u)}{\epsilon^2}\frac{R_{\hat{\nabla}} \epsilon}{6 \sqrt{\min (n_x, n_u)}}}_{T_2}\right] \log \left[\frac{n_x+n_u}{\delta}\right] 
\end{aligned}
$$ with $R_{\hat{\nabla}}=\frac{2 n_x n_u \bar{C}_{\max}}{r}+\frac{\epsilon}{2}+\bar{h}_1.$ In the analysis throughout the paper, we only keep the dominant term $T_1$ in $n_s$, since $T_1$ is in the order $\mathcal{O}({\epsilon}^{-2})$ while $T_2$ is in the order $\mathcal{O}({\epsilon}^{-1}).$

By taking the maximum over $K$ inside $\mathcal{G}^0,$ we make the local parameters become the global parameters, e.g., $\bar{C}_{\max} :=\sup_{K \in \mathcal{G}^{0}, i\in [M]} C^{(i)}(K).$ 
\end{proof}

\begin{lemma}\label{lemma:finitely_many_finite_horizon_estimation}(Estimating $\nabla C^{(i)}(K)$ with finitely many finite-horizon rollouts): Given an arbitrary tolerance $\epsilon$ and probability $\delta$, suppose that the smoothing radius $r$ satisfies, 
$$
r \leq h_r\left(\frac{\epsilon}{4}\right)=\min \left\{\bar{h}_{\Delta}, \frac{\bar{C}_{\max}}{\bar{h}_{\text {cost}}}, \frac{\epsilon}{4 \bar{h}_{\text {grad }}}\right\},
$$
and the trajectory length $\tau$ satisfies
$$
\tau \geq h_{\tau}\left(\frac{r \epsilon}{4 n_x n_u}\right)=\frac{4 n_u n_x^2 (C_{\max}(K))^2\left(\|Q\|+\|R\|\|K\|^2\right)}{r \epsilon \mu {\sigma_{\min}(Q)}^2}.
$$
According to Assumption~\ref{assume:initial_point}, the distribution of the initial states satisfies $x_0^{(i)} \sim$ $\mathcal{D}$ and $\left\|x_0^{(i)}\right\| \leq H$ almost surely. Thus, for any given realization $x^{(i)}_{0,s}$\footnote{The notation $x^{(i)}_{0,s}$ denotes $s$-th sample of the initial state from $i$-th system.} of $x^{(i)}_0$, and for any system $i \in [M]$, we have $$\left\|x^{(i)}_{0,s}\right\| \leq H, \quad \left(x^{(i)}_{0,s}\right)\left(x^{(i)}_{0,s}\right)^{\top} \preceq \frac{H^2}{\mu} \mathbb{E}\left[x^{(i)}_0 x^{(i)\top}_0\right].$$ As a result, the summation over the finite-time horizon
$$
\sum_{t=0}^{\tau-1}\left(x^{(i)\top}_{t,j} Q x^{(i)}_{t,j}+u^{(i)\top}_{t,j} R u^{(i)}_{t,j}\right) \leq \frac{H^2}{\mu} C^{(i)}\left(K+U^{(i)}_j\right).
$$

Furthermore, suppose the number of samples $n_{s}$ satisfies
$$
\begin{aligned}
& n_{s} \geq h_{\text {sample,trunc }}\left(\frac{\epsilon}{4}, \delta, \frac{H^2}{\mu}\right) :=\frac{32 \sigma_{\tilde{\nabla}}^2\min (n_x, n_u)}{\epsilon^2} \log \left[\frac{n_x+n_u}{\delta}\right],
\end{aligned}
$$
where
$$
\begin{aligned}
& \sigma_{\tilde{\nabla}}^2:=\left(\frac{2 n_x n_u H^2 \bar{C}_{\max}}{r \mu}\right)^2+\left(\frac{\epsilon}{2}+\bar{h}_1\right)^2,
\end{aligned}
$$
then, with a high probability of at least $1-\delta$, the estimated gradient
$$
\tilde{\nabla} C^{(i)}(K):=\frac{1}{n_{s}} \sum_{s=1}^{n_{s}} \frac{n_u n_x}{r^2} \tilde{C}^{(i),(\tau)}\left(K+U^{(i)}_s\right) U_s^{(i)}
$$
satisfies $$\|\tilde{\nabla} C^{(i)}(K)-\nabla C^{(i)}(K)\|_F \leq \epsilon$$ 
for any system $i \in [M]$ and $K
\in \mathcal{G}^0$.
\end{lemma}
\begin{proof} The proof for this lemma is detailed in Lemma B.7 of \citep{gravell2020learning}. As in Lemma~\ref{lemma:infinite}, we only keep the dominant term in the requirement of sample size $n_s$. By taking the maximum over $K$ inside $\mathcal{G}^0,$ all the local parameters inside the polynomials such as  $h_r(\frac{\epsilon}{4})$ become global parameters. 
\end{proof}

\subsubsection{Proof of Lemma~\ref{lem:variance_reduction}}
\label{sec:proof_var_red_lem}
\begin{proof}
For our subsequent analysis, we will use $\mathcal{F}^n_l$
to denote the filtration that captures all the randomness
up to the $l$-th local step in round $n$. We have
\begin{align*}
&\left\|\frac{1}{ML} \sum_{i=1}^M\sum_{l=0}^{L-1}\left[\widehat{\nabla C^{(i)}(K^{(i)}_{n,l})}-\nabla C^{(i)}(K^{(i)}_{n,l})\right]\right\|_F\\
&=\left\|\frac{1}{ML} \sum_{i=1}^M\sum_{l=0}^{L-1}\left[\tilde{\nabla} C^{(i)}(K^{(i)}_{n,l})-\nabla C^{(i)}(K^{(i)}_{n,l})\right]\right\|_F\notag\\
&\le  \Big\lVert \underbrace{\frac{1}{ML}\sum_{i=1}^M\sum_{l=0}^{L}\left[\tilde{\nabla} C^{(i)}(K^{(i)}_{n,l})-\nabla^{\prime} C^{(i)}(K^{(i)}_{n,l})\right] \Big\rVert_F}_{T_1}\notag\\
 &+ \underbrace{\Big\lVert\frac{1}{ML}\sum_{i=1}^M\sum_{l=0}^{L-1}\left[\nabla^{\prime} C^{(i)}(K^{(i)}_{n,l}) -\hat{\nabla} C^{(i)}(K^{(i)}_{n,l}) \right]\Big\rVert_F}_{T_2}\notag\\
 &+   \underbrace{\Big\lVert\frac{1}{ML}\sum_{i=1}^M\sum_{l=0}^{L-1}\left[\hat{\nabla} C^{(i)}(K^{(i)}_{n,l}) -\nabla C^{(i)}(K^{(i)}_{n,l})\right]\Big\rVert_F}_{T_3}.
\end{align*}

Next, we will bound $T_1$, $T_2$, and $T_3$, respectively.
\paragraph{Bounding $T_2$:} From the proof of Lemma B.7 in~\citep{gravell2020learning}, we have
\begin{align}\label{eq:bound_T2}
T_2 \le \frac{1}{ML}\sum_{i=1}^M \sum_{l=0}^{L-1} \Big\lVert \nabla^{\prime} C^{(i)}(K^{(i)}_{n,l}) -\hat{\nabla} C^{(i)}(K^{(i)}_{n,l}) \Big\rVert_F \le \frac{\epsilon}{4}
% = \left(\frac{n_u n_x^2 \Bar{C}_{\max}^2\left(\|Q\|+\|R\|\Bar{\|K\|}_{\max}^2\right)}{r \tau \mu {\sigma_{\min}(Q)}^2}\right)^2
\end{align}
holds as long as $\tau \geq h_{\tau}\left(\frac{r \epsilon}{4 n_x n_u}\right).$

\paragraph{Bounding $T_3:$}

To precede, we bound $T_3$ as
\begin{align}
T_3 &=\Big\lVert\frac{1}{ML}\sum_{i=1}^M \sum_{l=0}^{L-1}\left[\hat{\nabla} C^{(i)}(K^{(i)}_{n,l}) -\nabla C^{(i)}(K^{(i)}_{n,l})\right]\Big\rVert_F\notag\\
&\le \frac{1}{ML}\sum_{i=1}^M \sum_{l=0}^{L-1} \underbrace{\Big\lVert \nabla C^{(i)}_r(K^{(i)}_{n,l})-\nabla C^{(i)}(K^{(i)}_{n,l})\Big\rVert_F}_{\text{Bias term }\circled{1}}\notag\\
& +  \underbrace{\Big\lVert \frac{1}{ML}\sum_{i=1}^M \sum_{l=0}^{L-1} \left[ \hat{\nabla} C^{(i)}(K^{(i)}_{n,l})-\nabla C^{(i)}_r(K^{(i)}_{n,l})\right]\Big\rVert_F}_{\text{Variance term }\circled{2}}
\end{align}
where $\nabla C^{(i)}_r(K^{(i)}_{n,l}) := \mc{E}_{U^{(i)}_{n,l} \sim \mathbb{B}_r}\left[\nabla C^{(i)}(K^{(i)}_{n,l} +U^{(i)}_{n,l}) \right].$ 

For the bias term \circled{1}, since the smoothing radius $r \leq h_r\left(\frac{\epsilon}{4}\right)$, we have that 
\begin{align}\label{eq:smoothing_error_variance}
\circled{1}= \Big\lVert \nabla C^{(i)}_r(K^{(i)}_{n,l})-\nabla C^{(i)}(K^{(i)}_{n,l})\Big\rVert_F \leq h_{\text{grad}}(K^{(i)}_{n,l}) r \leq \Bar{h}_{\text{grad}} r \le \frac{\epsilon}{4}.
\end{align}
For the variance term, \circled{2}, we will exploit the matrix Freedman inequality (Lemma~\ref{lem:Freedman}) to bound it. For simplicity, we denote 
$$e_l^{(i)}:= \frac{1}{ML}\left[\hat{\nabla} C^{(i)}(K^{(i)}_{n,l})-\nabla C^{(i)}_r(K^{(i)}_{n,l})\right], \quad e_l:=\sum_{i=1}^M e_l^{(i)},$$

Then, we have 
$$\frac{1}{ML}\sum_{i=1}^M \sum_{l=0}^{L-1} \left[ \hat{\nabla} C^{(i)}(K^{(i)}_{n,l})-\nabla C^{(i)}_r(K^{(i)}_{n,l})\right] = \sum_{l=0}^{L-1} e_l.$$

Next, we aim to prove the following claims:

\textbf{Claim I:} $Y_t: = \sum_{l=0}^t e_l $ is a martingale w.r.t $\mathcal{F}_{t-1}^n$ for $t=1,\cdots, L-1$ and $e_l:=\sum_{i=1}^M e_l^{(i)}$ is a martingale difference sequence.

\textbf{Proof:} Note that $\mc{E}[\hat{\nabla} C^{(i)}(K^{(i)}_{n,l})]=\nabla C^{(i)}_r(K^{(i)}_{n,l}).$ Then we can easily have $\mc{E}[e_l]=0$ for $l =0,\cdots, L-1.$ As a result, we have $\mc{E}[Y_t \mid \mathcal{F}_{t-1}^n] = Y_{t-1}$ since $Y_t = Y_{t-1} +e_t.$ In other words, $Y_t: = \sum_{l=0}^t e_l $ is a martingale w.r.t $\mathcal{F}_{t-1}^n$ for $t=1,\cdots, L-1$.

\textbf{Claim II:} $\left\|\mc{E}\left[e_l e_l^{\top}\biggm| \mathcal{F}_{l-1}^n\right]\right\| \le \frac{\sigma_{\hat{\nabla}}^2}{n_s ML^2}$ where $\sigma_{\hat{\nabla}}^2$ is as defined in Lemma~\ref{lemma:infinite}.

\textbf{Proof:} From Lemma B.7 in~\citep{gravell2020learning}, we can write 

\begin{align*}
\left\|\mc{E}\left[e_l^{(i)}e_l^{(i)\top}\biggm| \mathcal{F}_{l-1}^n\right]\right\|\le \frac{\sigma_{\hat{\nabla}}^2}{n_s M^2L^2}, \left\|\mc{E}\left[e_l^{(i)\top}e_l^{(i)}\biggm| \mathcal{F}_{l-1}^n\right]\right\|\le \frac{\sigma_{\hat{\nabla}}^2}{n_s M^2L^2},
\end{align*}
and based on this fact, we have
\begin{align*}
\left\|\mc{E}\left[e_l e_l^{\top}\biggm| \mathcal{F}_{l-1}^n\right]\right\| &=\left\|\mc{E}\left[\left(\sum_{i=1}^M e^{(i)}_l\right) \left(\sum_{i=1}^M e^{(i)\top}_l\right)\biggm| \mathcal{F}_{l-1}^n\right]\right\| \nonumber \\
&\le \sum_{i=1}^M \left\|\mc{E}\left[e_l^{(i)} e_l^{(i)\top}\biggm| \mathcal{F}_{l-1}^n\right]\right\| + \underbrace{\sum_{i\neq j}^M \left\|\mc{E}\left[e_l^{(i)} e_l^{(j)\top}\biggm| \mathcal{F}_{l-1}^n\right]\right\|}_{T_4=0}
\le \frac{\sigma_{\hat{\nabla}}^2}{n_s M L^2},
\end{align*}
where we use the fact that $T_4 =0$ because $e_l^{(i)}$ and $e_l^{(j)}$ are independent, if we conditioned on $\mathcal{F}_l^n.$ An identical argument holds for $\left\|\mc{E}\left[e_l^{\top}e_l\biggm| \mathcal{F}_{l-1}^n\right]\right\|$.

Define
$W_{\text{col},t}:=\sum_{l=0}^t \mc{E}\left[e_l e_l^{\top}\biggm| \mathcal{F}_{l-1}^n\right]$ and $W_{\text{row},t}:=\sum_{l=0}^t \mc{E}\left[e_l^{\top} e_l\biggm| \mathcal{F}_{l-1}^n\right]$, then we have $$\|W_{\text{col},t}\| \le \frac{\sigma_{\hat{\nabla}}^2}{n_s M L}, \quad \|W_{\text{row},t}\| \le \frac{\sigma_{\hat{\nabla}}^2}{n_s M L}.$$

\textbf{Claim III:} $\left\|e_l\right\| \le \frac{R_{\hat{\nabla}}}{n_s L}$ where $R_{\hat{\nabla}}=\frac{2 n_x n_u \bar{C}_{\max}}{r}+\frac{\epsilon}{2}+\bar{h}_1$.

\textbf{Proof:} From Lemma B.7 in~\citep{gravell2020learning}, we have $\|e_l^{(i)}\| \le \frac{n_s R_{\hat{\nabla}}}{ML}.$ With this fact, we have 
$$\left\| e_l\right\|\le \sum_{i=1}^M\left\| e^{(i)}_l\right\|\le \frac{R_{\hat{\nabla}}}{n_s L}.$$

With \textbf{Claim I, II} and \textbf{Claim III} and the matrix Freedman inequality~(\ref{lem:Freedman}), we have, for all $\epsilon \geq 0$,
\begin{align}\label{eq:freedman_inequality}
&\mathbb{P}\left\{\exists t \geq 0: \lambda_{\max }\left({Y}_t\right) \geq \epsilon  \text { and } \max \left\{\left\|{W}_{\text {col }, t}\right\|,\left\|{W}_{\text {row }, t}\right\|\right\}\le\frac{\sigma_{\hat{\nabla}}^2}{n_s M L} \right\} \notag\\
&\leq  (n_x +n_u)  \exp \left\{-\frac{\epsilon^2 / 2}{\frac{\sigma_{\hat{\nabla}}^2}{n_s M L}+\frac{R_{\hat{\nabla}}\epsilon}{3n_s L}}\right\}.
\end{align}
Therefore, rephrasing Eq.\eqref{eq:freedman_inequality}, if \begin{equation}\label{eq:sample_requirementt_variance}
n_s \ge \left(\underbrace{\frac{32\sigma_{\hat{\nabla}}^2 \min \left(n_x, n_u\right)}{ML\epsilon^2}}_{T_5}+\underbrace{\frac{32LR_{\hat{\nabla}}\sqrt{\min \left(n_x, n_u\right)}}{12L\epsilon}}_{T_6}\right) \log \left[\frac{ML(n_x+n_u)}{\delta}\right],\end{equation}
we have that
\begin{equation}
\begin{aligned}
\|Y_L\|_F =\|\frac{1}{ML}\sum_{i=1}^M \sum_{l=0}^{L-1} \left[ \hat{\nabla} C^{(i)}(K^{(i)}_{n,l})-\nabla C^{(i)}_r(K^{(i)}_{n,l})\right]\|_F \le \frac{\epsilon}{4},
\end{aligned}
\end{equation}
holds with probability $1-\delta$. As we discussed in Lemma~\ref{lemma:infinite}, we only keep the dominant term $T_5$ in the requirement of the sample size $n_s$ (as in Eq.\eqref{eq:sample_requirementt_variance}). Because $T_5$ is in the order $\mathcal{O}({\epsilon}^{-2})$ while $T_6$ is in the order $\mathcal{O}({\epsilon}^{-1})$. Moreover, we note that by making $\epsilon = \mathcal{O}(1/M)$, $T_6$ when compared to $T_5$ is negligible. We also note that, as an alternative to imposing the restriction $\epsilon = \mathcal{O}(1/M)$, one can apply the Freedman lemma (Lemma~\ref{lem:Freedman}) to $ML$ terms rather than to $L$ terms.

In summary, if $n_s \geq \frac{32\sigma_{\hat{\nabla}}^2 \min \left(n_x, n_u\right)}{ML\epsilon^2}\left[\frac{ML(n_x+n_u)}{\delta}\right]=\frac{h_{\text {sample }}\left(\frac{\epsilon}{4}, \frac{\delta}{ML} \right)}{ML} ,$
\begin{align}\label{eq:bound_T3_2}
\circled{2} = \Big\lVert \frac{1}{ML}\sum_{i=1}^M \sum_{l=0}^{L-1} \left[ \hat{\nabla} C^{(i)}(K^{(i)}_{n,l})-\nabla C^{(i)}_r(K^{(i)}_{n,l})\right]\Big\rVert_F \le \frac{\epsilon}{4}
\end{align} 
holds with probability $1-\delta.$

As a result, we have 
$
T_3 \le \frac{\epsilon}{2}
$
holds with probability $1-\delta,$ when $r \leq h_r\left(\frac{\epsilon}{4}\right)$ and $n_s \geq \frac{h_{\text {sample }}\left(\frac{\epsilon}{4}, \frac{\delta}{ML} \right)}{ML}.$ In what follows, we will provide an upper bound on the term $T_1$.

\paragraph{Bounding $T_1$:} We can follow the same analysis of bounding $\circled{2}$ in $T_3$ to bound $T_1.$ Different from the filtration we define in analyzing $\circled{2}$, we need to define a new filtration $\tilde{\mathcal{F}}_{l-1}^n$, where $\tilde{\mathcal{F}}_{l-1}^n:= \mathcal{F}_{l-1}^n \cup U_l^n$ and  $U_l^n:=\left\{U_{n,l,s}^{(i)}\right\}^{i=1,\cdots, N}_{s=1,\cdots,n_s}$. Note that $U_l^n$ is the sigma-field generated by the randomness of all random smoothing matrices $U_{n,l,s}^{(i)}$ \footnote{Here we use the index $s$ to denote $s$-th sample. Note that in each local iteration $l$, we need to generate the random smoothing matrices $n_s$ times.} from all the systems at the $n$-th global iteration and $l$-th local iteration. Replacing $\sigma_{\hat{\nabla}}^2$  Eq.\eqref{eq:sample_requirementt_variance} with into $\sigma_{\tilde{\nabla}}^2$ and $R_{\hat{\nabla}}$ with $R_{\tilde{\nabla}},$ we have that
\begin{align}\label{eq:bound_T1}
T_1= \Big\lVert\frac{1}{ML}\sum_{i=1}^M\sum_{l=0}^{L}\left[\tilde{\nabla} C^{(i)}(K^{(i)}_{n,l})-\nabla^{\prime} C^{(i)}(K^{(i)}_{n,l})\right] \Big\rVert_F \le \frac{\epsilon}{4}
\end{align}
holds with probability $1-\delta$ when 
\begin{align*}
n_s \geq \frac{32 \sigma_{\tilde{\nabla}}^2\min (n_x, n_u)}{ML\epsilon^2} \log \left[\frac{ML(n_x+n_u)}{\delta}\right]=\frac{h_{\text {sample,trunc }}\left(\frac{\epsilon}{4}, \frac{\delta}{ML}, \frac{H^2}{\mu}\right)}{ML}.
\end{align*}

Combing the upper bound of $T_1$ (Eq.\eqref{eq:bound_T1}), $T_2$ (Eq.\eqref{eq:bound_T2})  and $T_3$ (Eq.\eqref{eq:smoothing_error_variance} and \eqref{eq:bound_T3_2}), we have 
\begin{align}
\left\|\frac{1}{ML} \sum_{i=1}^M\sum_{l=0}^{L-1}\left[\widehat{\nabla C^{(i)}(K^{(i)}_{n,l})}-\nabla C^{(i)}(K^{(i)}_{n,l})\right]\right\|_F\notag \le T_1 +T_2 +T_3 \le \epsilon
\end{align}
when the trajectory length $\tau$  satisfies $\tau \geq h_{\tau}\left(\frac{r \epsilon}{4 n_x n_u}\right),$ the smoothing radius satisfies $r \leq h_r\left(\frac{\epsilon}{4}\right)$ and the size of samples satisfies $n_s \geq \max\left\{\frac{h_{\text {sample,trunc }}\left(\frac{\epsilon}{4}, \frac{\delta}{ML}, \frac{H^2}{\mu}\right)}{ML},\frac{h_{\text {sample }}\left(\frac{\epsilon}{4}, \frac{\delta}{ML} \right)}{ML} \right\}=\frac{h_{\text {sample,trunc }}\left(\frac{\epsilon}{4}, \frac{\delta}{ML}, \frac{H^2}{\mu}\right)}{ML}.$
Thus, we complete the proof of Lemma~\ref{lem:variance_reduction}.
\end{proof}

\subsubsection{Proof of Theorem~\ref{thm:main_fedlqr}}

\paragraph{Outline:} To prove Theorem~\ref{thm:main_fedlqr}, we first introduce some lemmas: Lemma \ref{lem:stability_local_model_free} establishes stability of the local policies; Lemma \ref{lem：drift_term_model_free} provides the drift analysis; Lemma \ref{lem:model_free_prp} quantifies the per-round progress of our \texttt{FedLQR} algorithm. As a result, we are able to present the iterative stability guarantees and convergence analysis of \texttt{FedLQR} in the model-free setting.

\begin{lemma} \label{lem:stability_local_model_free} (Stability of the local policies)  Suppose $K_n \in \mathcal{G}^0$ and the heterogeneity level satisfies $(\epsilon_1 \Bar{h}^1_{\text{het}} + \epsilon_2 \Bar{h}^2_{\text{het}})^2 \le\Bar{h}^3_{\text{het}}$, where $\Bar{h}^3_{\text{het}}$ is as defined in Eq.\eqref{eq:het_requirement}. If the local step-size $\eta_l$ satisfies $$\eta_l \le \min\left\{ \frac{\underline{h}_{\Delta} \mu}{H^2 \left(h_1 +\sqrt{\bar{\epsilon}}\right)},\frac{1}{9 \Bar{h}_{\text{grad}}} \right\},$$ the smoothing radius satisfies $$r\le \min\left\{\frac{\min_{i \in [M]} C^{(i)}(K_0)}{\Bar{h}_{\text{cost}}}, \underline{h}_{\Delta}, h_r\left(\frac{\sqrt{\bar{\epsilon}}}{4}\right)\right\},$$ the trajectory length satisfies $
\tau \geq h_{\tau}\left(\frac{r \sqrt{\bar{\epsilon}}}{4 n_x n_u}\right),$ and the number of the sample size satisfies $$n_{s} \geq \max \left\{h_{\text {sample,trunc }}\left(\frac{\sqrt{\bar{\epsilon}}}{4}, \frac{\delta}{L}, \frac{H^2}{\mu}\right), h_{\text {sample}}\left(\frac{\sqrt{\bar{\epsilon}}}{2}, \frac{\delta}{L}\right)\right\}$$ where we choose a fixed error tolerance $\bar{\epsilon}$ to be $$\bar{\epsilon}:=\min_{j \in [M]} \left\{\frac{3 \mu^2 \sigma_{\min }(R)\left(C^{(j)}(K_0) - C^{(j)}(K_j^*)\right)}{5||\Sigma_{K_j^*}||}\right\},$$   then with probability $1-\delta,$, where $\delta \in (0,1)$, $K_{n,l}^{(i)} \in \mathcal{G}^0$ holds for all $i \in [M]$ and $l = 0,1, \cdots, L-1$.
\end{lemma}

\begin{proof}
For any $i, j \in [M],$ according to the local Lipschitz property in Lemma~\ref{Lemma:lipschitz}, we have that 
\begin{align*}%\label{eq:stability_decrease}
C^{(j)}(K_{n,1}^{(i)}) - C^{(j)}(K_n)&\le \left\langle \nabla C^{(j)}(K_n), K_{n,1}^{(i)} -K_n   \right\rangle\\
&+ \frac{h_{\text{grad}(K_n)}}{2}\Big\lVert K_{n,1}^{(i)} -K_n \Big\rVert_F^2 \quad (\text{Local lipschitz})\\
& = -\left\langle \nabla C^{(j)}(K_n), \eta_l\tilde{\nabla} C^{(i)}(K_n) \right\rangle 
+ \frac{h_{\text{grad}(K_n)}}{2}\Big\lVert \eta_l\tilde{\nabla} C^{(i)}(K_n) \Big\rVert_F^2,
\end{align*}
holds if $\Big\lVert \eta_l\tilde{\nabla} C^{(i)}(K_n) \Big\rVert_F  \le \underline{h}_{\Delta} \le h_{\Delta}(K_n)$. Note that this inequality holds when $\eta_l$ satisfies
\begin{align}\label{eq:step_model_free}
\Big\lVert \eta_l\tilde{\nabla} C^{(i)}(K_n) \Big\rVert_F & =\eta_l \Big\lVert \frac{1}{n_s}\sum_{s=1}^{n_s}\frac{n_u n_x}{r^2} \tilde{C}^{(i),(\tau)}\left(K_n+U^{(i)}_{s}\right) U^{(i)}_{s}\Big\rVert_F\notag\\
&\stackrel{(a)}{\le}\eta_l  \frac{H^2 }{\mu}\Big\lVert  \frac{1}{n_{s}} \sum_{i=1}^{n_{s}} \frac{n_x n_u}{r^2} C^{(i)}\left(K_n+U^{(i)}_{s}\right)U^{(i)}_{s} \Big\rVert_F\notag\\
& = \frac{\eta_l H^2}{\mu} \Big\lVert\hat{\nabla} C^{(i)}(K)\Big\rVert_F\notag\\
& \le  \frac{\eta_l H^2}{\mu} \left[\Big\lVert \nabla C^{(i)}(K)\Big\rVert_F+\Big\lVert \hat{\nabla} C^{(i)}(K)- \nabla C^{(i)}(K)\Big\rVert_F \right]\notag\\
& \stackrel{(b)}{\le} \frac{\eta_l H^2}{\mu} \left[\Big\lVert \nabla C^{(i)}\left(K_n\right) \Big\rVert_F +\sqrt{\bar{\epsilon}}\right]\notag\\
& \le  \frac{\eta_l H^2}{\mu}\left(h_1 +\sqrt{\bar{\epsilon}}\right)
\end{align}

where\footnote{For sake of the notation, we ignore the dependence on the local iteration $l$ and global iteration $n$ when we index $U_s^{(i)}$ in this part.} (a) is due to Lemma~\ref{lemma:finitely_many_finite_horizon_estimation}; according to Lemma~\ref{lemma:infinite}, (b) holds with high probability, when the number of the sample size satisfies $n_{s} \geq  h_{\text {sample}}\left(\frac{\sqrt{\bar{\epsilon}}}{2}, \frac{\delta}{L}\right)$. The last inequality follows from the uniform upper gradient bound in Lemma~\ref{lemma:uniform_bounds}. Then we can easily conclude that $\Big\lVert \eta_l\tilde{\nabla} C^{(i)}(K_n) \Big\rVert_F  \le \underline{h}_{\Delta}$ holds when $\eta_l \le \frac{\underline{h}_{\Delta} \mu}{H^2 \left(h_1 +\sqrt{\bar{\epsilon}}\right)}.$

Following the analysis in Eq~\eqref{eq:stability_decrease},  we have
\begin{align*}%\label{eq:stability_decrease_II}
C^{(j)}(K_{n,1}^{(i)}) &- C^{(j)}(K_n)\le - \eta_l\left\langle \nabla C^{(j)}(K_n), \nabla C^{(j)}(K_n) \right\rangle \notag \\ 
&- \eta_l\underbrace{\left\langle \nabla C^{(j)}(K_n), \nabla C^{(i)}(K_n)-\nabla C^{(j)}(K_n)  \right\rangle}_{T_{1}} \notag\\
& \underbrace{-\eta_l\left\langle \nabla C^{(j)}(K_n), \tilde{\nabla} C^{(i)}(K_n)- \nabla C^{(i)}(K_n) \right\rangle}_{T_{2}}+ \frac{h_{\text{grad}(K_n)}}{2}\Big\lVert \eta_l\tilde{\nabla} C^{(i)}(K_n) \Big\rVert_F^2,
\end{align*}
where $T_{1}$ can be upper bounded as
\begin{align*}
T_{1} &\le \eta_l \Big\lVert \nabla C^{(j)}(K_n)\Big\rVert_F \Big\lVert \nabla C^{(i)}(K_n)- \nabla C^{(j)}(K_n) \Big\rVert_F\notag\\
&\le  \eta_l \sqrt{\min\{n_x, n_u\}} \Big\lVert \nabla C^{(j)}(K_n)\Big\rVert_F (\epsilon_1\bar{h}^{1}_{\text{het}} +\epsilon_2\Bar{h}^{2}_{\text{het}}),
\end{align*}
where we use the policy gradient heterogeneity bound in Lemma~\ref{lem:gradient_het} and the fact that $K_n \in \mathcal{G}^0$.

We can bound $T_2$ as follows  
\begin{align*}
T_{2} &\le \eta_l \Big\lVert \nabla C^{(j)}(K_n)\Big\rVert_F \Big\lVert \tilde{\nabla} C^{(i)}(K_n)- \nabla C^{(i)}(K_n) \Big\rVert_F\notag\\
&\le  \eta_l \Big\lVert \nabla C^{(j)}(K_n)\Big\rVert_F \sqrt{\bar{\epsilon}},
\end{align*}
where it holds with probability $1-\delta$. Here we use the Cauchy-Schwarz inequality in the first inequality, and the second inequality is due to Lemma~\ref{lemma:finitely_many_finite_horizon_estimation} since $n_{s} \geq h_{\text {sample,trunc }}\left(\frac{\sqrt{\bar{\epsilon}}}{4}, \delta, \frac{H^2}{\mu}\right)$, the smoothing radius satisfies $r \leq h_r\left(\frac{\sqrt{\bar{\epsilon}}}{4}\right)$ and the length of trajectories satisfies $
\tau \geq h_{\tau}\left(\frac{r \sqrt{\bar{\epsilon}}}{4 n_x n_u}\right).$ 

Plugging the upper bounds of $T_1$ and $T_2$ in Eq~\eqref{eq:stability_decrease_II}, we have:
\begin{align*}%\label{eq:stability_decrease_part3}
C^{(j)}(K_{n,1}^{(i)})&- C^{(j)}(K_n)\stackrel{(a)}{\le}\\
&- \eta_l \Big\lVert  \nabla C^{(j)}(K_n)\Big\rVert_F^2 + \eta_l \sqrt{\min\{n_x, n_u\}} \Big\lVert \nabla C^{(j)}(K_n)\Big\rVert_F (\epsilon_1\bar{h}^{1}_{\text{het}} +\epsilon_2\Bar{h}^{2}_{\text{het}})\notag\\
&+\eta_l \Big\lVert \nabla C^{(j)}(K_n)\Big\rVert_F\sqrt{\bar{\epsilon}} + \frac{3h_{\text{grad}(K_n)}\eta_l^2}{2}\Big\lVert \tilde{\nabla} C^{(i)}(K_n) -\nabla C^{(i)}(K_n) \Big\rVert_F^2\notag\\
 & + \frac{3h_{\text{grad}(K_n)}\eta_l^2}{2} \Big\lVert \nabla C^{(i)}(K_n)-\nabla C^{(j)}(K_n)\Big\rVert_F^2 \notag\\
 &+\frac{3h_{\text{grad}(K_n)}\eta_l^2}{2}\Big\lVert \nabla C^{(j)}(K_n)\Big\rVert_F^2\notag\\
  & \stackrel{(b)}{\le} - \eta_l \Big\lVert  \nabla C^{(j)}(K_n)\Big\rVert_F^2 + \eta_l \sqrt{\min\{n_x, n_u\}}\Big\lVert \nabla C^{(j)}(K_n)\Big\rVert_F (\epsilon_1 \Bar{h}^{1}_{\text{het}}+\epsilon_2 \Bar{h}^{2}_{\text{het}})\notag\\
 & +\eta_l \Big\lVert \nabla C^{(j)}(K_n)\Big\rVert_F \sqrt{\bar{\epsilon}} + \frac{3\Bar{h}_{\text{grad}}\eta_l^2}{2} \bar{\epsilon}\notag\\
 &+ \frac{3\Bar{h}_{\text{grad}}\eta_l^2\min\{n_x, n_u\}}{2}(\epsilon_1\bar{h}^{1}_{\text{het}} +\epsilon_2\Bar{h}^{2}_{\text{het}})^2 +\frac{3\Bar{h}_{\text{grad}}\eta_l^2}{2}\Big\lVert \nabla C^{(j)}(K_n)\Big\rVert_F^2,
\end{align*}
where $(a)$ follows from Eq.\eqref{eq:sum_expand}; $(b)$ follows from the same reasoning as we bound $T_{1}$ and $T_{2}$ and the fact that $K_n \in \mathcal{G}^0$. If we choose the local step-size $\eta_l$ satisfies $\eta_l \le \frac{1}{9 \Bar{h}_{\text{grad}}}$, i.e., $\frac{3\Bar{h}_{\text{grad}}\eta_l^2}{2} \le \frac{\eta_l}{6}$, we have
\begin{align*}%\label{eq:stability_decrease_part3}
&C^{(j)}(K_{n,1}^{(i)}) - C^{(j)}(K_n) \stackrel{(a)}{\le} - \eta_l \Big\lVert  \nabla C^{(j)}(K_n)\Big\rVert_F^2 + \frac{\eta_l}{6} \Big\lVert \nabla C^{(j)}(K_n)\Big\rVert_F^2 \notag\\
 &+\frac{3\eta_l\min\{n_x, n_u\}}{2}(\epsilon_1 \Bar{h}^{1}_{\text{het}}+\epsilon_2 \Bar{h}^{2}_{\text{het}})^2 +\frac{\eta_l \Big\lVert \nabla C^{(j)}(K_n)\Big\rVert_F^2}{6}+ \frac{3\eta_l \bar{\epsilon}}{2} + \frac{\eta_l }{6} \bar{\epsilon}\notag\\
 &+ \frac{\eta_l\min\{n_x, n_u\}}{6}(\epsilon_1\bar{h}^{1}_{\text{het}} +\epsilon_2\Bar{h}^{2}_{\text{het}})^2 +\frac{\eta_l}{6}\Big\lVert \nabla C^{(j)}(K_n)\Big\rVert_F^2 \notag\\
 & \le - \frac{\eta_l}{2} \Big\lVert  \nabla C^{(j)}(K_n)\Big\rVert_F^2 +\frac{5\eta_l\min\{n_x, n_u\}}{3}(\epsilon_1 \Bar{h}^{1}_{\text{het}}+\epsilon_2 \Bar{h}^{2}_{\text{het}})^2 + \frac{5\eta_l}{3} \bar{\epsilon}\notag\\
 & \stackrel{(b)}{\le}  -\frac{2\eta_l\sigma_{\text{min}}(R)\mu^2}{||\Sigma_{K_j^*}||} (C^{(j)}(K_n) - C^{(j)}(K_j^*))  +\frac{5\eta_l\min\{n_x, n_u\}}{3}(\epsilon_1 \Bar{h}^{1}_{\text{het}}+\epsilon_2 \Bar{h}^{2}_{\text{het}})^2 + \frac{5\eta_l}{3} \bar{\epsilon},
\end{align*}
where $(a)$ follows from the Young's inequality in Eq.\eqref{eq:vector_youngs}; and $(b)$ follows from the gradient domination in Lemma~\ref{lemma:gradient_domination}. 

Therefore, if the heterogeneity satisfies $(\epsilon_1 \Bar{h}^1_{\text{het}} + \epsilon_2 \Bar{h}^2_{\text{het}})^2 \le\Bar{h}^3_{\text{het}}$, then we have $$(\epsilon_1 \Bar{h}^1_{\text{het}} + \epsilon_2 \Bar{h}^2_{\text{het}})^2 \le \min_{j \in [M]} \left\{\frac{3 \mu^2 \sigma_{\min }(R)\left(C^{(j)}(K_0) - C^{(j)}(K_j^*)\right)}{5||\Sigma_{K_j^*}||\min\{n_x, n_u\}}\right\}.$$

Since the error tolerance $$\bar{\epsilon}=\min_{j \in [M]} \left\{\frac{3 \mu^2 \sigma_{\min }(R)\left(C^{(j)}(K_0) - C^{(j)}(K_j^*)\right)}{5||\Sigma_{K_j^*}||}\right\},$$ we have 
\begin{align*}
&C^{(j)}(K_{n,1}^{(i)}) - C^{(j)}(K_j^*) \le  \left(1-\frac{2\eta \mu^2 \sigma_{\min }(R)}{\left\|\Sigma_{K_j^*}\right\|}\right) (C^{(j)}(K_n) - C^{(j)}(K_j^*))\notag\\
&+\frac{\eta_l \mu^2 \sigma_{\min }(R)\left(C^{(j)}(K_0) - C^{(j)}(K_j^*)\right)}{||\Sigma_{K_j^*}||} +\frac{\eta_l\mu^2 \sigma_{\min }(R)\left(C^{(j)}(K_0) - C^{(j)}(K_j^*)\right)}{||\Sigma_{K_j^*}||} \notag\\
&\stackrel{(a)}{\le}  \left(1-\frac{2\eta \mu^2 \sigma_{\min }(R)}{\left\|\Sigma_{K_j^*}\right\|}\right) (C^{(j)}(K_0) - C^{(j)}(K_j^*))  +\frac{2\eta_l\mu^2 \sigma_{\min }(R)\left(C^{(j)}(K_0) - C^{(j)}(K_j^*)\right)}{||\Sigma_{K_j^*}||}\notag\\
&= C^{(j)}(K_0) - C^{(j)}(K_j^*), \forall j \in [M],
\end{align*}
where we use the fact that $K_n \in \mathcal{G}^0$ in $(a)$. The above inequality implies $K_{n,1}^{(i)}\in \mathcal{G}^0$ with high probability $1-\delta$ when $K_n \in \mathcal{G}^0.$ Then we can use the induction method to obtain that $K_{n,2}^{(i)}\in \mathcal{G}^0$, since $K_{n,1}^{(i)}\in \mathcal{G}^0$. By repeating this step for $L$ times, we have that all the local polices $K_{n,l}^{(i)} \in \mathcal{G}^0$ holds for all $i \in [M]$ and $l = 0,1, \cdots, L-1$, when the global policy $K_n \in \mathcal{G}^0.$ 
\end{proof}

\begin{lemma}\label{lem：drift_term_model_free}
(Drift term analysis)  Suppose $K_n \in \mathcal{G}^0.$ If $\eta_l \le \min \left\{\frac{1}{4\Bar{h}_{\text {grad}}}, \frac{1}{4},\frac{\log 2}{L(3\Bar{h}_{\text {grad}}+2)}\right\},$ the number of the sample size $n_s$ satisfies $$n_{s} \geq \frac{h_{\text {sample,trunc }}\left(\frac{\sqrt{\epsilon}}{4}, \frac{\delta}{L}, \frac{H^2}{\mu}\right)}{ML},$$ the smoothing radius satisfies $r \leq h_r\left(\frac{\sqrt{\epsilon}}{4}\right)$ and the length of trajectories satisfies $$
\tau \geq h_{\tau}\left(\frac{r \sqrt{\epsilon}}{4 n_x n_u}\right),$$ given any $\delta \in (0,1),$ the difference between the local policy and global policy can be bounded by
\begin{align*}
\Big\lVert K_{n, l}^{(i)}-K_{n}\Big\rVert_F^2 \le 2\eta_l L \left[ \Big\lVert\nabla C^{(i)}(K_{n})\Big\rVert_F^2 +ML\epsilon \right] = \frac{2\eta}{\eta_g} \left[ \Big\lVert\nabla C^{(i)}(K_{n})\Big\rVert_F^2 +ML\epsilon \right]
\end{align*}
holds, with probability $1-\delta$, for all $i \in [M]$ and $l =0,1,\cdots, L-1.$
\end{lemma}
\begin{proof} \begin{align}\label{eq:model_free_drift_bound}
&\Big\lVert K_{n, l}^{(i)}-K_{n}\Big\rVert_F^2 = \Big\lVert K^{(i)}_{n,l-1}-K_{n} - \eta_l\tilde{\nabla} C^{(i)}(K^{(i)}_{n,l-1})\Big\rVert_F^2 \notag\\
& = \Big\lVert K^{(i)}_{n,l-1}-K_{n} \Big\rVert_F^2 -2\eta_l  \left[\left\langle \tilde{\nabla} {C^{(i)}(K^{(i)}_{n,l-1})}, K^{(i)}_{n,l-1}-K_{n}\right\rangle\right] \notag\\
& + \Big\lVert \eta_l\tilde{\nabla} C^{(i)}(K^{(i)}_{n,l-1})\Big\rVert_F^2 \notag\\
&= \Big\lVert K^{(i)}_{n,l-1}-K_{n} \Big\rVert_F^2 -2\eta_l  \left[\left\langle \tilde{\nabla} {C^{(i)}(K^{(i)}_{n,l-1})} -\nabla C^{(i)}(K^{(i)}_{n,l-1}), K^{(i)}_{n,l-1}-K_{n}\right\rangle\right] \notag\\
& -2\eta_l \left[\left\langle \nabla C^{(i)}(K^{(i)}_{n,l-1})-\nabla C^{(i)}(K_{n}), K^{(i)}_{n,l-1}-K_{n}\right\rangle\right] - 2\eta_l  \left[\left\langle\nabla C^{(i)}(K_{n}), K^{(i)}_{n,l-1}-K_{n}\right\rangle\right]\notag\\
&+ \Big\lVert \eta_l\tilde\nabla {C^{(i)}(K^{(i)}_{n,l-1})}\Big\rVert_F^2 \notag\\
&\stackrel{(a)}{\le} \Big\lVert K^{(i)}_{n,l-1}-K_{n} \Big\rVert_F^2 -2\eta_l  \left[\left\langle \tilde{\nabla} {C^{(i)}(K^{(i)}_{n,l-1})} -\nabla C^{(i)}(K^{(i)}_{n,l-1}), K^{(i)}_{n,l-1}-K_{n}\right\rangle\right] \notag\\
& +2\eta_l \Big\lVert\nabla C^{(i)}(K^{(i)}_{n,l-1})-\nabla C^{(i)}(K_{n})\Big\rVert_F \Big\lVert K^{(i)}_{n,l-1}-K_{n}\Big\rVert_F+ 2\eta_l  \Big\lVert\nabla C^{(i)}(K_{n})\Big\rVert_F\Big\lVert K^{(i)}_{n,l-1}-K_{n}\Big\rVert_F\notag\\
&+ \Big\lVert \eta_l\tilde{\nabla} {C^{(i)}(K^{(i)}_{n,l-1})}\Big\rVert_F^2\notag\\
& \stackrel{(b)}{\le} \Big\lVert K^{(i)}_{n,l-1}-K_{n} \Big\rVert_F^2 -2\eta_l  \left[\left\langle \tilde{\nabla} {C^{(i)}(K^{(i)}_{n,l-1})} -\nabla C^{(i)}(K^{(i)}_{n,l-1}), K^{(i)}_{n,l-1}-K_{n}\right\rangle\right] \notag\\
& +2\eta_l h_{\text {grad}}(K_n) \Big\lVert K^{(i)}_{n,l-1}- K_{n}\Big\rVert_F \Big\lVert K^{(i)}_{n,l-1}-K_{n}\Big\rVert_F+ \eta_l  \Big\lVert\nabla C^{(i)}(K_{n})\Big\rVert_F^2 +\eta_l\Big\lVert K^{(i)}_{n,l-1}-K_{n}\Big\rVert_F^2\notag\\
&+ \Big\lVert \eta_l\tilde{\nabla} {C^{(i)}(K^{(i)}_{n,l-1})}\Big\rVert_F^2
\end{align}
where we use Cauchy–schwarz inequality for $(a)$; and for (b), we use Eq.~\eqref{eq:vector_youngs}.

Following the analysis in Eq.\eqref{eq:model_free_drift_bound}, we have
\begin{align*}
&\Big\lVert K_{n, l}^{(i)}-K_{n}\Big\rVert_F^2 \le \left(1+2\eta_l h_{\text {grad}}(K_n)+  \eta_l\right)\Big\lVert K^{(i)}_{n,l-1}-K_{n} \Big\rVert_F^2 +\eta_l  \Big\lVert\nabla C^{(i)}(K_{n})\Big\rVert_F^2 \notag\\
&-2\eta_l  \left[\left\langle \tilde{\nabla} {C^{(i)}(K^{(i)}_{n,l-1})} -\nabla C^{(i)}(K^{(i)}_{n,l-1}), K^{(i)}_{n,l-1}-K_{n}\right\rangle\right] + \Big\lVert \eta_l\tilde{\nabla} {C^{(i)}(K^{(i)}_{n,l-1})}\Big\rVert_F^2 \notag\\
&\stackrel{(a)}{\le} \left(1+2\eta_l h_{\text {grad}}(K_n)+  \eta_l\right)\Big\lVert K^{(i)}_{n,l-1}-K_{n} \Big\rVert_F^2 +\eta_l  \Big\lVert\nabla C^{(i)}(K_{n})\Big\rVert_F^2 \notag\\
&+2\eta_l  \left[\Big\lVert \tilde{\nabla} {C^{(i)}(K^{(i)}_{n,l-1})} -\nabla C^{(i)}(K^{(i)}_{n,l-1})\Big\rVert_F \Big\lVert K^{(i)}_{n,l-1}-K_{n}\Big\rVert_F\right] \notag\\
&+ 2 \eta_l^2\Big\lVert \tilde{\nabla} {C^{(i)}(K^{(i)}_{n,l-1})} -\nabla C^{(i)}(K^{(i)}_{n,l-1})\Big\rVert_F^2 + 2 \eta_l^2\Big\lVert \nabla C^{(i)}(K^{(i)}_{n,l-1})\Big\rVert_F^2 \notag\\
& \stackrel{(b)}{\le} \left(1+2\eta_l h_{\text {grad}}(K_n)+  \eta_l\right)\Big\lVert K^{(i)}_{n,l-1}-K_{n} \Big\rVert_F^2 +\eta_l  \Big\lVert\nabla C^{(i)}(K_{n})\Big\rVert_F^2 \notag\\
&+\eta_l  \Big\lVert \tilde{\nabla} {C^{(i)}(K^{(i)}_{n,l-1})} -\nabla C^{(i)}(K^{(i)}_{n,l-1})\Big\rVert_F^2 + \eta_l  \Big\lVert K^{(i)}_{n,l-1}-K_{n}\Big\rVert_F^2 \notag\\
&+ 2 \eta_l^2\Big\lVert \tilde{\nabla} {C^{(i)}(K^{(i)}_{n,l-1})} -\nabla C^{(i)}(K^{(i)}_{n,l-1})\Big\rVert_F^2\notag\\
&+ 4 \eta_l^2\Big\lVert \nabla C^{(i)}(K^{(i)}_{n,l-1}) - \nabla C^{(i)}(K_{n})\Big\rVert_F^2 + 4 \eta_l^2\Big\lVert \nabla C^{(i)}(K_{n})\Big\rVert_F^2 \notag\\
& \stackrel{(c)}{\le}  \left(1+2\eta_l h_{\text {grad}}(K_n)+ \eta_l\right)\Big\lVert K^{(i)}_{n,l-1}-K_{n} \Big\rVert_F^2 +(\eta_l + 4 \eta_l^2) \Big\lVert\nabla C^{(i)}(K_{n})\Big\rVert_F^2 \notag\\
&+\eta_l  \Big\lVert \tilde{\nabla} {C^{(i)}(K^{(i)}_{n,l-1})} -\nabla C^{(i)}(K^{(i)}_{n,l-1})\Big\rVert_F^2 + \eta_l  \Big\lVert K^{(i)}_{n,l-1}-K_{n}\Big\rVert_F^2 \notag\\
&+ 2 \eta_l^2\Big\lVert \tilde{\nabla} {C^{(i)}(K^{(i)}_{n,l-1})} -\nabla C^{(i)}(K^{(i)}_{n,l-1})\Big\rVert_F^2
+ 4 \eta_l^2 h_{\text {grad}}(K_n)^2\Big\lVert K^{(i)}_{n,l-1} -  K_{n}\Big\rVert_F^2  \notag\\
& \stackrel{(d)}{=} \left(1+2\eta_l h_{\text {grad}}(K_n)+  2\eta_l+4 \eta_l^2 h_{\text {grad}}(K_n)^2\right)\Big\lVert K^{(i)}_{n,l-1}-K_{n} \Big\rVert_F^2 +(\eta_l + 4 \eta_l^2) \Big\lVert\nabla C^{(i)}(K_{n})\Big\rVert_F^2 \notag\\
&+\left(\eta_l + 2 \eta_l^2\right) \Big\lVert \tilde{\nabla} {C^{(i)}(K^{(i)}_{n,l-1})} -\nabla C^{(i)}(K^{(i)}_{n,l-1})\Big\rVert_F^2\notag\\
& \le \left(1+2\eta_l \bar{h}_{\text {grad}}+  2\eta_l+4 \eta_l^2 \bar{h}^2_{\text {grad}}\right)\Big\lVert K^{(i)}_{n,l-1}-K_{n} \Big\rVert_F^2 +(\eta_l + 4 \eta_l^2) \Big\lVert\nabla C^{(i)}(K_{n})\Big\rVert_F^2 \notag\\
&+\left(\eta_l + 2 \eta_l^2\right) \underbrace{\Big\lVert \tilde{\nabla} {C^{(i)}(K^{(i)}_{n,l-1})} -\nabla C^{(i)}(K^{(i)}_{n,l-1})\Big\rVert_F^2}_{T_1},%\label{eq:drift_term}
\end{align*}
where we use Cauchy-Schwarz inequality and Eq.\eqref{eq:youngs} for $(a)$; for $(b)$, we use Eq.\eqref{eq:youngs} and \eqref{eq:sum_expand}; for $(c)$, we use the gradient smoothness lemma in Lemma~\ref{Lemma:lipschitz}; and for $(d)$, we use the fact that $K_n \in \mathcal{G}^0$.

From Lemma~\ref{lemma:finitely_many_finite_horizon_estimation}, we can bound $T_1$ term as follows:
\begin{align*}
T_1 = \Big\lVert \tilde{\nabla} {C^{(i)}(K^{(i)}_{n,l-1})} -\nabla C^{(i)}(K^{(i)}_{n,l-1})\Big\rVert_F^2 \le ML\epsilon,
\end{align*}
where it holds with probability $1-\delta$, since $n_{s} \geq \frac{h_{\text {sample,trunc }}\left(\frac{\sqrt{\epsilon}}{4}, \delta, \frac{H^2}{\mu}\right)}{ML}$, the smoothing radius satisfies $r \leq h_r\left(\frac{\sqrt{\epsilon}}{4}\right)$ and the length of trajectories satisfies $
\tau \geq h_{\tau}\left(\frac{r \sqrt{\epsilon}}{4 n_x n_u}\right).$

Hence, we obtain
\begin{align*}
\Big\lVert K_{n, l}^{(i)}-K_{n}\Big\rVert_F^2 & \le \left(1+2\eta_l \Bar{h}_{\text {grad}}+  2\eta_l+4 \eta_l^2 \Bar{h}_{\text {grad}}^2\right)\Big\lVert K^{(i)}_{n,l-1}-K_{n} \Big\rVert_F^2 +(\eta_l + 4 \eta_l^2) \Big\lVert\nabla C^{(i)}(K_{n})\Big\rVert_F^2 \notag\\
&+\left(\eta_l + 2 \eta_l^2\right)ML\epsilon \notag\\
&\stackrel{(a)}{\le} (1+ 3\eta_l\Bar{h}_{\text {grad}} +2\eta_l)\Big\lVert K^{(i)}_{n,l-1}-K_{n} \Big\rVert_F^2+2\eta_l \Big\lVert\nabla C^{(i)}(K_{n})\Big\rVert_F^2  + 2\eta_l ML \epsilon \notag\\
& \le (1+ 3\eta_l\Bar{h}_{\text {grad}} +2\eta_l)^{l}\underbrace{\Big\lVert K^{(i)}_{n,0}-K_{n} \Big\rVert_F^2}_{= 0}\notag\\
&+2 \eta_l\sum_{j=0}^{l-1}\left(1+ 3\eta_l\Bar{h}_{\text {grad}} +2\eta_l \right)^j\left[ \Big\lVert\nabla C^{(i)}(K_{n})\Big\rVert_F^2 +ML\epsilon \right] \notag\\
& \le 2\eta_l\times\frac{\left(1+ 3\eta_l\Bar{h}_{\text {grad}} +2\eta_l \right)^l -1}{\left(1+ 3\eta_l\Bar{h}_{\text {grad}} +2\eta_l \right)-1} \left[ \Big\lVert\nabla C^{(i)}(K_{n})\Big\rVert_F^2 +ML\epsilon \right] \notag\\
& \le 2\times\frac{\left(1+ 3\eta_l\Bar{h}_{\text {grad}} +2\eta_l \right)^l -1}{3\Bar{h}_{\text {grad}} +2}\left[ \Big\lVert\nabla C^{(i)}(K_{n})\Big\rVert_F^2 +ML\epsilon \right]\notag\\
& \stackrel{(b)}{\le} 2\times\frac{1+ l(3\eta_l\Bar{h}_{\text {grad}} +2\eta_l)-1}{3\Bar{h}_{\text {grad}} +2}\left[ \Big\lVert\nabla C^{(i)}(K_{n})\Big\rVert_F^2 +ML\epsilon \right] \notag\\
& \le 2\eta_l L \left[ \Big\lVert\nabla C^{(i)}(K_{n})\Big\rVert_F^2 +ML\epsilon \right],
\end{align*}
where $(a)$ is due to the choice of local step-size which satisfies $2\eta_l \Bar{h}_{\text {grad}} +2\eta_l + 4\eta_l^2\Bar{h}^2_{\text {grad}} \le 3\eta_l\Bar{h}_{\text {grad}} +2\eta_l$ and $\eta_l + 2 \eta_l^2 \le \eta_l + 4 \eta_l^2 \le 2\eta_l$, i.e., $\eta_l \le \min \left\{\frac{1}{4\Bar{h}_{\text {grad}}}, \frac{1}{4}\right\}.$ For $(b)$, we used the fact that $(1+x)^{\tau+1} \leq 1+2 x(\tau+1)$ holds for $x \leq \frac{\log 2}{\tau}$. In other words, $\left(1+ 3\eta_l\Bar{h}_{\text {grad}} +2\eta_l \right)^l \le 1+ l(3\eta_l\Bar{h}_{\text {grad}} +2\eta_l)$ when $3\eta_l\Bar{h}_{\text {grad}} +2\eta_l \le \frac{\log 2}{l},$ i.e., $\eta_l \le \frac{\log 2}{L(3\Bar{h}_{\text {grad}}+2)}.$ 
\end{proof}

\begin{lemma}\label{lem:model_free_prp}
(Per round progress) Suppose $K_n \in\mathcal{G}^0.$ If we choose the local step-size as $$\eta_l =\frac{1}{2} \min\left\{\frac{\underline{h}_{\Delta} \mu}{H^2 \left(h_1 +\sqrt{\epsilon}\right)},\frac{1}{9 \Bar{h}_{\text{grad}}}, \frac{1}{4},\frac{\log 2}{L(3\Bar{h}_{\text {grad}}+2)},\frac{1}{256L\bar{h}^2_{\text {grad}}} \right\},$$ 
with step-size $\eta:= L\eta_l \eta_g = \frac{1}{2}\min\{\frac{\underline{h}_{\Delta} \mu}{H^2 \left(h_1 +\sqrt{{\epsilon}}\right)},1,\frac{1}{32\Bar{h}_{\text {grad}}} \},$ and the smoothing radius\footnote{The exact requirement of $r$ is $r \le \min\left\{\frac{\min_{i \in [M]} C^{(i)}(K_0)}{\Bar{h}_{\text{cost}}}, \underline{h}_{\Delta}, h_r\left(\frac{\sqrt{\epsilon}}{4}\right), h_r\left(\frac{\sqrt{\bar{\epsilon}}}{4}\right)\right\}$. Here, without loss of generality, we drop the  $h_r\left(\frac{\sqrt{\bar{\epsilon}}}{4}\right)$ term from the $\min$ expression. This can be done because the error tolerance $\epsilon$ is usually small, and so $h_r\left(\frac{\sqrt{\epsilon}}{4}\right) \le h_r\left(\frac{\sqrt{\bar{\epsilon}}}{4}\right)$ holds. The assumptions on $\tau$ and $n_s$ follow similarly.} $$r\le \min\left\{\frac{\min_{i \in [M]} C^{(i)}(K_0)}{\Bar{h}_{\text{cost}}}, \underline{h}_{\Delta}, h_r\left(\frac{\sqrt{\epsilon}}{4}\right)\right\},$$ where the trajectory length satisfies $
\tau \geq h_{\tau}\left(\frac{r \sqrt{\epsilon}}{4 n_x n_u}\right),$ and the number of the sample size satisfies $$n_{s} \geq \frac{h_{\text {sample,trunc }}\left(\frac{\sqrt{\epsilon}}{4}, \frac{\delta}{L}, \frac{H^2}{\mu}\right)}{ML},$$ then with probability $1-\delta$, for any small $\delta \in (0,1)$, the  \texttt{FedLQR} algorithm provides the following convergence guarantee: 
\begin{align}\label{eq:one_step_convergence}
C^{(i)}(K_{n+1}) - C^{(i)}(K_i^*)&\le \left(1-\frac{\eta \mu^2 \sigma_{\min }(R)}{\left\|\Sigma_{K_i^*}\right\|}\right) (C^{(i)}(K_n) - C^{(i)}(K_i^*)) + 2\eta \epsilon\notag\\
&+ 2\eta\min\{n_x,n_u\}(\epsilon_1 \bar{h}^1_{\text{het}} + \epsilon_1 \bar{h}^2_{\text{het}}).^2
\end{align}\end{lemma}
\begin{proof} For any $i \in [M],$ according to the local Lipschitz property in Lemma~\ref{Lemma:lipschitz}, we have that 
\begin{align}\label{eq:model_free_prp}
&C^{(i)}(K_{n+1}) - C^{(i)}(K_n) \le \langle \nabla C^{(i)}(K_n), K_{n+1}-K_n \rangle +\frac{h_{\text {grad}}(K_n)}{2}\lVert K_{n+1}-K_n\rVert_F^2 \notag\\
& = -\left\langle\nabla C^{(i)}(K_n), \frac{\eta}{ML}\sum_{j=1}^M \sum_{l=0}^{L-1} \tilde{\nabla} C^{(j)}(K_{n, l}^{(j)})\right\rangle +\frac{h_{\text {grad}}(K_n)}{2} \Big\lVert \frac{\eta}{ML}\sum_{j=1}^M \sum_{l=0}^{L-1}\tilde{\nabla} C^{(j)}(K_{n, l}^{(j)})\Big\rVert_F^2,
\end{align}
holds when $\left\| \frac{\eta}{ML}\sum_{j=1}^M \sum_{l=0}^{L-1} \tilde{\nabla} C^{(j)}(K_{n, l}^{(j)})\right\|_F \le \underline{h}_{\Delta} \le h_{\Delta}(K_n)$. Following the same analysis as Eq.\eqref{eq:step_model_free}, this inequality holds when $$\eta \le \frac{\underline{h}_{\Delta} \mu}{H^2 \left(h_1 +\sqrt{{\epsilon}}\right)}, \quad r \le \min\left\{\frac{\min_{i \in [M]} C^{(i)}(K_0)}{\Bar{h}_{\text{cost}}}, \underline{h}_{\Delta} \right\}.$$
Following the analysis in Eq.\eqref{eq:model_free_prp}, we have
\begin{align*}
&C^{(i)}(K_{n+1}) - C^{(i)}(K_n)  \\
&\le -\left\langle\nabla C^{(i)}(K_n), \frac{\eta}{ML}\sum_{j=1}^M \sum_{l=0}^{L-1} \tilde{\nabla} C^{(j)}(K_{n, l}^{(j)})-\frac{\eta}{ML}\sum_{j=1}^M \sum_{l=0}^{L-1} \nabla C^{(j)}(K_{n, l}^{(j)})\right\rangle\notag\\
& -\left\langle\nabla C^{(i)}(K_n), \frac{\eta}{ML}\sum_{j=1}^M \sum_{l=0}^{L-1} \nabla C^{(j)}(K_{n, l}^{(j)}) -\nabla C^{(j)}(K_{n})\right\rangle\notag\\
& -\left\langle\nabla C^{(i)}(K_n), \frac{\eta}{M}\sum_{j=1}^M  \nabla C^{(j)}(K_{n})-\nabla C^{(i)}(K_n)\right\rangle - \eta\left\lVert C^{(i)}(K_n)\right\rVert^2\notag\\
&+\frac{h_{\text {grad}}(K_n)}{2} \Big\lVert \frac{\eta}{ML}\sum_{j=1}^M \sum_{l=0}^{L-1}\tilde{\nabla} C^{(j)}(K_{n, l}^{(j)})\Big\rVert_F^2 \notag\\
& \stackrel{(a)}{\le} \eta \Big\lVert \nabla C^{(i)}(K_n)\Big\rVert_F \Big\lVert\frac{1}{ML}\sum_{j=1}^M \sum_{l=0}^{L-1}\left[\tilde{\nabla} C^{(j)}(K_{n, l}^{(j)})- \nabla C^{(j)}(K_{n, l}^{(j)})\right]\Big\rVert_F\notag\\
& +\eta \Big\lVert\nabla C^{(i)}(K_n)\Big\rVert_F \Big\lVert\frac{1}{ML}\sum_{j=1}^M \sum_{l=0}^{L-1} \left[\nabla C^{(j)}(K_{n, l}^{(j)})-\nabla C^{(j)}(K_{n})\right]\Big\rVert_F \notag\\
& +\eta \Big\lVert\nabla C^{(i)}(K_n)\Big\rVert_F \Big\lVert\frac{1}{M}\sum_{j=1}^M \left[\nabla C^{(j)}(K_{n})-\nabla C^{(i)}(K_{n})\right]\Big\rVert_F-\eta \Big\lVert \nabla C^{(i)}(K_n)\Big\rVert_F^2\notag\\
&+\frac{h_{\text {grad}}(K_n)}{2} \Big\lVert \frac{\eta}{ML}\sum_{j=1}^M \sum_{l=0}^{L-1}\tilde{\nabla} C^{(j)}(K_{n, l}^{(j)})\Big\rVert_F^2 \notag\\
& \stackrel{(b)}{\le} \frac{\eta}{4} \Big\lVert \nabla C^{(i)}(K_n)\Big\rVert_F^2 + \eta\Big\lVert\frac{1}{ML}\sum_{j=1}^M \sum_{l=0}^{L-1}\left[\tilde{\nabla} C^{(j)}(K_{n, l}^{(j)})- \nabla C^{(j)}(K_{n, l}^{(j)})\right]\Big\rVert_F^2\notag\\
& +\frac{\eta}{8} \Big\lVert \nabla C^{(i)}(K_n)\Big\rVert_F^2 +  \frac{2\eta h_{\text {grad}}(K_n)^2}{ML}\sum_{j=1}^M \sum_{l=0}^{L-1}\Big\lVert K_{n, l}^{(j)}- K_{n}\Big\rVert_F^2\notag\\
& +\frac{\eta}{4} \Big\lVert \nabla C^{(i)}(K_n)\Big\rVert_F^2 +  \frac{\eta}{M}\sum_{j=1}^M \Big\lVert  \nabla C^{(j)}(K_{n})- C^{(i)}(K_{n})\Big\rVert_F^2\notag\\
&-\eta \Big\lVert \nabla C^{(i)}(K_n)\Big\rVert_F^2
+\frac{h_{\text {grad}}(K_n)}{2} \Big\lVert \frac{\eta}{ML}\sum_{j=1}^M \sum_{l=0}^{L-1}\tilde{\nabla} C^{(j)}(K_{n, l}^{(j)})\Big\rVert_F^2, 
\end{align*}
where $(a)$ is due to Cauchy–Schwarz inequality; and $(b)$ is due to Cauchy–Schwarz inequality and Eq.\eqref{eq:youngs_inner}.
Moreover, we have
\begin{align}\label{eq:main_lqr_analysis}
&C^{(i)}(K_{n+1}) - C^{(i)}(K_n) \notag\\
&\stackrel{(b)}{\le} \frac{\eta}{4} \Big\lVert \nabla C^{(i)}(K_n)\Big\rVert_F^2 + \eta\Big\lVert\frac{1}{ML}\sum_{j=1}^M \sum_{l=0}^{L-1}\left[\tilde{\nabla} C^{(j)}(K_{n, l}^{(j)})- \nabla C^{(j)}(K_{n, l}^{(j)})\right]\Big\rVert_F^2\notag\\
& +\frac{\eta}{8} \Big\lVert \nabla C^{(i)}(K_n)\Big\rVert_F^2 +  \frac{2\eta h_{\text {grad}}(K_n)^2}{ML}\sum_{j=1}^M \sum_{l=0}^{L-1}\Big\lVert K_{n, l}^{(j)}- K_{n}\Big\rVert_F^2\notag\\
& +\frac{\eta}{4} \Big\lVert \nabla C^{(i)}(K_n)\Big\rVert_F^2 +  \frac{\eta}{M}\sum_{j=1}^M \Big\lVert  \nabla C^{(j)}(K_{n})- C^{(i)}(K_{n})\Big\rVert_F^2\notag\\
&-\eta \Big\lVert \nabla C^{(i)}(K_n)\Big\rVert_F^2
+\frac{h_{\text {grad}}(K_n)}{2} \Big\lVert \frac{\eta}{ML}\sum_{j=1}^M \sum_{l=0}^{L-1}\tilde{\nabla} C^{(j)}(K_{n, l}^{(j)})\Big\rVert_F^2 \notag\\
&\stackrel{(c)}{\le} -\frac{3\eta}{8}\Big\lVert \nabla C^{(i)}(K_n)\Big\rVert_F^2 + \eta \epsilon + \frac{4\eta^2 \bar{h}_{\text {grad}}^2}{\eta_g M}\sum_{j=1}^M \left[\Big\lVert\nabla C^{(j)}(K_{n})\Big\rVert_F^2 +ML\epsilon \right]\notag\\
&+ \eta \min\{n_x,n_u\}(\epsilon_1 \bar{h}^1_{\text{het}} + \epsilon_2 \bar{h}^2_{\text{het}})^2+\frac{h_{\text {grad}}(K_n)}{2} \Big\lVert \frac{\eta}{ML}\sum_{j=1}^M \sum_{l=0}^{L-1}\tilde{\nabla} C^{(j)}(K_{n, l}^{(j)})\Big\rVert_F^2,
\end{align}
where $(b)$ follows from the gradient Lipschitz property in Lemma~\ref{Lemma:lipschitz}; and  $(c)$ follows from the policy gradient heterogeneity property in Lemma~\ref{lem:gradient_het}, Lemma~\ref{lem:variance_reduction} and Lemma~\ref{lem：drift_term_model_free}.

Following the analysis in Eq.\eqref{eq:main_lqr_analysis}, we have 
\begin{align}\label{eq:main_lqr_middle}
&C^{(i)}(K_{n+1}) - C^{(i)}(K_n) \notag\\
&\stackrel{(d)}{\le} -\frac{3\eta}{8}\Big\lVert \nabla C^{(i)}(K_n)\Big\rVert_F^2 + \eta \epsilon + \frac{4\eta^2 \bar{h}_{\text {grad}}^2}{\eta_g M}\sum_{j=1}^M \left[\Big\lVert\nabla C^{(j)}(K_{n})\Big\rVert_F^2 +ML\epsilon \right]\notag\\
&+ \eta \min\{n_x,n_u\}(\epsilon_1 h^1_{\text{het}} + \epsilon_1 h^2_{\text{het}})^2\notag\\
&+\frac{4 \eta^2 \bar{h}_{\text {grad}}}{2} \Big\lVert \frac{1}{ML}\sum_{j=1}^M \sum_{l=0}^{L-1}\tilde{\nabla} C^{(j)}(K_{n, l}^{(j)})-\nabla C^{(j)}(K_{n, l}^{(j)})\Big\rVert_F^2 \notag\\
&+\frac{4 \eta^2 \bar{h}_{\text {grad}}}{2} \Big\lVert \frac{1}{ML}\sum_{j=1}^M \sum_{l=0}^{L-1}\nabla C^{(j)}(K_{n, l}^{(j)})-\nabla C^{(j)}(K_{n})\Big\rVert_F^2 \notag\\
&+\frac{4 \eta^2 \bar{h}_{\text {grad}}}{2M} \sum_{j=1}^M\Big\lVert  \nabla C^{(j)}(K_{n})-\nabla C^{(i)}(K_{n})\Big\rVert_F^2 +\frac{4 \eta^2 \bar{h}_{\text {grad}}}{2}\Big\lVert\nabla C^{(i)}(K_{n})\Big\rVert_F^2\notag\\
&\stackrel{(e)}{\le} -\frac{3\eta}{8}\Big\lVert \nabla C^{(i)}(K_n)\Big\rVert_F^2 + (\eta +2\eta^2 \bar{h}_{\text {grad}} ) \epsilon + \frac{4\eta^2 \bar{h}_{\text {grad}}^2}{\eta_g M}\sum_{j=1}^M \left[\Big\lVert\nabla C^{(j)}(K_{n})\Big\rVert_F^2 +ML\epsilon \right]\notag\\
&+ (\eta +2\eta^2 \bar{h}_{\text {grad}} ) \min\{n_x,n_u\}(\epsilon_1 h^1_{\text{het}} + \epsilon_2 h^2_{\text{het}})^2 +2\eta^2 \bar{h}_{\text {grad}} \Big\lVert \nabla C^{(i)}(K_n)\Big\rVert_F^2 \notag\\
&+\frac{4 \eta^2 \bar{h}_{\text {grad}}}{2} \Big\lVert \frac{1}{ML}\sum_{j=1}^M \sum_{l=0}^{L-1}\nabla C^{(j)}(K_{n, l}^{(j)})-\nabla C^{(j)}(K_{n})\Big\rVert_F^2\notag\\
&\stackrel{(f)}{\le} -\left(\frac{3\eta}{8}+ 2\eta^2 \bar{h}_{\text {grad}}\right)\Big\lVert \nabla C^{(i)}(K_n)\Big\rVert_F^2 + (\eta +2\eta^2 \bar{h}_{\text {grad}} ) \epsilon \notag\\
&+ \frac{4\eta^2 \bar{h}_{\text {grad}}^2}{\eta_g M}\sum_{j=1}^M \left[\Big\lVert\nabla C^{(j)}(K_{n})\Big\rVert_F^2 +ML\epsilon \right]
+ (\eta +2\eta^2 \bar{h}_{\text {grad}} ) \min\{n_x,n_u\}(\epsilon_1 h^1_{\text{het}} + \epsilon_2 h^2_{\text{het}})^2 \notag\\
&+\frac{4 \eta^2 \bar{h}^2_{\text {grad}}}{2ML} \sum_{j=1}^M \sum_{l=0}^{L-1}\Big\lVert K_{n, l}^{(j)}-K_{n}\Big\rVert_F^2\notag\\
&\stackrel{(g)}{\le} -\left(\frac{3\eta}{8}+ 2\eta^2 \bar{h}_{\text {grad}}\right)\Big\lVert \nabla C^{(i)}(K_n)\Big\rVert_F^2 + (\eta +2\eta^2 \bar{h}_{\text {grad}} ) \epsilon \notag\\
&+ \frac{4\eta^2 \bar{h}_{\text {grad}}^2+4\eta^3 \bar{h}_{\text {grad}}^2}{\eta_g M}\sum_{j=1}^M \left[\Big\lVert\nabla C^{(j)}(K_{n})\Big\rVert_F^2 +ML\epsilon \right]\notag\\
&+ (\eta +2\eta^2 \bar{h}_{\text {grad}} ) \min\{n_x,n_u\}(\epsilon_1 \bar{h}^1_{\text{het}} + \epsilon_2 \bar{h}^2_{\text{het}})^2,
\end{align}
where $(d)$ is due to Eq.\eqref{eq:sum_expand}; $(e)$ is due to variance reduction property in Lemma~\ref{lem:variance_reduction} and policy gradient heterogeneity in Lemma~\ref{lem:gradient_het}; $(f)$ is due to gradient Lipschitz property in Lemma~\ref{Lemma:lipschitz}; $(g)$ is due to drift term analysis in Lemma~\ref{lem：drift_term_model_free}.

Continuing the analysis in Eq.\eqref{eq:main_lqr_middle}, we have that
\begin{align}\label{eq:final_result}
&C^{(i)}(K_{n+1}) - C^{(i)}(K_n)\le-\left(\frac{3\eta}{8}+ 2\eta^2 \bar{h}_{\text {grad}}\right)\Big\lVert \nabla C^{(i)}(K_n)\Big\rVert_F^2 + (\eta +2\eta^2 \bar{h}_{\text {grad}} ) \epsilon \notag\\
&+ \frac{4\eta^2 \bar{h}_{\text {grad}}^2+4\eta^3 \bar{h}_{\text {grad}}^2}{\eta_g M}\sum_{j=1}^M \left[\Big\lVert\nabla C^{(j)}(K_{n})\Big\rVert_F^2 +ML\epsilon \right]\notag\\
&+ (\eta +2\eta^2 \bar{h}_{\text {grad}} ) \min\{n_x,n_u\}(\epsilon_1 \bar{h}^1_{\text{het}} + \epsilon_1 \bar{h}^2_{\text{het}})^2\notag\\
& \stackrel{(a)}{\le} -\left(\frac{3\eta}{8}+ 2\eta^2 \bar{h}_{\text {grad}}\right)\Big\lVert \nabla C^{(i)}(K_n)\Big\rVert_F^2 + (\eta +2\eta^2 \bar{h}_{\text {grad}} ) \epsilon \notag\\
&+ \frac{4\eta^2 \bar{h}_{\text {grad}}^2+4\eta^3 \bar{h}_{\text {grad}}^2}{\eta_g M}\sum_{j=1}^M \left[2\Big\lVert\nabla C^{(j)}(K_{n})-\nabla C^{(i)}(K_{n})\Big\rVert_F^2+ 2\Big\lVert\nabla C^{(i)}(K_{n})\Big\rVert_F^2 +ML\epsilon \right]\notag\\
&+ (\eta +2\eta^2 \bar{h}_{\text {grad}} ) \min\{n_x,n_u\}(\epsilon_1 \bar{h}^1_{\text{het}} + \epsilon_1 \bar{h}^2_{\text{het}})^2\notag\\
& \stackrel{(b)}{\le} -\left(\frac{3\eta}{8}+ 2\eta^2 \bar{h}_{\text {grad}}+\frac{8\eta^2 \bar{h}_{\text {grad}}^2+8\eta^3 \bar{h}_{\text {grad}}^2}{\eta_g }\right)\Big\lVert \nabla C^{(i)}(K_n)\Big\rVert_F^2\notag\\
&+ \left(\eta +2\eta^2 \bar{h}_{\text {grad}}+ \frac{4\eta^2 \bar{h}_{\text {grad}}^2+4\eta^3 \bar{h}_{\text {grad}}^2}{\eta_g }ML\right) \epsilon \notag\\
&+\left(\eta +2\eta^2 \bar{h}_{\text {grad}} +\frac{8\eta^2 \bar{h}_{\text {grad}}^2+8\eta^3 \bar{h}_{\text {grad}}^2}{\eta_g }\right) \min\{n_x,n_u\}(\epsilon_1 \bar{h}^1_{\text{het}} + \epsilon_2 \bar{h}^2_{\text{het}})^2\notag\\
& \stackrel{(c)}{\le} -\frac{\eta}{4}\Big\lVert \nabla C^{(i)}(K_n)\Big\rVert_F^2\notag+ 2\eta  \epsilon+2\eta \min\{n_x,n_u\}(\epsilon_1 \bar{h}^1_{\text{het}} + \epsilon_2 \bar{h}^2_{\text{het}})^2\notag\\
& \stackrel{(d)}{\le} -\frac{\eta \mu^2 \sigma_{\min }(R)}{\left\|\Sigma_{K_i^*}\right\|} (C^{(i)}(K_n) - C^{(i)}(K_i^*))+ 2\eta  \epsilon+2\eta\min\{n_x,n_u\}(\epsilon_1 \bar{h}^1_{\text{het}} + \epsilon_2 \bar{h}^2_{\text{het}})^2,
\end{align}
where $(a)$ is due to Eq.\eqref{eq:sum_expand};  $(b)$ is due to policy gradient heterogeneity in Lemma~\ref{lem:gradient_het}; and $(c)$ is due to the choice of step-size such that $\frac{3\eta}{8}+ 2\eta^2 \bar{h}_{\text {grad}}+\frac{8\eta^2 \bar{h}_{\text {grad}}^2+8\eta^3 \bar{h}_{\text {grad}}^2}{\eta_g }\le \frac{\eta}{4}$ and $$\eta +2\eta^2 \bar{h}_{\text {grad}} +\frac{8\eta^2 \bar{h}_{\text {grad}}^2+8\eta^3 \bar{h}_{\text {grad}}^2}{\eta_g}\le 2\eta,$$ which holds when $\eta \le \min \{\frac{1}{32\bar{h}_{\text {grad}}},1\}$ and $\eta_l \le \frac{1}{256L\bar{h}^2_{\text {grad}}};$ for $(d)$ we use the gradient domination lemma in Lemma~\ref{lemma:gradient_domination}.

In conclusion, we have that
\begin{align*}
C^{(i)}(K_{n+1}) - C^{(i)}(K_i^*)&\le \left(1-\frac{\eta \mu^2 \sigma_{\min }(R)}{\left\|\Sigma_{K_i^*}\right\|}\right) (C^{(i)}(K_n) - C^{(i)}(K_i^*)) + 2\eta \epsilon\notag\\
&+ 2\eta\min\{n_x,n_u\}(\epsilon_1 \bar{h}^1_{\text{het}} + \epsilon_2 \bar{h}^2_{\text{het}})^2,
\end{align*} 
holds when the step-size, smoothing radius, trajectory length, and sample size satisfy the requirements mentioned above and those in Lemma~\ref{lem:stability_local_model_free} and Lemma~\ref{lem：drift_term_model_free}.
\end{proof}

With this lemma, we are now ready to provide the convergence guarantees for the \texttt{FedLQR} under the model-free setting. 
\paragraph{Proof of the iterative stability guarantees of \texttt{FedLQR}:}  Here, we leverage the method of induction to prove 
\texttt{FedLQR}'s iterative stability guarantees. First, we start from an initial policy $K_0 \in \mathcal{G}^0$. At round $n$, we assume $K_n \in \mathcal{G}^0.$ According to Lemma~\ref{lem:stability_local_model_free}, we have that all the local policies $K_{n,l}^{(i)} \in \mathcal{G}^0$. Furthermore, frame the hypotheses of in Lemma~\ref{lem:model_free_prp}, we have that 
\begin{align*}
C^{(i)}(K_{n+1}) - C^{(i)}(K_i^*)&\le \left(1-\frac{\eta \mu^2 \sigma_{\min }(R)}{\left\|\Sigma_{K_i^*}\right\|}\right) (C^{(i)}(K_n) - C^{(i)}(K_i^*)) + 2\eta \epsilon\notag\\
&+ 2\eta\min\{n_x,n_u\}(\epsilon_1 \bar{h}^1_{\text{het}} + \epsilon_2 \bar{h}^2_{\text{het}})^2.
\end{align*}

Since $(\epsilon_1 \Bar{h}^1_{\text{het}} + \epsilon_2 \Bar{h}^2_{\text{het}})^2 \le\Bar{h}^3_{\text{het}},$ we have
\begin{equation}
\begin{aligned}
C^{(i)}(K_{n+1}) &- C^{(i)}(K_i^*) 
\le\left(1- \frac{\eta \mu^2 \sigma_{\min }(R)}{\left\|\Sigma_{K_i^*}\right\|}\right) (C^{(i)}(K_{0}) - C^{(i)}(K_i^*)) +2\eta \epsilon \notag\\
&+\frac{\eta \mu^2 \sigma_{\min }(R)}{2\left\|\Sigma_{K_i^*}\right\|}(C^{(i)}(K_{0}) - C^{(i)}(K_i^*))\notag\\
& \stackrel{(a)}{\le} C^{(i)}(K_{0}) - C^{(i)}(K_i^*),
\end{aligned}
\end{equation}
where $(a)$ follows from the fact that $\epsilon$ can be arbitrarily small by choosing a small smoothing radius, sufficient long trajectory length, and enough samples. 

With this, we can easily have that the global policy $K_{n+1}$ at the next round $n+1$ is also stabilizing, i.e., $K_{n+1} \in \mathcal{G}^0.$ Therefore, we can finish proving \texttt{FedLQR}'s iterative stability property by inductively reasoning.

\paragraph{Proof of \texttt{FedLQR}'s convergence:} From Eq.\eqref{eq:one_step_convergence}, we have 
\begin{align*}
C^{(i)}(K_{n+1}) - C^{(i)}(K_i^*)&\le \left(1-\frac{\eta \mu^2 \sigma_{\min }(R)}{\left\|\Sigma_{K_i^*}\right\|}\right) (C^{(i)}(K_n) - C^{(i)}(K_i^*)) + 2\eta \epsilon + 2\eta\min\{n_x,n_u\}(\epsilon_1 \bar{h}^1_{\text{het}} + \epsilon_2\bar{h}^2_{\text{het}})^2,
\end{align*} 

Using the above inequality recursively,  \texttt{FedLQR} enjoys the following convergence guarantee after $N$ rounds:
\begin{align*}
C^{(i)}(K_{N}) - C^{(i)}(K_i^*)&\le \left(1-\frac{\eta \mu^2 \sigma_{\min }(R)}{\left\|\Sigma_{K_i^*}\right\|}\right)^N (C^{(i)}(K_0) - C^{(i)}(K_i^*)) + \frac{2\left\|\Sigma_{K_i^*}\right\|}{\mu^2 \sigma_{\min }(R)} \epsilon\notag\\
&+ \frac{2\min\{n_x,n_u\}\left\|\Sigma_{K_i^*}\right\|}{\mu^2 \sigma_{\min }(R)}(\epsilon_1 \bar{h}^1_{\text{het}} + \epsilon_2 \bar{h}^2_{\text{het}})^2.
\end{align*}
Suppose the trajectory length satisfies $
\tau \geq h_{\tau}\left(\frac{r \epsilon^{\prime}}{4 n_x n_u}\right),$ the smoothing radius satisfies $r\le  h^{\prime}_r\left(\frac{\epsilon^{\prime}}{4}\right),$ where $$h^{\prime}_r\left(\frac{\epsilon^{\prime}}{4}\right):= \min\left\{\frac{\min_{i \in [M]} C^{(i)}(K_0)}{\Bar{h}_{\text{cost}}}, \underline{h}_{\Delta}, h_r\left(\frac{\epsilon^{\prime}}{4}\right)\right\},$$ 
and the number of the sample size of each agent $n_s$ satisfies $$n_{s} \geq \frac{h_{\text {sample,trunc }}\left(\frac{\epsilon^{\prime}}{4}, \frac{\delta}{ML}, \frac{H^2}{\mu}\right)}{ML},$$ with $  \epsilon^{\prime} = \frac{4 \left\| \Sigma_{K_i^\star}  \right\|}{\mu^2 \sigma_{\min}(R)} \epsilon.$

When the number of rounds $N \ge \frac{c_{\text{uni},4}\left\|\Sigma_{K_i^*}\right\|}{\eta\mu^2 \sigma_{\min}(R)}\log\left(\frac{2(C^{(i)}(K_0) - C^{(i)}(K_i^*))}{\epsilon^{\prime}}\right)$,  our \texttt{FedLQR} algorithm enjoys the following convergence guarantee:
\begin{align*}
C^{(i)}(K_{N}) - C^{(i)}(K_i^*)  &\le \left(1 -\frac{\eta \mu^2 \sigma_{\min }(R)}{\left\|\Sigma_{K_i^*}\right\|} \right)^N(C^{(i)}(K_0) - C^{(i)}(K_i^*)) +\frac{\epsilon^{\prime}}{2}\\
&+ \frac{2\min\{n_x, n_u\} \left\|\Sigma_{K_i^*}\right\|}{\mu^2 \sigma_{\min}(R)} (\epsilon_1 h^1_{\text{het}} + \epsilon_2 h^2_{\text{het}})^2 \notag\\
& \le \epsilon^{\prime}+\frac{2\min\{n_x, n_u\} \left\|\Sigma_{K_i^*}\right\|}{\mu^2 \sigma_{\min}(R)} (\epsilon_1 h^1_{\text{het}} + \epsilon_2 h^2_{\text{het}})^2,
\end{align*}
which completes the proof.

\end{document}